\documentclass[twoside,11pt]{article}

\usepackage{jmlr2e,natbib}

\usepackage{xr-hyper}
\makeatletter
\newcommand*{\addFileDependency}[1]{
  \typeout{(#1)}
  \@addtofilelist{#1}
  \IfFileExists{#1}{}{\typeout{No file #1.}}
}
\makeatother

\usepackage{amsfonts, amsmath, graphicx, float,
  subcaption, enumerate, color, verbatim ,algorithm,
  algorithmic, url, hyperref, epsfig, amssymb, bm, natbib}
\usepackage[shortlabels]{enumitem}
\graphicspath{{figures/}}
\allowdisplaybreaks

\newcommand{\0}{\bm{0}}
\newcommand{\1}{\bm{1}}
\renewcommand{\H}{\mathbf{H}}

\newcommand{\F}{\mathbf{F}}

\renewcommand{\l}{\bm{l}}
\newcommand{\M}{\mathbf{M}}
\newcommand{\p}{\bm{p}}

\renewcommand{\u}{\bm{U}}

\newcommand{\V}{\mathbf{V}}

\newcommand{\x}{\bm{x}}

\newcommand{\y}{\bm{y}}

\newcommand{\z}{\bm{z}}

\newcommand{\btheta}{\bm{\theta}}
\newcommand{\bbeta}{\bm{\btheta}} %

\newcommand{\hbeta}{\hat{\bm{\btheta}}} %
\newcommand{\tbeta}{\tilde{\bm{\btheta}}} %

\newcommand{\bSigma}{\bm{\Sigma}}

\newcommand{\eeta}{\bm{\eta}}%
\newcommand{\tPsi}{\widetilde{\Psi}}

\newcommand{\ttheta}{{\tilde{\btheta}}}
\newcommand{\htheta}{{\hat{\btheta}}}
\newcommand{\bvtheta}{{\bm{\vartheta}}}

\newcommand{\tmu}{{\tilde{\mu}}}

\newcommand{\onen}{\frac{1}{n}}
\newcommand{\oneN}{\frac{1}{N}}

\newcommand{\op}{o_{P}(1)}
\newcommand{\Op}{O_{P}(1)}
\newcommand{\ud}{\mathrm{d}}

\newcommand{\sumn}{\sum_{i=1}^{n}}
\newcommand{\sumK}{\sum_{k=1}^{K}}

\newcommand{\sumN}{\sum_{i=1}^{N}}

\newcommand{\tp}{^{\mathrm{T}}}
\newcommand{\tr}{\mathrm{tr}}
\newcommand{\cvp}{\xrightarrow{P}}
\newcommand{\cvd}{\xrightarrow{D}}
\renewcommand{\Pr}{\mathbb{P}}
\newcommand{\plt}{\mathrm{plt}}
\newcommand{\Exp}{\mathbb{E}}
\newcommand{\Nor}{\mathbb{N}}
\newcommand{\Var}{\mathbb{V}}

\newcommand{\lus}{{\textnormal{\tiny LUS}}}%
\newcommand{\mle}{{\textnormal{\tiny MLE}}}%

\newtheorem{assumption}{\bf Assumption}
\newenvironment{assumptionp}[1]{%
  \renewcommand\theassumption{#1}
  \assumption
}{\endtheorem}

\usepackage{lastpage}
\jmlrheading{23}{2022}{1-\pageref{LastPage}}{5/21; Revised
9/22}{10/22}{21-0506}{HaiYing Wang and Jae Kwang Kim}
\ShortHeadings{Maximum sampled conditional likelihood for informative subsampling}{Wang and Kim}

\begin{document}

\title{Maximum sampled conditional likelihood for informative subsampling}

\author{\name HaiYing Wang \email haiying.wang@uconn.edu \\
       \addr Department of Statistics\\
       University of Connecticut\\
       Storrs, CT 06269, USA
       \AND
       \name Jae Kwang Kim \email jkim@iastate.edu \\
       \addr Department of Statistics\\
       Iowa State University\\
       Ames, IA 50011, USA
       }

\editor{Tong Zhang}

\maketitle

\begin{abstract}
Subsampling is a computationally effective approach to extract information from massive data sets when computing resources are limited. After a subsample is taken from the full data, most available methods use an inverse probability weighted (IPW) objective function to estimate the model parameters. The IPW estimator does not fully utilize the information in the selected subsample. In this paper, we propose to use the maximum sampled conditional likelihood estimator (MSCLE) based on the sampled data. We established the asymptotic normality of the MSCLE and prove that its asymptotic variance covariance matrix is the smallest among a class of asymptotically unbiased estimators, including the IPW estimator. We further discuss the asymptotic results with the L-optimal subsampling probabilities and illustrate the estimation procedure with generalized linear models. Numerical experiments are provided to evaluate the practical performance of the proposed method.
\end{abstract}


\section{Introduction}
\label{sec:introduction}
In the era of big data, many data sets have  huge volumes. If the data sets  are too voluminous  then  traditional data processing software products are not capable of  processing the data within a reasonable amount of time.
In this case, a subsample or coreset of the full data is often used to alleviate
the computational burden. Subsampling is an emerging area of research that
balances the trade-off between computational efficiency and statistical
efficiency by developing an efficient subsampling design and estimation
strategy.
For this purpose, existing research  focuses more on designing the subsampling probabilities and less on improving the estimator based on the selected subsample, e.g., \cite{drineas2006sampling,yang2015explicit,WangZhuMa2018}, among others. 

In the linear regression model setup, optimal subsampling designs are well studied in the literature.  Specifically,  statistical leverage scores or their variants are often recommended to construct subsampling probabilities, see
\cite{drineas2006sampling,dhillon2013new,mcwilliams2014fast,PingMa2014-JMLR,yang2015explicit,nie2018minimax}, and the references therein. {Instead of calculating exact leverage scores directly on the  %
full data,} \cite{Drineas:12} proposed fast algorithms to approximate them. The
aforementioned sampling probabilities are  not dependent on the response
variable, and this type of sampling schemes is referred to as non-informative
subsampling. For this scenario, \cite{WangYangStufken2019} proposed a
deterministic selection algorithm that has high estimation efficiency.

Beyond linear regressions, \cite{WangZhuMa2018} proposed an optimal
subsampling method under the A-optimality criterion for logistic regression,
which defines subsampling probabilities that minimize the asymptotic mean
squared error of the resulting subsample estimator. They further considered the
L-optimality to further improve the computational efficiency. This method has
been extended to include other models such as generalized linear models
\citep{ai2020optimal}, quantile regressions \citep{wang2021optimal},
quasi-likelihood models \citep{yu2020quasi}, and etc. \cite{influence2018}
suggested using the influence function to define optimal
probabilities. \cite{shen2021surprise} proposed the surprise sampling method
that gives optimal forms to a variety of objectives. \cite{WangZou2019}
systematically compared with-replacement sampling with Poisson sampling and
recommended Poisson sampling for its higher estimation efficiency and
computational feasibility.  Readers are refered to \cite{YaoWang2021JDS} for a
systematic review on this topic.

Unlike non-informative subsampling, the optimal subsampling probabilities depend on the response variable as well as the covariates.
If the sampling probabilities depend on the response variable in addition to the
auxiliary (covariate) variable, the sampling mechanism is called informative
\citep{pfe98}. The selection probabilities in the 
informative sampling  utilize information in both the covariates and the responses, so the resulting subsample often contain more relevant information compared with non-informative sampling. Under informative sampling, the selection probabilities are often inversely applied to obtain the IPW estimator \citep{cha03}. 
However, the inverse probability weighting scheme %
may not achieve efficient estimation. %

Generally speaking, the informative subsampling can be viewed as a biased sampling problem in statistics, as discussed in \cite{cox69} and \cite{qin2017biased}. 
To understand the biased sampling problem, it is useful to consider selection bias in the context of two-phases of sampling. In the first phase, we have a random sample of size $N$.  In the second phase, we select a subset of the original sample with known selection probability $\pi(\x,y)$ which is a function of the observations in the first-phase sample.  If the selection probability depends on the outcome variable $y$, it is also called outcome-dependent (two-phase)   sampling. The two-phase sampling design is commonly used in many disciplines. \cite{kim2006} and \cite{saegusa2013} developed some theory for two-phase sampling. 
Some examples of two-phase sampling can be found in \cite{hsi85}, \cite{kal88}, \cite{wil91}, \cite{sco91}, \cite{hu96}, \cite{sco97}, \cite{hu97}, \cite{bre97}, and \cite{whi97}. The case-control study is a popular example of the outcome-dependent two-phase sampling. If the outcome is binary and the case with $y=1$ is rare, it is sensible to oversample  the cases with $y=1$ in the final sample. 
Such outcome-dependent two-phase sampling is used
because it is either more efficient or cost effective.
Removing or reducing the selection bias in such sampling is a crucial part of the estimation
problem.

For logistic regressions with case-control sampling, \cite{scott1986fitting}
showed that the bias for the unweighted estimator only appears in the intercept
estimator and they provided the expression of the bias term. For
extremely imbalanced data, %
\cite{Wang2020RareICML} proved that with a case-control subsample, the IPW
estimator and the unweighted estimator with bias correction have the
same convergence rate and asymptotic distribution as the full data estimator if
sufficient number of controls are selected; otherwise the latter is more efficient than the former. The
case-control sampling probabilities depend on the responses only and do not
utilize the information in the covariates. For binary logistic regression,
\cite{fithian2014local} proposed the local case-control subsampling
probabilities which depend on both the responses and the covariates. They showed
that the bias in the regression coefficient estimator is corrected by the pilot
estimator used to calculated the sampling probabilities. \cite{wang2019more}
extended this bias correction idea to the optimal subsampling probabilities
under the A- and L-optimality, and proved that an unweighted estimator with bias
correction has a higher estimation efficiency. For multi-class logistic
regression, \cite{Han2019} developed the local uncertainty sampling (LUS) that
focuses on correcting the unweighted objective function instead of correcting
the bias in the resulting estimator. This approach is using the sampled data
conditional likelihood and is not restricted by a specific form of the
subsampling probabilities. \cite{WangZhangWang2021} adopted the idea and
investigated optimal negative subsampling for the case of extremely imbalanced
binary data. 

The aforementioned investigations exclusively focus on logistic regression. In this paper, we will show that there is a general approach to extract more information from an informative subsample without using inverse probability weighting. The basic strategy is to treat the subsample estimation as a missing data problem and obtain the conditional likelihood of the sampled data, which is based on the conditional density function of the study variable given the covariate variable for the sampled data. The sampled conditional likelihood has been used in the context of biased sampling problem, but to our best knowledge,  it has not been addressed in the subsampling area.
The sampled conditional likelihood approach is applicable for general parametric models and sampling probabilities. The investigations of \cite{fithian2014local, wang2019more,Han2019} are all specific cases of this general approach. We first establish the consistency and the asymptotic normality of the maximum sampled conditional likelihood estimator under some regularity conditions. Thus, statistical inference such as normal-based confidence intervals can be developed. We also show that the resulting estimator has the highest estimation efficiency {among a class of asymptotically unbiased estimators}, and it is more efficient than the IPW estimator. The  maximum sampled conditional likelihood estimator can be computed by applying the Fisher-scoring algorithm with the closed-form formula for the Hessian matrix. Thus, the computation is relatively simple and fast. As illustrated in the simulation study in Section 6, the efficiency gains of using the sampled conditional likelihood over the IPW estimator are  substantial.

If the subsampling probabilities are unknown, we may use an independent pilot sample to estimate the parameters in the subsampling probabilities. Because the subsampling probabilities are under our control, unlike the missing data problem, we can always use the correct subsampling probabilities in constructing the sampled conditional likelihood. Thus, even if the pilot samples are systematically different from the original sample, the statistical properties of the maximum sampled conditional likelihood estimator remain valid. 

The rest of the paper is organized as follows. We present the proposed idea in
Section~\ref{sec:sampl-data-likel} and discuss its asymptotic properties in
Section~\ref{sec:asymptotic-results}. Section~\ref{sec:estim-param-subs}
considers the practical situation that informative subsampling probabilities
depend on unknowns. Section~\ref{sec:subs-prob-based} illustrates structural
results using a version of widely used informative subsampling probabilities and
Section~\ref{sec:gener-line-models} demonstrates the proposed method in the
context of generalized linear models. Section~\ref{sec:numer-exper} provides
numerical results based on simulated and real data sets. Some concluding remarks
are made in Section 8. Proofs, technical details, and additional numerical
results are given in the appendix.

\section{Sampled data conditional likelihood estimator}
\label{sec:sampl-data-likel}
Let $(\x_i,y_i)$, $i=1, ..., N$, be an independent sample from the distribution of $(\x,y)$, %
where $\x$ is the covariate variable and $y$ is the main response variable. Denote the density function of $y_i$ given $\x_i$ as $f(y_i\mid\x_i,\btheta)$, where $\btheta$ is the parameter of interest, and it  is often estimated by the maximum likelihood estimator (MLE),
\begin{equation}
  \htheta_\mle
  =\arg\max_{\btheta}\sumN\ell(\btheta;\x_i,y_i),
\end{equation}
where $\ell(\btheta;\x_i,y_i)=\log f(y_i\mid\x_i,\btheta)$ is the log-likelihood function. 
For massive data, the computational cost in the maximization for $\htheta_\mle$ can be high, especially when there is no closed-form solution and an iterative algorithm has to be used. To solve this computational issue, approximation using a small subsample or a coreset of the data 
is regarded as an effective solution.

To emphasize the fact that informative subsampling probabilities depend on the responses, we denote them as $\pi(\x_i,y_i)= \Pr(\delta_i=1\mid\x_i,y_i)\in(0,1]$, $i=1, ..., N$. Here, $\delta_i$ is the indicator variable signifying if the $i$-th data point is selected, i.e., $\delta_i=1$ if $(\x_i,y_i)$ is in the subsample and $\delta_i=0$ otherwise. We assume that the distribution of $\delta_i$ is Bernoulli with parameter
\begin{equation}
  \Pr(\delta_i=1\mid\x_i,y_i)=\pi(\x_i,y_i),
  \quad\text{for}\quad i=1, ..., N.
\end{equation}

A commonly used subsample estimator is based on a IPW objective function, 
\begin{equation}\label{eq:4}
\htheta_{W}=\arg\max_{\btheta}\ell_{W}(\btheta)
=\arg\max_{\btheta}\sumN\frac{\delta_i\ell(\btheta;\x_i,y_i)}{\pi(\x_i,y_i)}.
\end{equation}
Here, the inverse probability weighting is necessary for $\htheta_{W}$ to be consistent. %
An unweighted estimator is biased and inconsistent to neither the full data MLE nor the true parameter. %
However, the weighting scheme in~(\ref{eq:4}) does not fully  extract the information in the subsample as it down-weights more informative data points.
Intuitively, one wants to assign a larger $\pi(\x_i,y_i)$ to a data point if it contains more information about $\btheta$ so that we sample it with a higher probability. 
If $\pi_i$'s are non-informative such as in the leverage-based sampling, then an unweighted estimator can still be consistent to the true parameter
\citep[e.g.,][]{PingMa2014-JMLR}. Actually, for non-informative sampling, an unweighted estimator is the best estimator with the smallest asymptotic variance.
If the subsampling mechanism is non-informative in the sense that the subsampling probability satisfies $\pi( \x, y) = \pi( \x)$, according to \cite{rubin1976}, the MLE of $\btheta$ is obtained by maximizing the complete-case (CC) log-likelihood
$$\htheta_{CC} 
=\arg\max_{\btheta}\sumN \delta_i\ell(\btheta;\x_i,y_i) . $$
Thus, the IPW estimator in (\ref{eq:4}) is inefficient.  
However, non-informative subsampling probabilities may not be as effective as informative subsampling probabilities in identifying informative data points in the first place. Thus we focus on informative subsampling only. 

To avoid the inverse probability weighting, we propose to use the sampled data conditional likelihood to obtain the subsample estimator. By Bayes' theorem the conditional density function of $y_i$ given $\x_i$ for sampled data is
\begin{equation}\label{eq:6}
  f(y_i\mid\x_i,\delta_i=1;\btheta)=
  \frac{f(y_i\mid\x_i;\btheta)\pi(\x_i,y_i)}
  {\int f(y\mid\x_i;\btheta)\pi(\x_i,y)\ud y},
\end{equation}
where $\ud y$ is the Lebesgue measure for continuous responses and it is the counting measure for discrete responses.

The sampled data conditional log-likelihood function from~(\ref{eq:6}) has a general form of
\begin{equation}\label{eq:5}
  \ell_{S}(\btheta)=\sumN\ell_{S}(\btheta;\x_i,y_i)
  =\sumN\delta_i\big[\log f(y_i\mid\x_i;\btheta)
    -\log\{\bar{\pi}(\x_i; \btheta)\}\big]+C,
\end{equation}
where 
\begin{equation*}
   \bar{\pi}(\x_i;\btheta)=
  \Exp\{\pi(\x_i,y_i)\mid\x_i\}
  =\int f(y\mid\x_i;\btheta)\pi(\x_i,y)\ud y,
\end{equation*}
and $C=\sumN\delta_i\log\{\pi(\x_i,y_i)\}$ does not contain $\btheta$.
The proposed estimator $\htheta_{S}$ is the maximizer of (\ref{eq:5}), namely,
\begin{equation}
  \htheta_{S}=\arg\max_{\btheta}\ell_{S}(\btheta).
  \label{msle} 
\end{equation}
{Note that the density in (\ref{eq:6}) is not the joint density of $(\x_i,y_i)$ for sampled data; it is the conditional density of $y_i$ given $\x_i$ for sampled data. Therefore we call the corresponding likelihood function the conditional likelihood, and 
call our estimator the maximum sampled conditional likelihood estimator (MSCLE).} If the subsampling probability $\pi(\x_i,y_i)$ depends on $\x_i$ only, then $\bar{\pi}(\x_i;\btheta)=\pi(\x_i)$ does not contain $\btheta$ and thus the MSCLE reduces to the MLE that maximizes the complete-case likelihood function.

Note that $\htheta_{W}$ uses $\pi^{-1}(\x_i,y_i)$'s as weights while
$\htheta_{S}$ uses $\log\{\bar{\pi}(\x_i; \btheta)\}$ to correct the
log-likelihood function. 
In computing $\htheta_{W}$, if a data point with a very
small value of $\pi(\x_i,y_i)$ is selected, then the objective function may be
dominated by this data point. Although this scenario happens with a very small
probability, the asymptotic variance of $\htheta_{W}$ will be greatly inflated
\citep{hesterberg1995weighted,owen2000safe,PingMa2014-JMLR}. However, if we use
$\htheta_{S}$, this problem will be ameliorated for the following two
reasons.
1) %
Since $\bar{\pi}(\x_i; \btheta)$ is an weighted average of $\pi(\x_i,y_i)$
across the conditional distribution of $y_i$ given $\x_i$, we know that
$\bar{\pi}(\x_i; \btheta)$'s are less variable than $\pi(\x_i,y_i)$'s. As a
result, even when $\pi(\x_i, y_i)$ is very close to zero,
$\bar{\pi}(\x_i; \btheta)$ can be bounded away from zero. 2) Even if
$\bar{\pi}(\x_i; \btheta)$ and $\pi(\x_i,y_i)$ approach to zero at the same
rate, $-\log\{\bar{\pi}(\x_i; \btheta)\}$ approaches to infinity much slower
than $\pi^{-1}(\x_i,y_i)$ does.  Compared with $\htheta_{W}$, $\htheta_{S}$ is
based on the conditional likelihood of the sampled data, and thus it has a
higher estimation efficiency. The price to pay for using $\htheta_{S}$ is
that the integration in $\bar{\pi}(\x_i; \btheta)$ depends on the model
structure, and $\bar{\pi}(\x_i; \btheta)$ may have complicated or no explicit
expressions for some models. On the other hand, $\htheta_{W}$ directly uses
$\pi^{-1}(\x_i,y_i)$ and {does not requires any extra model assumption and thus
it is  easier to implement}. We will show in
Section~\ref{sec:gener-line-models} that $\bar{\pi}(\x_i; \btheta)$'s have
a closed-form solution for many popular  generalized linear models, and
Newton's method is applicable with closed-form expressions for the score
functions and Hessian matrices.

Assuming that derivatives can pass integration, the score function associated with $\ell_{S}( \btheta)$ is
\begin{equation*}
  \dot\ell_{S} ( \btheta) = \frac{\partial}{ \partial \btheta}  \ell_{S}( \btheta) = \sum_{i =1}^N \delta_i \dot\ell( \btheta; \x_i, y_i) - \sum_{i=1}^N \delta_i \frac{ \partial \bar{\pi} (\x_i; \btheta)/ \partial \btheta }{ \bar{\pi} ( \x_i ; \btheta) },
\end{equation*}
where $\dot\ell( \btheta; \x, y) = \partial \log f( y \mid \x; \btheta)/ \partial \btheta$ and 
\begin{equation}
  \frac{ \partial}{ \partial \btheta} \bar{\pi}(\x_i;\btheta)
  = \int\frac{ \partial}{ \partial \btheta} f(y \mid \x_i; \btheta)  \pi(\x_i,y) d y
  = \Exp\{ \dot\ell( \btheta; \x_i, y_i) \pi(\x_i,y_i) \mid \x_i\}.
  \label{eq:36}
\end{equation}
Thus, we can express 
\begin{equation} 
  \dot\ell_{S}(\btheta) = \sumN\delta_i\left[
    \dot\ell(\btheta;\x_i,y_i)
    -\Exp\{\dot\ell(\btheta;\x_i,y_i)\mid\x_i,\delta_i=1\} \right],
\label{eqn:10} 
\end{equation} 
where 
\begin{equation}
 \Exp\{\dot\ell(\btheta;\x_i,y_i)\mid\x_i,\delta_i=1\} = 
  \frac{\Exp\{\dot\ell(\btheta;\x_i, y_i)\pi(\x_i, y_i)\mid\x_i\}}
  {\bar{\pi}(\x_i;\btheta)}.
\end{equation}
The second term in (\ref{eqn:10}) {can be called the bias-adjustment term for the score function.  %
}

\section{Asymptotic Results}
\label{sec:asymptotic-results}
We now presents some asymptotic results of the MSCLE of $\btheta$ proposed in (\ref{msle}).  To do this, we need the following regularity assumptions. As a convention in this paper, we use $\dot\ell$ and $\ddot\ell$ to denote gradient vector (of the first derivatives) and Hessian matrix (of the second derivatives) of {a function $\ell$ with respect to $\btheta$, respectively.}
\begin{assumption}\label{asmp:1}
  Assume that $\|\dot\ell(\btheta;\x,y)\|^2$ and $\|\ddot\ell(\btheta;\x,y)\|$ are integrable, 
  where $\|A\|=\tr^{1/2}(A\tp A)$ is the norm for a vector or matrix $A$. {Here, $\ell(\btheta;\x,y)=\log f(y\mid\x,\btheta)$ is the log-likelihood function of the original data distribution.}
\end{assumption}

\begin{assumption}\label{asmp:2}
  The parameter space $\bm\Theta$ is compact and the third order partial derivative of $\ell(\btheta;\x,y)$ and $\log\{\bar{\pi}(\x;\btheta)\}$ with respect to any components of $\btheta$ is bounded in absolute value by an integrable function $B(\x,y)$ that does not depend on $\btheta$.
\end{assumption}

\begin{assumption}\label{asmp:3}
  The matrix
  \begin{align}
    \bSigma_{\btheta}
    &=\Exp\bigg[\dot\ell^{\otimes2}(\btheta;\x,y)\pi(\x,y)
    -\frac{\Exp^{\otimes2}\{\dot\ell(\btheta;\x,y)\pi(\x,y)\mid\x\}}
    {\bar{\pi}(\x;\btheta)}\bigg]
      \label{eq:38}
  \end{align}
  is finite and positive definite, where $A^{\otimes2}=AA\tp$ for a vector or matrix $A$.
\end{assumption}

Compared with commonly used assumptions in maximum likelihood theory, the
above Assumptions~\ref{asmp:2} and \ref{asmp:3} impose additional constraints on
the subsampling probability $\pi(\x,y)$.
Assumption~\ref{asmp:2} imposes an integrable bound on
$\log\{\bar{\pi}(\x;\btheta)\}$ and its derivatives. This is to prevent the
distribution of $\pi(\x,y)$ from having a large probability around zero. Since
$\bar{\pi}(\x; \btheta)$ is an weighted average of $\pi(\x,y)$ across the
conditional distribution of $y$, it put
less probability around the boundary zero, so the required condition here is
less restrictive compared with those required by the IPW estimators for
which $\pi(\x,y)$ is in the denominator. For the IPW 
estimator, it is often assumed that $\pi(\x,y)$ is bounded away from zero
\citep[e.g.][]{ai2020optimal,wang2021optimal,yu2020quasi}, which is a much
stronger condition. 
If the subsampling probability is non-informative, i.e. $\pi(\x,y)=\pi(\x)$, then $\bar{\pi}(\x;\btheta)$ does not dependent on $\btheta$, and thus Assumption~\ref{asmp:2} reduces to the common condition that the third derivative of the log-likelihood is bounded in absolute value by an integrable random variable. 
Assumption~\ref{asmp:3} ensures that the variance covariance matrix for the sampled data conditional score function is positive definite. 
The key restriction here is the integrability of
$\bar{\pi}^{-1}(\x;\btheta)\Exp^{\otimes2}\{\dot\ell(\btheta;\x,y)\pi(\x,y)\mid\x\}$.
Again, this is less restrictive compared with the
conditions required by the IPW estimator, and a specific example will be
provided in Section~\ref{sec:subs-prob-based}. 
For non-informative subsampling, a sufficient condition for
  Assumption~\ref{asmp:3} is %
that $\Exp\{\dot\ell^{\otimes2}(\btheta;\x,y)\pi(\x)\}$ is finite and positive definite.  

We can also see that our conditions on the sampling probability in
Assumptions~\ref{asmp:2} and \ref{asmp:3} are less restrictive from another
angle. These conditions imposes restrictions on the distribution of $\x$ only,
while for IPW estimator with informative sampling probability the required
conditions impose restrictions on the distribution of $y$ as well
\citep[e.g.][etc]{WangZhuMa2018,wang2021optimal,WangZou2019,yu2020quasi}.

\begin{theorem}\label{thm:1}
Let $\{(\x_i, y_i), i=1, \cdots, N\}$ be $N$ independent realizations of $(X,Y)$ with joint density $f(y \mid \x; \btheta)f_X(\x)$ for some $\btheta \in \bm\Theta$ and $f_X(\x)$ is completely unspecified. Let $\htheta_{S}$ be the MSCLE of $\btheta$ defined in (\ref{msle}). 
  Under Assumptions~\ref{asmp:1}-\ref{asmp:3}, as $N$ goes to infinity, 
  \begin{equation}
    \sqrt{N}(\htheta_{S}-\btheta)\cvd\Nor(\0,\;\bSigma_{\btheta}^{-1}),
  \end{equation}
  where $\bSigma_{\btheta}$ is defined in (\ref{eq:38}), and $\cvd$ means convergence in distribution.
\end{theorem}
{In Theorem \ref{thm:1}, the subsample size $n^*=\sumN\delta_i$ is random, and the average subsample size $n=\Exp(n^*)=N\Exp\{\pi(\x,y)\}$ goes to infinity as $N \rightarrow \infty$.}

The estimator $\htheta_{S}$ is based on the conditional likelihood of the
sampled data, so it is expected to be more efficient than the IPW estimator.  
Actually, it is the most efficient estimator in a class of asymptotically
unbiased estimators.
To see this, consider the following class of estimating equations
\begin{equation} 
 \sumN \delta_i\u(\btheta;\x_i, y_i) =\0, 
 \label{ee}
\end{equation}
where  $\u(\btheta;\x,y)$ satisfies 
$\Exp\{\delta\u(\btheta;\x,y)\mid\x\}=\0$. Let $\htheta_u$ be the class of estimators
obtained through solving the class of estimating equations in (\ref{ee}). The
class of estimators defined via solving (\ref{ee}) includes the IPW estimator
and the MSCLE as special cases.  If $\u(\btheta;\x,
y)=\pi^{-1}(\x,y)\dot\ell(\btheta;\x,y)$, then $\htheta_u$ becomes the IPW
estimator defined in~(\ref{eq:4}); if $\u(\btheta;\x,y)=\dot\ell_{S}(\btheta;\x,y)= 
\dot\ell(\btheta;\x,y)
- \Exp\{\dot\ell(\btheta;\x,y)\mid\x,\delta=1\}$, then $\htheta_u$
becomes the MSCLE defined in~(\ref{msle}).

The following theorem shows that the MSCLE is the most efficient among the class of estimators defined through solving (\ref{ee}). 

\begin{theorem}\label{thm:2}
Assume that the partial derivatives of $\Exp\{\delta\u(\btheta;\x,y)\mid\x\}$ with respect to $\btheta$ can be passed under the integration sign, and that the regularity conditions for the following standard asymptotic expansion holds.   
\begin{equation} 
 \htheta_u = \btheta - \M_{\btheta}^{-1} N^{-1}\sumN\delta_i\u(\btheta; \x_i, y_i) + o_P(N^{-1/2}),
\label{result1}
\end{equation} 
where $\M_{\btheta}=\Exp\{\delta\dot\u(\btheta; \x, y)\}$ is full rank and $\dot\u(\btheta;\x, y)=\partial\u(\btheta;\x, y)/\partial\btheta\tp$. Assume that $\Exp\{\delta\u^{\otimes2}(\btheta;\x,y)\}$ exists. The asymptotic variance covariance matrix of $\htheta_u$ scaled by $N$ is $\M_{\btheta}^{-1}\Exp\{\delta\u^{\otimes2}(\btheta;\x,y)\}(\M_{\btheta}^{-1})\tp$, and it satisfies that
\begin{equation}\label{eq:59}
  \M_{\btheta}^{-1}\Exp\{\delta\u^{\otimes2}(\btheta;\x,y)\}(\M_{\btheta}^{-1})\tp
  \ge\bSigma_{\btheta}^{-1},
\end{equation}
in the Loewner ordering for any $\btheta$, where the equality holds if $\u(\btheta;\x, y)$ is a linear function of the subsampled data conditional score function, namely $\u(\btheta;\x, y)=-
    \M_{\btheta}\bSigma_{\btheta}^{-1}\dot\ell_{S}(\btheta;\x,y)$.
\end{theorem}

\begin{remark}\label{remark:3}
For the IPW estimator $\htheta_{W}$ in (\ref{eq:4}), let $\V_W$ denote the asymptotic variance scaled by $N$, which typically %
has a form of \citep[see][]{yu2020quasi}
\begin{equation}\label{eq:26}
  \V_{W}=\F^{-1}\Exp\bigg\{\frac{\dot\ell^{\otimes2}(\btheta;\x,y)}
  {\pi(\x,y)}\bigg\}\F^{-1},%
\end{equation}
where $\F=\Exp\{\dot\ell^{\otimes2}(\btheta;\x,y)\}$ is the Fisher information matrix of the original data distribution. 
From Theorem~\ref{thm:2}, %
$\V_{W}\ge\bSigma_{\btheta}
^{-1}$. %
{ Additionally, for $\V_{W}$ in (\ref{eq:26}), it requires $\pi^{-1}(\x,y)\dot\ell^{\otimes2}(\btheta;\x,y)$ to be integrable, which may be violated if $\pi(\x,y)$ has a high density in the neighborhood of zero. On the other hand, for $\bSigma_{\btheta}^{-1}$ in (\ref{eq:38}), it requires that $\bar{\pi}^{-1}(\x;\btheta)\Exp^{\otimes2}\{\dot\ell(\btheta;\x,y)
    \pi(\x,y)\mid\x\}$ is integrable. It is the average probability $\bar{\pi}(\x;\btheta)=\Exp\{\pi(\x,y)\mid\x\}$ that is in the denominator. The MSCLE is less restrictive compared with the IPW estimator because even if $\pi(\x,y)$ has a high density in the neighborhood of zero, $\bar{\pi}(\x;\btheta)$ may not be small.}
\end{remark}

\section{Estimated subsampling probabilities}
\label{sec:estim-param-subs}
In practice, informative subsampling probabilities may depend on unknown parameters, say $\bvtheta$, and a pilot subsample is often used to estimate it. Here $\bvtheta$ may be the same as $\btheta$ or contain $\btheta$ as its components. We denote the pilot estimator of $\bvtheta$ as $\tilde\bvtheta_{\plt}$, and assume that $\tilde\bvtheta_{\plt}$ is independent of the data to sample and converges to a limit, i.e., $\tilde\bvtheta_{\plt}\cvp\bvtheta_p$. Here $\cvp$ means convergence in probability. 

The assumption that a pilot estimator is independent of the data is commonly used in the literature \citep[e.g.,][]{fithian2014local,Han2019}, and it is reasonable in the context of subsampling. If one uses simple random sampling to obtain a pilot subsample or simply use the first certain number of observations as a pilot sample, then the pilot subsample is independent of the rest of the data. Thus, taking the rest of the data as the full data and performing informative subsampling, the pilot estimator is independent of the full data. Since one is likely to assign a larger pilot subsample size for a larger original data sample size, it is reasonable to assume that the pilot estimator converges to some limit. Here we do not have to assume that the pilot estimator converges to the true parameter, namely, the pilot estimator can be misspecified.
As discussed in Section~\ref{sec:introduction}, because the subsampling probabilities are under our control,  we can always use 
the realized subsampling probabilities in constructing the sampled likelihood. Thus, even if the pilot samples are systematically different from the original sample, we can still use the realized subsampling probabilities and the  
 statistical properties of the maximum sampled conditional likelihood estimator remain valid. {If the pilot sample is systematically different from the original sample, the pilot estimator may not be consistent to the true parameter and the resulting subsampling probabilities may not be optimal anymore. We will show numerically in Section~\ref{sec:numer-exper} that the proposed MSCLE is more robust to pilot misspecification compared with the IPW estimator.}

Another practical consideration is the subsample size. In Section~\ref{sec:sampl-data-likel}, the average subsample size is $n=N\Exp\{\pi(\x,y)\}$ which is the same order of $N$. Since the intended average subsample size may be much smaller than the full data sample size in practice, in this section we use $n$ to denote the average subsample size, namely $\Exp\{\pi_{N}(\x_i,y_i;\bvtheta)\}=n/N$, 
and allow $n=o(N)$. This means the subsampling probabilities are dependent on $N$ and may go to zero. 

Given an estimated pilot $\tilde{\bvtheta}_{\plt}$, the sampled data conditional log-likelihood function is written as
\begin{equation}\label{eq:60}
  \ell_{S}(\btheta\mid\tilde\bvtheta_{\plt})
  =\sumN\delta_i\big[\log f(y_i\mid\x_i;\btheta)
  -\log\{\bar{\pi}_{N}(\x_i;\btheta\mid\tilde\bvtheta_{\plt})\}\big]+C,
\end{equation}
where 
\begin{equation}\label{eq:60b} 
   \bar{\pi}_{N}(\x_i;\btheta\mid \bvtheta)
  =\int f(y\mid\x_i;\btheta)\pi_{N}(\x_i,y; \bvtheta )\ud y,
\end{equation}
and $C=\sumN\delta_i\log\{\pi_{N}(\x_i,y_i;\tilde\bvtheta_{\plt})\}$ does not contain $\btheta$. Here, we use notation conditional on $\tilde\bvtheta_{\plt}$ to emphasize its dependence on $\tilde\bvtheta_{\plt}$.

Denote the sampled data estimator through maximizing $\ell_{S}(\btheta\mid\tilde\bvtheta_{\plt})$ as $\htheta_{s,\tilde\bvtheta_{\plt}}$. We need the following regularity assumptions to investigate the asymptotic distribution of $\htheta_{s,\tilde\bvtheta_{\plt}}$.

\begin{assumptionp}{1'}\label{asmp:1p}
  For $\bvtheta$ in an neighborhood of $\bvtheta_p$
  \begin{align}
    &\limsup_{n,N\rightarrow\infty}\ \frac{N}{n}\Exp\{\pi_{N}(\x,y;\bvtheta)
      \|\dot\ell(\btheta;\x,y)\|^4\}<\infty \label{eq:34}
      \quad\text{and}\quad\\
    &\limsup_{n,N\rightarrow\infty}\ \frac{N}{n}\Exp\{\pi_{N}(\x, y,\bvtheta)
      \|\ddot\ell(\btheta;\x,y)\|^2\}<\infty.\label{eq:33}
  \end{align}
\end{assumptionp}

\begin{assumptionp}{2'}\label{asmp:2p}
  The parameter space $\bm\Theta$ is compact and there exist a function $B_{\bvtheta}(\x,y)$ such that for any component of $\btheta$, say $\theta_{j_1}$, $\theta_{j_2}$, and $\theta_{j_3}$,
  \begin{align}\label{eq:61}%
    \sup_{\btheta\in\bm\Theta}\bigg|\frac{\partial^3\ell(\btheta;\x,y)}
    {\partial\theta_{j_1}\partial\theta_{j_2}\partial\theta_{j_3}}\bigg|
    \le&B_{\bvtheta}(\x,y),\\
    \sup_{\btheta\in\bm\Theta}\bigg|\frac{\partial^3
    \log\{\bar{\pi}_{N}(\x;\btheta\mid\bvtheta)\}}
    {\partial\theta_{j_1}\partial\theta_{j_2}\partial\theta_{j_3}}\bigg|
    \le&B_{\bvtheta}(\x,y),\label{eq:32}
  \end{align}
  where $\bar{\pi}_{N}(\x;\btheta\mid\bvtheta)$ is defined in (\ref{eq:60b}) 
  and $B_{\bvtheta}(\x,y)$ satisfies that for $\bvtheta$ in a neighborhood of $\bvtheta_p$
\begin{align}\label{eq:35}
  \limsup_{n,N\rightarrow\infty}\
  \frac{N}{n}\Exp\{\pi_{N}(\x,y;\bvtheta)B_{\bvtheta}(\x,y)\}<\infty.
\end{align}

\end{assumptionp}

\begin{assumptionp}{3'}\label{asmp:3p}
  As $n$ and $N$ goes to infinity, for $\bvtheta$ in an neighborhood of $\bvtheta_p$, the matrix
  \begin{align}\label{eq:50}
    \bSigma_{N,\btheta,\bvtheta}
    &=\frac{N}{n}\Exp\bigg[
    \dot\ell^{\otimes2}(\btheta;\x,y)\pi_{N}(\x,y;\bvtheta)
    -\frac{\Exp^{\otimes2}\{\dot\ell(\btheta;\x,y)
    \pi_{N}(\x,y;\bvtheta)\mid\x\}}
    {\Exp\{\pi_{N}(\x,y;\bvtheta)\mid\x\}}\bigg]
    \rightarrow\bSigma_{\btheta,\bvtheta},
  \end{align}
  where $\bSigma_{N,\btheta,\bvtheta}$ and $\bSigma_{\btheta,\bvtheta}$ are finite, positive definite, and continuous with respective to $\bvtheta$.
\end{assumptionp}

\begin{remark}\label{remark:4}
  Assumption~\ref{asmp:1p} is essentially moment conditions on the first and second derivatives of the log-likelihood. If $n^{-1}N\pi_{N}(\x,y;\bvtheta)$ are bounded, then the integrability of $\|\dot\ell(\btheta;\x,y)\|^4$ and $\|\ddot\ell(\btheta;\x,y)\|^2$ is sufficient for Assumption~\ref{asmp:1p}. We impose stronger moment conditions here compared with the independent and identically distributed (i.i.d.) case in Section~\ref{sec:sampl-data-likel}, because we allow $\pi_{N}(\x,y;\bvtheta)$ to depend on $N$ and $\bvtheta$ has to be estimated. Assumptions~\ref{asmp:2p}-\ref{asmp:3p} are the counterparts corresponding to the Assumptions~\ref{asmp:2}-\ref{asmp:3} in Section~\ref{sec:sampl-data-likel}. 
 Note that $\bSigma_{N,\btheta,\bvtheta}$ is always semi-positive definite because it can be written as
  \begin{equation}
    \bSigma_{N,\btheta,\bvtheta}=\frac{N}{n}\Exp\Bigg(
    \pi_{N}(\x,y;\bvtheta)\bigg[\dot\ell(\btheta;\x,y)
    -\frac{\Exp\{\dot\ell(\btheta;\x,y)\pi_{N}(\x,y;\bvtheta)\mid\x\}}
    {\Exp\{\pi_{N}(\x,y;\bvtheta)\mid\x\}}\bigg]^{\otimes2}\Bigg).
  \end{equation}
If $n^{-1}N\pi_{N}(\x,y;\bvtheta)$ are  bounded away from zero, then with
Assumption~\ref{asmp:1p}, a sufficient condition for Assumption~\ref{asmp:3p} is
that the sequence of random variables
\begin{equation*}
  \dot\ell(\btheta;\x,y)
  -\frac{\Exp\{\dot\ell(\btheta;\x,y)\pi_{N}(\x,y;\bvtheta)\mid\x\}}
  {\Exp\{\pi_{N}(\x,y;\bvtheta)\mid\x\}}
\end{equation*}
is full rank, for which a sufficient condition is that $\pi_{N}(\x,y;\bvtheta)$
and $\bm c\tp\dot\ell(\btheta;\x,y)$ has nonzero correlation for any nonzero and
nonrandom vector $\bm c$. %
If $\pi_{N}(\x,y;\bvtheta)$ is obtained by rescaling a function of $(\x,y)$
    i.e., $\pi_{N}(\x,y;\bvtheta)=n/N\pi(\x,y;\bvtheta)$ (e.g., the case in Section~\ref{sec:subs-prob-based}),
where $\pi(\x,y;\bvtheta)$
does not depend on $N$, then Assumption~\ref{asmp:3p} reduces to the same
requirement as in Assumption~\ref{asmp:3} except that $\pi(\x,y)$ is replaced by
$\pi(\x,y;\bvtheta)$. If $\pi(\x,y;\bvtheta)$ is bounded away from zero, then a
sufficient condition is the integrability of  
\begin{align*}
  \Exp^{-1}\{\pi(\x,y;\bvtheta)\mid\x\}
  \Exp^{2}\{\|\dot\ell(\btheta;\x,y)\|\pi(\x,y;\bvtheta)\mid\x\}
\end{align*}
together with a nonzero correlation between $\pi(\x,y;\bvtheta)$ and $\bm
c\tp\dot\ell(\btheta;\x,y)$ for any nonzero constant vector $\bm c$.
\end{remark}

\begin{theorem}\label{thm:1p}
  Under Assumptions~\ref{asmp:1p}-\ref{asmp:3p}, as $n$ and $N$ goes to infinity, 
  \begin{equation}
    \sqrt{n}(\htheta_{s,\tilde\bvtheta_{\plt}}-\btheta)\cvd\Nor(\0,\;\bSigma_{\btheta,\bvtheta_p}^{-1}),
  \end{equation}
  where $\bSigma_{\btheta,\bvtheta_p}$ is the limit of $\bSigma_{N,\btheta,\bvtheta_p}$.
\end{theorem}

Let $\htheta_{u,\tilde\bvtheta_{\plt}}$ be the estimator obtained through solving the following class of estimating equations
\begin{equation} 
 \sumN \delta_i\u_{\tilde\bvtheta_{\plt}}(\btheta;\x_i, y_i) =\0, 
 \label{ee2}
\end{equation}
where $\u_{\tilde\bvtheta_{\plt}}(\btheta;\x, y)$ satisfies that $\Exp\{\delta\u_{\tilde\bvtheta_{\plt}}(\btheta;\x,y)\mid\x,\tilde\bvtheta_{\plt}\}=\0$. The class of estimators defined in~(\ref{ee2}) includes the IPW estimator with $\u_{\tilde\bvtheta_{\plt}}(\btheta;\x, y)=\pi_{N}^{-1}(\x,y;\tilde\bvtheta_{\plt})\dot\ell(\btheta;\x,y)$ and the MSCLE with $\u_{\tilde\bvtheta_{\plt}}(\btheta;\x, y)=\dot\ell_{S}(\btheta;\x,y\mid\tilde\bvtheta_{\plt})$.    

Similar to Theorem~\ref{thm:2}, the following result shows that with estimated parameters in subsampling probabilities, the proposed estimator is the most efficient among the class of estimators defined through solving (\ref{ee2}). 

\begin{theorem}\label{thm:2p}
Assume that the partial derivatives of $\Exp\{\delta\u_{\tilde\bvtheta_{\plt}}(\btheta;\x,y)\mid\x,\tilde\bvtheta_{\plt}\}$ with respect to $\btheta$ can be passed under the integration sign, and the regularity conditions for the following standard asymptotic expansion holds.   
\begin{equation} 
 \htheta_{u,\tilde\bvtheta_{\plt}} = \btheta - n^{-1}\M_{\btheta,\bvtheta_p}^{-1}\sumN\delta_i\u_{\tilde\bvtheta_{\plt}}(\btheta; \x_i, y_i) + o_P(n^{-1/2}),
\end{equation} 
where $\M_{N,\btheta,\bvtheta}=n^{-1}N\Exp\{\delta\dot\u_{\bvtheta}(\btheta; \x, y)\}$
and $\dot\u(\btheta;\x, y)=\partial\u(\btheta;\x, y)/\partial\btheta\tp$.
Assume that $\M_{N,\btheta,\bvtheta}$ and $\V_{N,\btheta,\bvtheta}=n^{-1}N\Exp\{\delta\u_{\bvtheta}^{\otimes2}(\btheta;\x,y)\}$ are continuous in $\bvtheta$ and they converge to finite and full rank matrices $\M_{\btheta,\bvtheta}$ and $\V_{\btheta,\bvtheta}$, respectively. 
The asymptotic variance covariance matrix of $\htheta_u$ (multiplied by $n$) is $\M_{\btheta,\bvtheta_p}^{-1}\V_{\bvtheta_p}\M_{\btheta,\bvtheta_p}^{-1}$, and it satisfies that
\begin{equation}
  \M_{\btheta,\bvtheta_p}^{-1}\V_{\bvtheta_p}\M_{\btheta,\bvtheta_p}^{-1}
  \ge\bSigma_{\btheta,\bvtheta_p}^{-1},
\end{equation}
where the equality holds if $\u_{\tilde\bvtheta_{\plt}}(\btheta; \x, y)$ is asymptotically a linear function of $\dot\ell_{S}(\btheta;\x,y\mid\tilde\bvtheta_{\plt})$ in the sense that 
$n^{-1}N\Exp(\delta\|T_{N\tilde\bvtheta_{\plt}}\|^2\mid\tilde\bvtheta_{\plt})=\op$, where
$T_{N\tilde\bvtheta_{\plt}}=\u_{\tilde\bvtheta_{\plt}}(\btheta;\x,y)
    +\M_{\btheta,\bvtheta_p}\bSigma_{\btheta,\bvtheta_p}^{-1}
    \dot\ell_{S}(\btheta;\x,y\mid\tilde\bvtheta_{\plt})$.
\end{theorem}

\section{Subsampling probabilities based on gradient norms (GN)}
\label{sec:subs-prob-based}

We use a specific form of subsampling probabilities to illustrate our results. 
One way to specify subsampling probabilities is to use the norm of the per-observation score function, the gradient of the log-likelihood, by letting
\begin{equation}\label{eq:45}
  \pi(\x_i,y_i)\propto \|\dot\ell(\btheta;\x_i,y_i)\|,
  \quad\text{for}\quad i=1, ..., N,
\end{equation}
where $\propto$ means ``proportional to''. This option leads to the L-optimal subsampling probabilities that minimize the trace of the asymptotic variance covariance matrix for a linearly transformed IPW estimator \citep{WangZhuMa2018,ai2020optimal,yu2020quasi}. 
While the optimal probabilities in (\ref{eq:45}) are for the IPW estimator in (\ref{eq:4}), it can be used to obtain more efficient subsample for our MSCLE. 

Since (\ref{eq:45}) implies that the subsampling probabilities are dependent on the unknown $\btheta$, a pilot estimate is required. In addition, to control the average subsample size, the subsampling probabilities need to be re-scaled and the scaling value may also be unknown. Let $\bvtheta$ be the vector consisting of parameters that are required to calculate the subsampling probabilities. The subsampling probabilities are presented as 
\begin{equation}\label{eq:37}
  \pi_{N}(\x_i,y_i;\bvtheta)
  =\frac{n\|\dot\ell(\btheta;\x_i,y_i)\|}
  {N\Psi},%
  \quad i=1, ..., N,
\end{equation}
where $\bvtheta=(\btheta\tp,\Psi)\tp$ and $\Psi=\Exp\{\|\dot\ell(\btheta;\x,y)\|\}$. In practice, the unknown $\bvtheta$ is replaced by a pilot estimate $\tilde\bvtheta_{\plt}$. 

If the subsampling ratio $n/N$ is large (far from zero), then it is possible that some $\pi_{N}(\x_i,y_i;\bvtheta_p)>1$ in (\ref{eq:37}), and therefore the resulting average subsample size may be smaller than $n$. %
However, in subsampling, the subsample size is typically much smaller than the full data sample size, so it is reasonable to assume that $n=o(N)$. In this case, the truncation can be ignored, and the average subsample size is $n$ asymptotically, namely,
\begin{equation*}
  \frac{N}{n}\Exp\bigg\{
  \frac{n\|\dot\ell(\btheta;\x_i,y_i)\|}{N\Psi}\wedge1
  \bigg\}\rightarrow1,
\end{equation*}
where $a\wedge b$ means the smaller value of $a$ and $b$. 
For the rest of the paper, we assume that $n=o(N)$, and we call $\pi_{N}(\x_i,y_i;\bvtheta)$ in (\ref{eq:37}) subsampling probabilities since the number of cases that $\pi_{N}(\x_i,y_i;\bvtheta)>1$ is negligible in this scenario. 

With the specific form of $\pi_{N}(\x_i,y_i;\bvtheta)$ in (\ref{eq:37}), the sampled data conditional log-likelihood function given pilot estimate $\tilde\bvtheta_{\plt}=(\ttheta_{\plt}\tp,\tPsi_{\plt})\tp$ can be written specifically as
\begin{equation}
  \ell_{S}(\btheta\mid\tilde\bvtheta_{\plt})
  =\sumN\delta_i\big[\log f(y_i\mid\x_i;\btheta)
  -\log\Exp\{\|\dot\ell(\ttheta_{\plt};\x_i,y_i)\|\mid\x_i,\tilde\bvtheta_{\plt}\} \big]+C,
\end{equation}
where $C$ does not contain $\btheta$.
 Note that the second term contains $\btheta$ through the expectation.
Accordingly, a sufficient condition for Assumption~\ref{asmp:1p} is that $\|\dot\ell(\grave\btheta;\x,y)\|\|\ddot\ell(\btheta;\x,y)\|^2$ and $\|\dot\ell(\grave\btheta;\x,y)\|\|\dot\ell(\btheta;\x,y)\|^4$ are integrable for $\grave\btheta$ in the neighborhood of $\btheta_p$, and a sufficient condition for (\ref{eq:35}) in Assumption~\ref{asmp:2p} is that $\|\dot\ell(\grave\btheta;\x,y)\|B_{\btheta}(\x,y)$ is integrable for $\grave\btheta$ in the neighborhood of $\btheta_p$. These sufficient conditions for Assumptions~\ref{asmp:1p} and \ref{asmp:2p} with the specific subsampling probabilities $\pi_{N}(\x_i,y_i;\bvtheta)$ are essentially moments conditions on the log-likelihood of the original data distribution and its derivatives. %
Assumption~\ref{asmp:3p} requires that %
\begin{align}
    \bSigma_{\btheta,\bvtheta_p}
    &=\frac{1}{\Psi_p}
    \Exp\bigg[\dot\ell^{\otimes2}(\btheta;\x,y)
    \|\dot\ell(\btheta_p;\x,y)\|
    -\frac{\Exp^{\otimes2}\{\dot\ell(\btheta;\x,y)
    \|\dot\ell(\btheta_p;\x,y)\|\mid\x\}}
    {\Exp\{\|\dot\ell(\btheta_p;\x,y)\|\mid\x\}}\bigg]
        \label{eq:39}
  \end{align}
  is finite and positive definite, %
  and a sufficient condition for $\bSigma_{\btheta,\bvtheta_p}$ in (\ref{eq:39}) to be positive definite is that
  the quantity $\l\tp\dot\ell(\btheta;\x,y)$ is dependent on $y$ for any constant vector $\l\neq\0$. %

For the IPW estimator, 
under some regularity conditions, the asymptotic variance-covariance matrix (multiplied by $n$) is  
\begin{equation}\label{eq:53}
  \V_{W,\bvtheta_p}=\Psi_p\F^{-1}
  \Exp\bigg\{\frac{\dot\ell^{\otimes2}(\btheta;\x,y)}
  {\|\dot\ell(\btheta_p;\x,y)\|}\bigg\}\F^{-1}.
\end{equation}
{Comparing the expressions of $\bSigma_{\btheta,\bvtheta_p}^{-1}$ in~(\ref{eq:39}) and $\V_{W,\bvtheta_p}$ in (\ref{eq:53}), we see that the MSCLE requires weaker assumptions than the IPW estimator. For example, for a normal linear regression model with a known error variance, say $\sigma^2=1$,  $\dot\ell(\btheta;\x,y)=-0.5(y-\x\tp\btheta)\x$. If $\btheta_p\neq\btheta$, i.e., the pilot estimate is not consistent, then given $\x$,
\begin{equation}
  \frac{\dot\ell_{j}^2(\btheta;\x,y)}{\|\dot\ell(\btheta_p;\x,y)\|}
  =\frac{(y-\x\tp\btheta)^2x_j^2}{2|y-\x\tp\btheta_p|\|\x\|}
\end{equation}
is not integrable; its expectation exists and equals $+\infty$. Thus, any element of $\V_{W,\bvtheta_p}$ (if the integral exists) will be $\pm\infty$. This means that for a subsample taken according to (\ref{eq:37}), the IPW estimator is not applicable for linear regression.  %
On the other hand for the term in the denominator of (\ref{eq:39}), by direct calculation we know that
\begin{equation}
  \Exp\{\|\dot\ell(\btheta_p;\x,y)\|\mid\x\}
  =\Exp(|y-\x\tp\btheta_p|\mid\x)\|\x\|
  =\{2(\mu_d)\Phi(\mu_d)+2\phi(\mu_d)-\mu_d\}\|\x\|,
\end{equation}
where $\mu_d=\x\tp(\btheta_p-\btheta)$, and $\Phi$ and $\phi$ are the cumulative distribution function and the probability density function, respectively, of the standard normal distribution. Thus under some integrability requirement on the covariate $\x$, the MSCLE is applicable.}

From Theorem~\ref{thm:2p}, we know that $\bSigma_{\btheta,\bvtheta_p}^{-1}\le\V_{W,\bvtheta_p}$. Actually, for the IPW estimator, we have
\begin{align}
  \V_{W,\bvtheta_p}-\bSigma_{\btheta,\bvtheta_p}^{-1}
  =\Psi_p\Exp\big\{\|\dot\ell(\btheta_p;\x,y)\|\bm\xi^{\otimes2}\big\},
\end{align}
where
\begin{align}
  \bm\xi=\bSigma_{\btheta,\bvtheta_p}^{-1}\bigg[
  \dot\ell(\btheta;\x,y)
  -\frac{\Exp\{\dot\ell(\btheta;\x,y)
  \|\dot\ell(\btheta_p;\x,y)\|\mid\x\}}
  {\Exp\{\|\dot\ell(\btheta_p;\x,y)\|\mid\x\}}\bigg]
  -\F^{-1}\frac{\dot\ell(\btheta;\x,y)}{\|\dot\ell(\btheta_p;\x,y)\|}.
\end{align}
This implies that $\bSigma_{\btheta,\bvtheta_p}^{-1}=\V_{W,\bvtheta_p}$ if and only if $\bm\xi=\0$ almost surely, which occurs if and only if $\dot\ell(\btheta_p;\x,y)=\0$ almost surely and this is not possible. Thus with the subsampling probabilities $\pi_{N}(\x_i,y_i;\bvtheta)$'s in (\ref{eq:37}) $\bSigma_{\btheta,\bvtheta_p}^{-1}<\V_{W,\bvtheta_p}$.

The $\pi_{N}(\x_i,y_i;\bvtheta)$'s in (\ref{eq:37}) are a version of the L-optimal subsampling probabilities. If the pilot $\tilde\bvtheta_{\plt}$ is consistent, then $\pi_{N}(\x_i,y_i;\tilde\bvtheta_{\plt})$'s minimize the trace of the asymptotic variance-covariance matrix of the IPW estimator of $\F\btheta$ among all subsampling probabilities with the same average subsample size. If the dimension of $\btheta$ is one, then $\pi_{N}(\x_i,y_i;\bvtheta)$'s minimize the asymptotic variance of the IPW estimator of $\btheta$. Our results show that the estimation efficiency can be further improved by using the proposed MSCLE.

\section{Generalized linear models} 
\label{sec:gener-line-models}
We provide more detailed discussions to illustrate the proposed estimator in the context of generalized linear models (GLMs). As before, we use $\dot\ell$ and $\ddot\ell$, respectively, to denote gradient vector and Hessian matrix of $\ell$ with respect to a vector variable; and we use $b'$ and $b''$, respectively, to denote the first derivative and the second derivative of a function $b$ with respect to a scalar variable.

Let $y_i$ be the response and $\x_i$ be the corresponding covariate. A GLM assumes that the conditional mean of the response $y_i$ give the covariate $\x_i$, $\mu_i=\Exp(y_i\mid\x_i)$, satisfies
\begin{equation*}
  g(\mu_i)=g\{\Exp(y_i\mid\x_i)\}=\x_i\tp\bbeta,
\end{equation*}
  where $g$ is the link function, $\x_i\tp\bbeta$ is the linear predictor, and $\bbeta$ is the regression coefficient. For %
  commonly used GLMs, it is assumed that the distribution of the response $y_i$ given the covariate $\x_i$ belongs to the exponential family, namely,
\begin{equation}\label{eq:10}
  f(y_i\mid\x_i;\bbeta,\phi)=a(y_i,\phi)\exp\Big\{
  \frac{y_ib(\x_i\tp\bbeta)-c(\x_i\tp\bbeta)}{\phi}\Big\},
\end{equation}
where $a$, $b$ and $c$ are known scalar functions, and $\phi$ is the dispersion parameter. In the framework of GLMs, if the link function $g$ is selected such that $b$ is the identity function, i.e., $b(\x_i\tp\bbeta)=\x_i\tp\bbeta$, then the link function is called the canonical link. With a canonical link function, $g\{\Exp(y_i\mid\x_i)\}=c'(\x_i\tp\bbeta)$ where $c'$ is the derivative function of $c$. For example, binary logistic regression is a special case of GLMs with the canonical link and the response variable follows the Bernoulli distribution. Specifically, in logistic regression, $a(y_i,\phi)=1$, $b(\x_i\tp\bbeta)=\x_i\tp\bbeta$, $c(\x_i\tp\bbeta)=\log\{1+\exp(\x_i\tp\bbeta)\}$, and $\phi=1$. 
Multi-class logistic regression can also be formulated into the family of GLMs with multivariate responses. We will provide more details in Section~\ref{sec:example:-multi-class}.

For a GLM from the exponential family, %
the MLE of $\bbeta$ is the maximizer of
\begin{equation*}
  \ell(\bbeta)=\sumN\ell(\bbeta;\x_i,y_i)
  =\frac{1}{\phi}\sumN\{y_ib(\x_i\tp\bbeta)-c(\x_i\tp\bbeta)\}+C,
\end{equation*}
where $C$ does not contain $\bbeta$, and the maximizer is the solution to the score equation
\begin{align}\label{eq:41}
  \dot\ell(\bbeta)=\sumN\dot\ell(\bbeta;\x_i,y_i)
  =\frac{1}{\phi}\sumN(y_i-\mu_i)b'(\x_i\tp\bbeta)\x_i=\0,
\end{align}
$\mu_i=c'(\x_i\tp\bbeta)/b'(\x_i\tp\bbeta)$, and $b'$ and $c'$ are the first derivative functions of $b$ and $c$, respectively. 

\subsection{Informative subsampling estimation}

Multiple informative subsampling designs are available such as the local case-control \citep{fithian2014local}, the A- and L-optimal subsampling \citep{ai2020optimal,wang2019more}, and the local uncertainty subsampling \citep{Han2019}. %
In GLMs with univariate responses, these probabilities have a unified expression and they satisfy that 
\begin{equation}\label{eq:7}
  \pi_{N}(\x_i,y_i;\bvtheta)\propto
  |y_i-\mu_i|\big|b'(\x_i\tp\bbeta)\big|h(\x_i),
  \quad i=1, ..., N,
\end{equation}
where $h(\x_i)>0$ is a criterion function that may or may not depend on
$\bbeta$. For example, in logistic regression, if $h(\x_i)=1$, then
$\pi_{N}(\x_i,y_i;\bvtheta)$ corresponds to the local case-control subsampling;
if $h(\x_i)=\|\x_i\|$, then $\pi_{N}(\x_i,y_i;\bvtheta)$ corresponds to the
L-optimal subsampling discussed in Section~\ref{sec:subs-prob-based};
if $h(\x_i)=\|\ddot\ell^{-1}(\bbeta)\x_i\|$ with $\ddot\ell(\bbeta)$ being the
Hessian matrix of $\ell(\bbeta)$, then $\pi_{N}(\x_i,y_i;\bvtheta)$ corresponds
to the A-optimal subsampling and $h(\x_i)$ depends on $\btheta$ in
general for this case.

In practical implementation, a pilot estimate of $\bvtheta$ is used, and to control the average subsample size as $n$ the subsampling probabilities are often taken as
\begin{align}\label{eq:3}
  \pi_{N}(\x_i,y_i;\tilde\bvtheta_{\plt})
  =\frac{n}{N\tPsi_{\plt}}
  |y_i-\tmu_i|\big|b'(\x_i\tp\tbeta_{\plt})\big|h(\x_i),
  \quad i=1, ..., N,
\end{align}
where $\tmu_i=c'(\x_i\tp\tbeta_{\plt})/b'(\x_i\tp\tbeta_{\plt})$, and $\tbeta_{\plt}$ and  $\tPsi_{\plt}$ are pilot estimates of parameters $\bbeta$ and $\Exp\{|(y-\mu)b'(\x\tp\bbeta)|h(\x)\}$, respectively.

For a subsample taken according to $\pi_{N}(\x_i,y_i;\tilde\bvtheta_{\plt})$'s, by direct calculation, the sampled data log-likelihood function for $\bbeta$ simplifies to
\begin{equation}\label{eq:11}
  \ell_{S}(\bbeta\mid\tilde\bvtheta_{\plt})
  =\frac{1}{\phi}\sumN\delta_i\big\{y_ib(\x_i\tp\bbeta)-c(\x_i\tp\bbeta)
    -\phi\log\Exp(|y_i-\tmu_i|\mid\x_i;\tbeta_{\plt})\big\}+C,
\end{equation}
where $C$ do not contain $\bbeta$. 
By direct calculation, we know that the MSCLE can be obtained by solving the score equation,
\begin{equation}\label{eq:17}
  \dot\ell_{S}(\bbeta\mid\tilde\bvtheta_{\plt})
  =\frac{1}{\phi}\sumN\delta_i
  \Big(y_i-\frac{\tilde\kappa_{1,i}}{\tilde\kappa_{0,i}}\Big)
  b'(\x_i\tp\bbeta)\x_i=\0,
\end{equation}
where $\tilde\kappa_{0,i}=\Exp(|y_i-\tmu_i|\mid\x_i;\tbeta_{\plt})$ and $\tilde\kappa_{1,i}=\Exp(y_i|y_i-\tmu_i|\mid\x_i;\tbeta_{\plt})$. 
The Hessian matrix is
\begin{align}
  \ddot\ell_{S}(\bbeta\mid\tilde\bvtheta_{\plt})
  &=\frac{1}{\phi}\sumN\delta_i
    \Big(y_i-\frac{\tilde\kappa_{1,i}}{\tilde\kappa_{0,i}}\Big)
    b''(\x_i\tp\bbeta)\x_i\x_i\tp\notag\\
  &\quad-\frac{1}{\phi^2}\sumN\delta_i
    \Big(\frac{\tilde\kappa_{2,i}}{\tilde\kappa_{0,i}}
    -\frac{\tilde\kappa_{1,i}^2}{\tilde\kappa_{0,i}^2}\Big)
    \{b'(\x_i\tp\bbeta)\}^2\x_i\x_i\tp.\label{eq:54}
\end{align}
where $\tilde\kappa_{2,i}=\Exp(y_i^2|y_i-\tmu_i|\mid\x_i;\tbeta_{\plt})$. 
Note that on the right-hand-side of (\ref{eq:54}), the second term is the dominating term, and the first term is often ignored in numerical optimization, namely, using $\ddot\ell_{S,\tilde\bvtheta_{\plt}}^{F}(\bbeta)$ instead of $\ddot\ell_{S}(\bbeta\mid\tilde\bvtheta_{\plt})$, where
\begin{equation}\label{eq:19}
  \ddot\ell_{S}^F(\bbeta\mid\tilde\bvtheta_{\plt})
  =-\frac{1}{\phi^2}\sumN\delta_i
    \Big(\frac{\tilde\kappa_{2,i}}{\tilde\kappa_{0,i}}
    -\frac{\tilde\kappa_{1,i}^2}{\tilde\kappa_{0,i}^2}\Big)
    \{b'(\x_i\tp\bbeta)\}^2\x_i\x_i\tp.
\end{equation}
Here we use the superscript $^F$ because the resulting form of Newton's method is often called the Fisher scoring algorithm. 
If the canonical link is used as mostly implemented in practice, then $b'(\cdot)=1$ and $b''(\cdot)=0$, and thus $\ddot\ell_{S}^F(\bbeta\mid\tilde\bvtheta_{\plt})=\ddot\ell_{S}(\bbeta\mid\tilde\bvtheta_{\plt})$. 
The MSCLE can be calculated from the Fisher scoring algorithm %
by iteratively applying
\begin{equation}
  \bbeta^{(t+1)} = \bbeta^{(t)}
  - \{\ddot\ell_{S}^{F}(\bbeta^{(t)}\mid\tilde\bvtheta_{\plt})\}^{-1}
  \dot\ell_{S}(\bbeta^{(t)}\mid\tilde\bvtheta_{\plt}).
\end{equation}
We will show the explicit expressions of $\tilde\kappa_{0,i}$, $\tilde\kappa_{1,i}$, and $\tilde\kappa_{2,i}$ for the examples in the following subsections. 

For GLMs with the subsampling probabilities defined in (\ref{eq:3}), the $\bSigma_{\btheta,\bvtheta_p}$ has the following specific form.
\begin{align*}
  \bSigma_{\btheta,\bvtheta_p}
  &=\frac{1}{\phi^2\Psi_p}\Exp\Bigg(
  b'(\x\tp\bbeta_p)\{b'(\x\tp\bbeta)\}^2h(\x)
    \bigg[|y-\mu_p|(y-\mu)^2
    -\frac{\Exp^2\{(y-\mu)|y-\mu_p|\mid\x\}}{\Exp\{|y-\mu_p|\mid\x\}}
  \bigg]\x^{\otimes2}\Bigg),
\end{align*}
where $\mu=c'(\x\tp\bbeta)/b'(\x\tp\bbeta)$ and $\mu_p=c'(\x\tp\bbeta_p)/b'(\x\tp\bbeta_p)$. 
With a canonical link and a consistent pilot estimate, $\bSigma_{\btheta,\bvtheta_p}$ simplifies to
\begin{align*}
  \bSigma_{\btheta,\bvtheta_p}
  &=\frac{1}{\phi^2\Exp\{|y-\mu|h(\x)\}}
    \Exp\Bigg(\bigg[|y-\mu|^3
    -\frac{\Exp^2\{(y-\mu)|y-\mu|\mid\x\}}{\Exp\{|y-\mu|\mid\x\}}
  \bigg]h(\x)\x^{\otimes2}\Bigg).
\end{align*}

\subsection{Examples}
\label{sec:examples}
\subsubsection{Example: binary response models}
\label{sec:exampl-binary-resp}
Binary data are ubiquitous in case-control studies and classifications. Let $y\in\{0,1\}$ be the binary response variable and $\x$ be the covariate vector. Assume that the probability for $y=1$ given $\x$ is
\begin{equation}\label{eq:55}
  \Pr(y=1\mid\x;\bbeta) = p(\x\tp\bbeta),
\end{equation}
where $p(\cdot)$ is a smooth and monotonic function. The model in~(\ref{eq:55}) is a GLM with $a(y,\phi)=1$, $b(\x\tp\bbeta)=\log\{p(\x\tp\bbeta)\}-\log\{1-p(\x\tp\bbeta)\}$, $c(\x\tp\bbeta)=-\log\{1-p(\x\tp\bbeta)\}$, and $\phi=1$. It is easy to obtain that $\mu=c'(\x\tp\bbeta)/b'(\x\tp\bbeta)=p(\x\tp\bbeta)$. %

For independent full data from the model in (\ref{eq:55}), $(\x_i,y_i)$, $i=1, ..., N$, the subsampling probabilities in~(\ref{eq:3}) reduce to
\begin{equation}\label{eq:12}
  \pi_{N}(\x_i,y_i;\tilde\bvtheta_{\plt})
  =\frac{n\big|\{y_i-p(\x_i\tp\tbeta)\}b'(\x_i\tp\tbeta)\big|h(\x_i)}
  {N\tPsi_{\plt}},
  \quad i=1, ..., N.
\end{equation}
For a subsample taken according to $\pi_{N}(\x_i,y_i;\tilde\bvtheta_{\plt})$ in~(\ref{eq:12}), $\tilde\kappa_{0,i}$, $\tilde\kappa_{1,i}$, and $\tilde\kappa_{3,i}$ have the following specific expressions
\begin{align}
  \tilde\kappa_{0,i}
  &=p(\x_i\tp\bbeta)+p(\x_i\tp\tbeta_{\plt})-2p(\x_i\tp\bbeta)p(\x_i\tp\tbeta_{\plt}),\\ 
  \tilde\kappa_{1,i}
  &=\tilde\kappa_{2,i}=p(\x_i\tp\bbeta)-p(\x_i\tp\bbeta)p(\x_i\tp\tbeta_{\plt}).
\end{align}
We can use these expressions in (\ref{eq:17}) and (\ref{eq:19}), and use the Fisher scoring algorithm to find the MSCLE $\hbeta_{s,\tilde\bvtheta_{\plt}}$, which is the solution to %
\begin{align}
  \dot\ell_{S}(\bbeta\mid\tilde\bvtheta_{\plt})
  &=\sumn\delta_i\bigg[
  {y_i-p(\x_i\tp\bbeta)}
  -\frac{\{1-2p(\x_i\tp\tbeta_{\plt})\}p(\x_i\tp\bbeta)\{1-p(\x_i\tp\bbeta)\}}
    {p(\x_i\tp\bbeta)+p(\x_i\tp\tbeta_{\plt})
    -2p(\x_i\tp\bbeta)p(\x_i\tp\tbeta_{\plt})}\bigg]b'(\x_i\tp\bbeta)\x_i\notag\\
  &=\0.\label{eq:57}
\end{align}

According to Theorem~\ref{thm:1p}, the asymptotic variance-covariance matrix of $\hbeta_{s,\tilde\bvtheta_{\plt}}$ (multiplied by $n$) has an expression of
\begin{align}\label{eq:56}
  \bSigma_{\btheta,\bvtheta_p}=\frac{1}{\Psi_p}\Exp\bigg(
  \frac{p'(\x\tp\bbeta_p)p'(\x\tp\bbeta)b'(\x\tp\bbeta)h(\x)\x^{\otimes2}}
    {\big[p(\x\tp\bbeta_p)\{1-p(\x\tp\bbeta)\}
    +p(\x\tp\bbeta)\{1-p(\x\tp\bbeta_p)\}\big]}\bigg),
\end{align}
where
\begin{equation}
  \Psi_p=\Exp\big[\big|\{y-p(\x\tp\bbeta_p)\}b'(\x\tp\bbeta_p)\big|h(\x)\big].
\end{equation}
If the pilot estimates are consistent, then the expression in~(\ref{eq:56}) simplifies to
\begin{equation}\label{eq:62}
  \bSigma_{\btheta,\bvtheta_p}
  =\frac{\Exp\big[p'(\x\tp\bbeta)\{b'(\x\tp\bbeta)\}^2h(\x)\x^{\otimes2}\big]}
  {4\Exp\{|p'(\x\tp\bbeta)|h(\x)\}}.
  \end{equation}

A special case of the model in (\ref{eq:55}) is the widely used binary logistic regression model, 
\begin{equation}
  \Pr(y=1\mid\x;\bbeta)
  = p(\x\tp\bbeta)
  =\frac{\exp(\x\tp\bbeta)}{1+\exp(\x\tp\bbeta)},
\end{equation}
for which $b(\x\tp\bbeta)=\x\tp\bbeta$, $c(\x\tp\bbeta)=\log\{1+\exp(\x\tp\bbeta)\}$, and $b'(\x\tp\bbeta)=1$. 
For this model, %
if $h(\x_i)=1$ and $\tPsi_{\plt}=n^{-1}N$, then $\pi_{N}(\x_i,y_i;\tbeta_{\plt})$ in (\ref{eq:12}) become the local case-control subsampling probabilities \citep{fithian2014local}; letting $\tilde{\F}$ be a pilot estimate of the Fisher information matrix $\Exp[p(\x,\bbeta)\{1-p(\x,\bbeta)\}\x\x\tp]$, 
if $h(\x_i)=\|\x_i\|$ or $\|\tilde{\F}^{-1}\x_i\|$, then
$\pi_{N}(\x_i,y_i;\tbeta_{\plt})$ in (\ref{eq:12}) become the L-optimal or
A-optimal subsampling probabilities for the IPW estimator \citep{WangZhuMa2018,wang2019more}.

By direct calculation, %
the sampled data conditional log-likelihood score equation in~(\ref{eq:57}) for logistic regression is simplified to
\begin{equation}\label{eq:13}
  \dot\ell_{S}(\bbeta\mid\tilde\bvtheta_{\plt})
  =\sumN\delta_i\big[y_i-p\{\x_i\tp(\bbeta-\tbeta_{\plt})\}\big]\x_i=\0.
\end{equation}
\cite{fithian2014local} and \cite{wang2019more} proposed to solve 
\begin{equation}
  \sumN\delta_i\big\{y_i-p(\x_i\tp\bbeta)\}\x_i=\0,
\end{equation}
and then add $\tbeta_{\plt}$ to the resulting estimator to correct the bias, which is identical to the solution to (\ref{eq:13}). This indicates that the method in \cite{fithian2014local} and \cite{wang2019more} is actually the MSCLE, and the score equation in~(\ref{eq:13}) explains where the magic bias correction term in their method comes from.

Specific for logistic regression, the $\bSigma_{\btheta,\bvtheta_p}$ in~(\ref{eq:56}) %
simplifies to
\begin{align}\label{eq:63}
    \bSigma_{\btheta,\bvtheta_p}
    &=\frac{\Exp\big[
      p(\x\tp\bbeta_p)\{1-p(\x\tp\bbeta)\}
      p\{\x\tp(\bbeta-\bbeta_p)\}h(\x)\x^{\otimes2}\big]}
      {\Exp\big[\{y-p(\x\tp\bbeta_p)\}h(\x)\big]},
\end{align}
which is the same as the result in Theorem 21 of \cite{wang2019more}. If the pilot estimate is consistent so that $\bbeta_p=\bbeta$, then the above expression simplifies to 
\begin{align}\label{eq:64}
    \bSigma_{\btheta,\bvtheta_p}
    &=\frac{\Exp\big[
      p(\x\tp\bbeta)\{1-p(\x\tp\bbeta)\}h(\x)\x^{\otimes2}\big]}
      {4\Exp\big[p(\x\tp\bbeta)\{1-p(\x\tp\bbeta)\}h(\x)\big]}. 
\end{align}
When $h(\x)=1$, we see that 
\begin{align}\label{eq:58}
    \bSigma_{\btheta,\bvtheta_p}
    &=\frac{\F}{4\Exp\big[p(\x\tp\bbeta)\{1-p(\x\tp\bbeta)\}\big]}>\F,
\end{align}
if $\bbeta\neq\0$, indicating that the per-observation information matrix of a local case-control subsample is larger than that of a uniform subsample. 

There is no general relationship between the
efficiency of the MSCLE and that of the full data MLE, except that the full data MLE
has a higher estimation efficiency than a subsample estimator.
  The interesting result that the asymptotic variance of the local case-control
  subsampling estimator is twice of that of the full data MLE \citep{fithian2014local} lies
  in the very unique structure of the logistic regression model and the special
  form of the local case-control inclusion probability. The MSCLE includes the local case-control subsampling as a special
case, so we can obtain some insights on this interesting result. Since a larger sample size
  results in a smaller variance in general, it is important to pay attention to
  the average sample size for a subsampling method. For the local case-control subsampling
  probability with a consistent pilot under correctly specified logistic
  regression model, the average sample size is asymptotically
  $n=N\Exp\{|y-p(\x\tp\bbeta)|\}=2N\Exp[p(\x\tp\bbeta)\{1-p(\x\tp\bbeta)\}]$ and
  the per-observation information matrix is presented as $\bSigma_{\btheta,\bvtheta_p}$ in \eqref{eq:58}.
  Thus the asymptotic variance of the local case-control subsampling is
  $n^{-1}\bSigma_{\btheta,\btheta}^{-1}=2N^{-1}\F$, while the asymptotic
  variance of the full data MLE is $N^{-1}\F$.

  By separating the average sample size and the per-observation information
  matrix, we see that the local case-control subsampling is more appealing for more imbalanced data. For
  the extreme case that $p(\x\tp\bbeta)=0.5$ almost surely, the local case-control subsampling estimator uses half of the
  full data on average to achieve the twice variance and it is just the uniform
  subsampling. However for highly imbalanced data, $n$ can be much smaller than
  $0.5N$ and the per-observation matrix is much larger than $\F$. 

  One key for the interesting result in logistic regression with local case-control subsampling is that
  elements of $\bSigma_{\btheta,\btheta}$ and $\F$ are 
  proportional to each other. Unfortunately, if the model is not the logistic
  regression, if the sampling probability is not the local case-control subsampling probability, or if the
  pilot is not consistent, then this proportional relationship does not hold in
  general, as seen in \eqref{eq:56}, \eqref{eq:62}, \eqref{eq:63}, and
  \eqref{eq:64} for different cases with binary response models.

  The local case-control subsampling sampling is also a quite unique example that the IPW estimator is
  asymptotically as efficient as the MSCLE as noticed in
  \cite{shen2021surprise} and \cite{wang2019more}, but again if the model or the
  sampling probability is changed, then the result does not hold anymore.

\subsubsection{Example: multi-class logistic regressions}
\label{sec:example:-multi-class}
Now we discuss the multi-class logistic regression, which assume that an experiment has $K$ possible outcomes. This is a GLM with a multinomial distribution, a distribution in the multivariate exponential family. 
For $k=1, ..., K$, let $y_{i,k}=1$ if the $k$-th outcomes occurs in the $i$-th experiment and $y_{i,k}=0$ otherwise. Thus $\sumK y_{i,k}=1$. %
The multinomial logistic regression assumes that
\begin{equation}\label{eq:16}
  \Pr(y_{i,k}=1\mid\x_i;\bbeta)
  = p_k(\x_i,\bbeta)
  =\frac{\exp(\x_i\tp\bbeta_k)}{\sum_{l=1}^K\exp(\x_i\tp\bbeta_l)},
  \quad k = 1,2,..., K,
\end{equation}
where $\bbeta=(\bbeta_1\tp, ..., \bbeta_K\tp)\tp$ is the regression coefficient vector. 
Let $\y_i=(y_{i,1}, ..., y_{i,K})\tp$. For an observation $(\x_i,\y_i)$, the density of the multivariate response $\y_i$ at $\x_i$ is 
\begin{equation}\label{eq:70}
  f(\y_i\mid\x_i;\bbeta)
  =\frac{\exp(\sumK y_{i,k}\x_i\tp \bbeta_k)}
  {\sum_{l=1}^K\exp(\x_i\tp \bbeta_l)}%
  =\frac{\exp\{\bbeta\tp(\y_i\otimes\x_i)\}}
  {\sum_{l=1}^K\exp(\x_i\tp \bbeta_l)},%
  \quad y_{i,k}=0,1,
\end{equation}
where $\otimes$ is the kronecker product.
The full data log-likelihood is
\begin{equation}
  \ell(\bbeta)
  =\sumN\bigg[(\y_i\otimes\x_i)\tp\bbeta - \log\bigg\{\sum_{l=1}^K\exp(\x_i\tp \bbeta_l)\bigg\}\bigg],
\end{equation}
with the corresponding score function as
\begin{equation}
  \dot\ell(\bbeta)
  =\sumN\{\y_i-\p(\x_i,\bbeta)\}\otimes\x_i,
\end{equation}
where $\p(\x_i,\bbeta)=\{p_{1}(\x_i,\bbeta), ..., p_{K}(\x_i,\bbeta)\}\tp$. Note that for this model $\sumK p_{k}(\x_i,\bbeta)=1$, so not all $\bbeta_k$'s are estimable and a common constrain is to assume that the regression coefficient corresponding to a baseline class is $\0$, e.g., $\bbeta_K=\0$. Thus the full vector of unknown regression parameter is $\bbeta_{-K}=(\bbeta_1\tp, \bbeta_2\tp, ..., \bbeta_{K-1}\tp)\tp$, whose dimension is $(K-1)d$.  %

The approximate L-optimal subsampling probabilities for the IPW estimator \citep{yao2019optimal} are
\begin{equation}\label{eq:15}
  \pi_{N}(\x_i,\y_i;\tilde\bvtheta_{\plt})
  =\frac{n\|\y_i-\p(\x_i,\tbeta_{\plt})\|\|\x_i\|}
  {N\tPsi_{\plt}}, \quad i=1, ..., N,
\end{equation}
where $\tPsi_{\plt}$ is a pilot estimate of $\Exp\{\|\y-\p(\x,\bbeta)\|\|\x\|\}$. 
For the sampled data, the conditional log-likelihood score equation is (detailed derivations in Section~\ref{sec:deriv-equat-refeq:14}).
\begin{equation}\label{eq:14}
  \dot\ell_{S}(\bbeta\mid\tilde\bvtheta_{\plt})
  =\sumN\delta_i\big\{\y_i-\p_i^g(\bbeta,\tbeta_{\plt})\big\}\otimes\x_i=\0,
\end{equation}
where
\begin{equation}
  \p_i^g(\bbeta,\tbeta_{\plt})=
  \begin{pmatrix}
    \tilde{p}_{i,1}^g\\
    \vdots\\
    \tilde{p}_{i,K}^g
  \end{pmatrix},\quad
  \tilde{p}_{i,k}^g=\frac{\exp(\x_i\tp\bbeta_k+\tilde{g}_{i,k})}
  {\sum_{l=1}^K\exp(\x_i\tp\bbeta_l+\tilde{g}_{i,l})},\quad
  \tilde{g}_{i,k}=\log\{\pi_{N}(\x_i,\1_k;\tilde\bvtheta_{\plt})\},
\end{equation}
and $\1_k$ is the $K$-dimensional unit vector with the $k$-th element being one
and other elements being zero. Here, the expression of
$\dot\ell_{S}(\bbeta\mid\tilde\bvtheta_{\plt})$ in (\ref{eq:14}) holds in
general and it is not restricted to the sampling probability
in~(\ref{eq:15}). For the sampling probability in~(\ref{eq:15}), the
$\tilde{g}_{i,k}$ in (\ref{eq:14}) has a specific form of
$\tilde{g}_{i,k}=0.5\log\big\{\sum_{l=1}^Kp_{l}^2(\x_i,\tbeta_{\plt})
+1-2p_{k}(\x_i,\tbeta_{\plt})\big\}$.

Note that $\tilde{g}_{i,k}$'s do not contain $\bbeta$; they depend only on the sample data and the pilot $\tbeta_{\plt}$. Thus, solving (\ref{eq:14}) is as easy as solving the score function without correction. 
When we use the constrain that $\bbeta_K=\0$, we only need to solve the first $(K-1)d$ components of (\ref{eq:14}) and the last $d$ equations are automatically satisfied. For the case of $K=2$, the first $d$ equations are the same as the last $d$ equations. In this special case, $\tilde{g}_{i,k}=-\x_i\tp\tbeta_{\plt}$, and the score equation in (\ref{eq:14}) reduce to  that in (\ref{eq:13}). 

The asymptotic variance-covariance matrix (multiplied by $n$) for the MSCLE of $\bbeta_{-K}$ with data taken according to (\ref{eq:15}) is
\begin{align}\label{eq:72}
  \bSigma_{\btheta,\bvtheta_p}
  &=\frac{\Exp\big[\|\y-\p(\x,\bbeta_p)\|\|\x\|
    \{\y_{-K} -\p_{-K}^g(\bbeta,\bbeta_p)\}^{\otimes2}
    \otimes\x^{\otimes2}\big]}{\Exp\{\|\y-\p(\x,\bbeta_p)\|\|\x\|\}},
\end{align}
where $(\x,\y)$ is an observation from the data distribution, $\y_{-K}$ is $\y$
with the $K$-th element removed, $\p_{-K}^g(\bbeta,\bbeta_p)$ is
$\p^g(\bbeta,\bbeta_p)$ with the $K$-th element removed, and
$\p^g(\bbeta,\bbeta_p)$ has the same expression as $\p^g(\bbeta,\tbeta_{\plt})$
except that $\tbeta_{\plt}$ is replaced by $\bbeta_p$.

Since the expression of
$\dot\ell_{S}(\bbeta\mid\tilde\bvtheta_{\plt})$ in (\ref{eq:14}) holds in
general, it holds for the LUS probability as well, so the LUS estimator
is a specific case of the MSCLE. The LUS probability is designed
for imbalanced multi-class response models, and it has an expression of
\begin{align*}
\pi_{N}^{\lus}(\x,\y;\tilde\bvtheta_{\plt})
    =&2I\{\gamma\ge2q(\tbeta_{\plt})\}
      \frac{q(\tbeta_{\plt})+\{1-2q(\tbeta_{\plt})\}\eta(\tbeta_{\plt})}{\gamma}\\
    &+ I\{\gamma<2q(\tbeta_{\plt})\}
      \frac{\{\gamma-q(\tbeta_{\plt})\}+(1-\gamma)\eta(\tbeta_{\plt})}{\gamma-q(\tbeta_{\plt})},
\end{align*}
where $q(\tbeta_{\plt})=\max\{0.5, p_1(\x,\tbeta_{\plt}), ..., p_K(\x,\tbeta_{\plt})\}$,
$\eta(\tbeta_{\plt})=I\{\y\tp\p(\x,\tbeta_{\plt})=q(\tbeta_{\plt})\}$ with $I()$ being the indicator function,
and $\gamma$ is a tuning parameter that corresponds to an upper bound of the
average subsample size, namely, $n\ge N/\gamma$.

Our Theorem~\ref{thm:1p} indicates that the inverse of the asymptotic variance
of the LUS estimator is proportional to
\begin{align}\label{eq:71}
    \sumK\Exp\Big[p_k(\x,\btheta)
    \pi^{\lus}(\x,\1_k;\bbeta_p)\{\1_k-\p^g(\bbeta,\bbeta_{p})\}^{\otimes2}
    \otimes\x^{\otimes2} \Big],
\end{align}
where $\p^g(\bbeta,\bbeta_{p})$ has the same expression as
$\p_i^g(\bbeta,\tbeta_{\plt})$ in (\ref{eq:14}) except that $\tbeta_{\plt}$ and
$\x_i$ are replaced by $\bbeta_{p}$ and $\x$, respectively.
Our results allow $n=o(N)$ which corresponds to the scenario that
$\gamma\rightarrow\infty$. In this case, the LUS probability reduces to satisfy
that $\pi_{N}^{\lus}(\x,\y;\tilde\bvtheta_{\plt}) \propto
q(\tbeta_{\plt})+\{1-2q(\tbeta_{\plt})\}\eta(\tbeta_{\plt})$, and (\ref{eq:71}) still holds with
$\pi^{\lus}(\x,\1_k;\bbeta_p)$ replaced by $q(\bbeta_p)+\{1-2q(\bbeta_p)\}I\{\1_k\tp\p(\x,\bbeta_p)=q(\bbeta_p)\}$. 
If the pilot is
consistent, i.e., $\bbeta_{p}=\bbeta$, and $\gamma$ is fixed and finite, then the expression in~(\ref{eq:71})
reduces to be proportional to the inverse of $\mathcal{V}_{\lus}$
in Corollary 4.1 of \cite{Han2019}.

Note that the parameter $\gamma$ only controls an upper bound of the average
subsample size, so the expressions in (\ref{eq:72}) and (\ref{eq:71}) are not
comparable because they are scaled by different average subsample sizes. We will
provide numerical comparisons by using the save average subsample sizes in
Section~\ref{sec:multi-class-logistic}. We will also implement the IPW
estimator to show that the MSCLE has a higher estimation efficiency than the
IPW estimator for the LUS probability as well.

\subsubsection{Example: Poisson regression}
\label{sec:poisson-regression}
Poisson regression models are commonly used for modeling count data. It assumes that given the covariate $\x_i$, the response $y_i$ follows a Poisson distribution with density
\begin{equation}
  f(y_i\mid\x_i;\bbeta)
  =\frac{e^{-\mu_i}\mu_i^{y_i}}{y_i!}
  =\frac{1}{y_i!}\exp(y_i\x_i\tp\bbeta-e^{\x_i\tp\bbeta}), \quad y_i= 0, 1, ...,
\end{equation}
where $\mu_i=\exp(\x_i\tp\bbeta)$. This is a specific case of (\ref{eq:10}) with $\phi=1$, $a(y_i,\phi)=(y_i!)^{-1}$, $b(\x_i\tp\bbeta)=\x_i\tp\bbeta$ and $c(\x_i\tp\bbeta)=\exp(\x_i\tp\bbeta)$. The link function in Poisson regression is $g(\mu_i)=\log(\mu_i)$ and it is the canonical link.

For this model, the subsampling probabilities in~(\ref{eq:3}) %
reduce to
\begin{equation}
  \pi_{N}(\x_i,y_i;\tilde\bvtheta_{\plt})
  =\frac{n|y_i-\tmu_i|h(\x_i)}
  {N\tPsi_{\plt}}, \quad i=1, ..., N,
\end{equation}
where $\tmu_i=\exp(\x_i\tp\tbeta_{\plt})$. For a subsample taken according to $\pi_{N}(\x_i,y_i;\tilde\bvtheta_{\plt})$, the sampled data conditional log-likelihood has the score equation and the Hessian matrix as %
\begin{equation}
  \dot\ell_{S}(\bbeta\mid\tilde\bvtheta_{\plt})
  =\sumN\delta_i\Big(y_i-\frac{\tilde\kappa_{1,i}}{\tilde\kappa_{0,i}}\Big)\x_i,
\end{equation}
and
\begin{equation}
  \ddot\ell_{S}(\bbeta\mid\tilde\bvtheta_{\plt})
  =-\sumN\delta_i\Big(\frac{\tilde\kappa_{2,i}}{\tilde\kappa_{0,i}}
  -\frac{\tilde\kappa_{1,i}^2}{\tilde\kappa_{0,i}^2}\Big)\x_i\x_i\tp,
\end{equation}
respectively.

Now we show the closed-form expressions of $\tilde\kappa_{0,i}$, $\tilde\kappa_{1,i}$, and $\tilde\kappa_{2,i}$. 
Let $m_i=\lfloor{\tmu_i}\rfloor$ be the largest integer that is smaller than or equal to $\tmu_i$, and let $F(\cdot;\mu_i)$ be the cumulative distribution function for a Poisson distribution with mean $\mu_i$. Here we show the results and give the detailed derivations in Section~\ref{sec:deriv-equat-refeq:20}. We have that
\begin{align}
  \tilde\kappa_{0,i}
  &=2\tmu_iF(m_i;\mu_i)-2\mu_iF(m_i-1;\mu_i)+\mu_i-\tmu_i,\label{eq:20}\\
  \tilde\kappa_{1,i}
  &=2\mu_i(\tmu_i-1)F(m_i-1;\mu_i)-2\mu_i^2F(m_i-2;\mu_i)
    +\mu_i+\mu_i^2-\mu_i\tmu_i,\label{eq:21}\\
  \tilde\kappa_{2,i}
  &=\tilde\kappa_{1,i}+\mu_i^2(\tmu_i-2)\{2F(m_i-2;\mu_i)-1\}
    -2\mu_i^3F(m_i-3;\mu_i) +\mu_i^3.\label{eq:22}
\end{align}
With the above expressions, we can find the MSCLE by applying the Newton's algorithm. We will demonstrate the performance of the resulting estimator in Section~\ref{sec:numer-exper}

\section{Model misspecification}
\label{sec:model-missp}

When the assumed working model is misspecified, the meaning of consistency of an
estimator needs to be carefully defined. One definition of the  consistency of a subsample
estimator is that it has the same probability limit as the full data estimator
\citep{fithian2014local,Han2019,shen2021surprise}. In this definition, the
IPW estimator $\htheta_{W}$ defined in~(\ref{eq:4}) is consistent, because the
conditional expectation of the selection indicator $\delta$ give the data
equals the inclusion probability and therefore the expectation of the objective
function for the IPW estimator is the same as that of the full data
MLE. Specifically, the full data MLE under a misspecified model estimates the
$\btheta_l$ that solves the following population estimation equation
\begin{equation}\label{eq:65}
  \Exp^t\{\dot\ell(\btheta_l;\x,y)\}=\0,
\end{equation}
where we use $\Exp^t$ to emphasize that the expectation is taken with respect to
the true data distribution. Since
$\Exp\{\delta\pi^{-1}(\x,y)\dot\ell(\btheta;\x,y)]=\Exp^t\{\dot\ell(\btheta;\x,y)]$,
the IPW estimator $\htheta_{W}$ defined in~(\ref{eq:4}) always estimates
the same $\btheta_l$ as the full data MLE. Actually, the
subsample IPW estimator $\htheta_{W}$ is often investigated as an estimator
of the full data estimator in the literature \citep[e.g.][]{WangZhuMa2018,yu2020quasi}.

For the MSCLE, we use the assumed model structure to derive the subsampled
conditional likelihood, so if the subsampling probabilities also depend on the
data then the MSCLE is not consistent to the same limit as the full data MLE
when the assumed model structure is misspecified. The MSCLE is consistent to the
$\btheta_{ls}$ that solves the following population estimation equation if the
assumed model is misspecified,
\begin{equation}
\label{eq:66}
  \Exp^t\bigg[\Exp^t\{\dot\ell(\btheta_{ls};\x,y)\pi(\x,y;\bvtheta)\mid\x\}
  - \frac{\Exp^t\{\pi(\x,y;\bvtheta)\mid\x\}}
  {\Exp^w\{\pi(\x,y;\bvtheta)\mid\x\}}
  \Exp^w\{\dot\ell(\btheta_{ls};\x,y)\pi(\x,y;\bvtheta)\mid\x\}\bigg]=\0,
\end{equation}
where $\Exp^w$ is the expectation under the assumed working model and $\bvtheta$
the fixed parameter in calculating the sampling probabilities. Clearly,
$\btheta_l$ and $\btheta_{ls}$ are different in general, although they can be
the same under the special case of the local case-control subsampling.

For an illustration,  consider the problem under the more specific class of binary response
models discussed in Section~\ref{sec:exampl-binary-resp}. For this class of
models, (\ref{eq:65}) and (\ref{eq:66}) are simplified as
\begin{equation}\label{eq:67}
  \Exp_{\x}\big[\{p_t(\x)-p(\x\tp\btheta_l)\}
  b'(\x\tp\btheta_l)\x\big]=\0,
\end{equation}
and
\begin{equation}\label{eq:68}
  \Exp_{\x}\bigg[\{p_t(\x)-p(\x\tp\btheta_{ls})\}
  \frac{\pi(\x,0;\bvtheta)\pi(\x,1;\bvtheta)}
  {p(\x\tp\btheta_{ls})\pi(\x,1;\bvtheta)+\{1-p(\x\tp\btheta_{ls})\}
    \pi(\x,0;\bvtheta)}b'(\x\tp\btheta_{ls})\x\bigg]=\0,
\end{equation}
respectively, where $\Exp_{\x}$ is the expectation respect to the distribution
of $\x$ and $p_t(\x)$ is the true probability of $y=1$ given $\x$.
If the subsampling probabilities do not depend on the data, i.e.,
$\pi(\x,y;\bvtheta)$ is a consistent across the data distribution, then
(\ref{eq:67}) and (\ref{eq:68}) are identical, meaning that the MSCLE is
consistent to the full data MLE.
In general for a
correctly specified model, $p_t(\x)=p(\x\tp\btheta_t)$ at the true parameter $\btheta_t$ and thus (\ref{eq:67})
and (\ref{eq:68}) are both true when $\btheta_l=\btheta_{ls}=\btheta_t$. For a
misspecified model, heuristically, the $\btheta_l$ in (\ref{eq:67}) tries to
make $p(\x\tp\btheta_l)$ to be close to $p_t(\x)$ and thus it is a reasonable
value to consider. On the other hand, the
$\btheta_{ls}$ in (\ref{eq:68}) also tries to make $p(\x\tp\btheta_{ls})$ to be
close to $p_t(\x)$, and thus it is also a reasonable value for the parameter
when the model is misspecified. It is difficult to argue which of $\btheta_l$
and $\btheta_{ls}$ is better, because the true $p_t(\x)$ is unknown and the
definition of the closeness of $p(\x\tp\btheta)$ to $p_t(\x)$ varies. For
example, the $\btheta$'s that minimize
$\Exp_{\x}[\{p_t(\x)-p(\x\tp\btheta)\}^2]$ and
$\Exp_{\x}\{|p_t(\x)-p(\x\tp\btheta)|\}$ are typically different although they
are both reasonable objectives, and there is no easy way to tell whether $\btheta_l$
or $\btheta_{ls}$ is closer to any of them.

Specifically, if the sampling probabilities in (\ref{eq:12}) are used for
logistic regression ($b'(\cdot)=1$) and assume
that the pilot $\bvtheta=\btheta_{ls}$, then (\ref{eq:68}) simplifies to
\begin{equation}\label{eq:69}
  \Exp_{\x}\big[\{p_t(\x)-p(\x\tp\btheta_{ls})\} h(\x)\x\big]=\0. 
\end{equation}
We see that $\btheta_{ls}$ in (\ref{eq:69}) is in general different from $\btheta_l$ unless
$h(\x)=1$. Again, the $\btheta_{ls}$ in the above population estimation equation
also seems to be a reasonable target if the model is misspecified. Existing
investigations on optimal subsampling focus on defining the $h(\x)$ to improve
the estimation efficiency under a correctly specified model. %
Finding $h(\x)$ to improve the subsample performance under model
misspecifications is also interesting and it requires future investigations.

\section{Numerical experiments}
\label{sec:numer-exper}

\subsection{Multi-class logistic regression}
\label{sec:multi-class-logistic}
We demonstrate the performance of the proposed MSCLE estimator using the multi-class logistic regression model discussed in Section~\ref{sec:example:-multi-class}. We set the full data sample size $N=10^6$, and let the subsample sizes be $n=500; 1000; 1500;$ and $2000$. We assume that the responses have three possible categories ($K=3$), and let the dimension of the covariates $\x_i=(1,\x_{-1,i}\tp)\tp$'s be $d=4$ where the first element of one is for the intercept parameters. For this setup, the dimension of the unknown regression coefficient vector is eight, and we set the true parameter to be $\bbeta_{-K}=(0.05, 0.1, 0.15, 0.2, 0.25, 0.3, 0.35, 0.4)\tp$ when generating the data. To generate the covariates corresponding to the slope parameters, $\x_{-1,i}$'s, we consider the following distributions. 
\begin{enumerate}[(a)]
\item Multivariate normal distribution
  $\x_{-1,i}\sim\Nor(\bm{0}, \bm\Omega)$, where the $(i,j)$-th element
  of $\bm\Omega$ is $\Sigma_{ij}=0.5^{|i-j|}$ and $I()$ is the indicator
  function. This distribution is symmetric with light tails. The
  resulting probabilities of the responses for the three possible
  categories are 0.3, 0.39, and 0.31, respectively.
\item Multivariate log-normal %
  distribution
  $\x_{-1,i}\sim\mathbb{LN}(\bm{0}, \bm\Omega)$, where
  $\x_{-1,i}=\exp(\z_i)$, $\z_i\sim\Nor(\bm{0}, \bm\Omega)$, and $\bm\Omega$
  is the same as defined in the above case a). This distribution is
  asymmetric and positively skewed. The resulting probabilities of the
  responses for the three possible categories are 0.22, 0.65, and
  0.13, respectively.
\item Multivariate $t$ distribution,
  $\x_{-1,i}\sim\mathbb{T}_3(\bm{0}, \bm\Omega)$, where $\bm\Omega$ is
  defined in case a). This distribution is symmetric, and response
  composition is similar to case a). However, this distribution have
  heavier tails.
\item Independent exponential distribution,
  $\x_{-1,i}\sim\mathbb{EXP}(1)$, where components of $\x_{-1,i}$
  independently follow the standard exponential distribution. This
  distribution is asymmetric and positively skewed. The response
  composition is similar to case b).
\end{enumerate}

\begin{figure}[htp]%
  \centering
  \begin{subfigure}{\textwidth}\centering
    \includegraphics[width=0.45\textwidth,page=1]{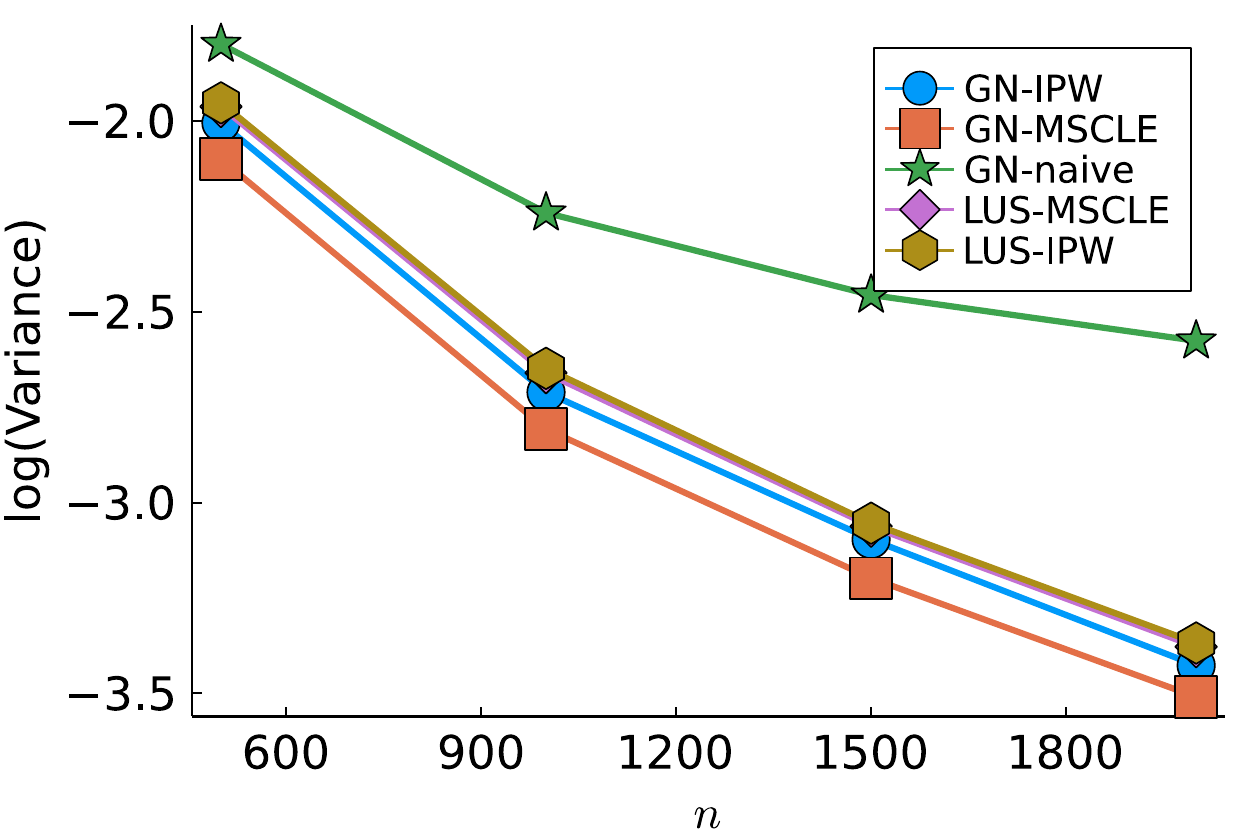}
    \includegraphics[width=0.45\textwidth,page=1]{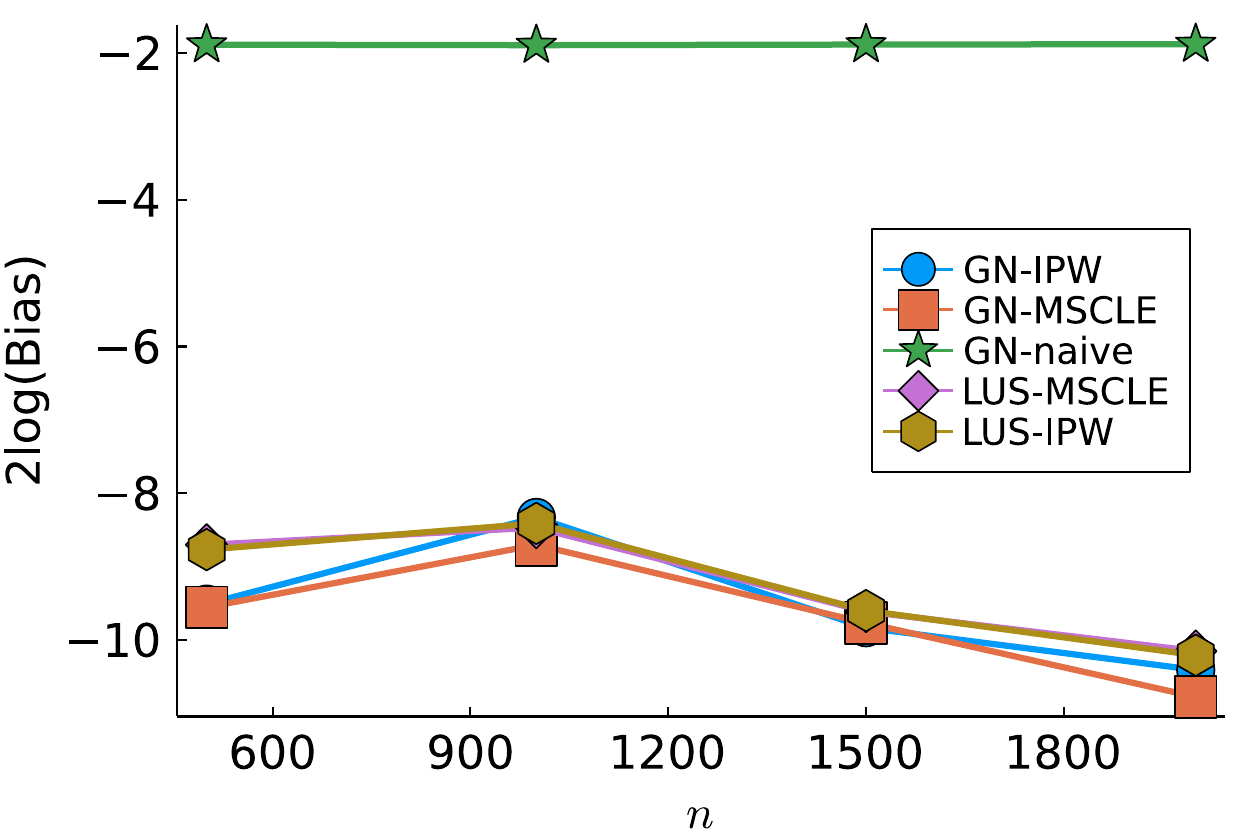}
    \caption{$\x_{-1,i}$'s are normal}
  \end{subfigure}
  \begin{subfigure}{\textwidth}\centering
    \includegraphics[width=0.45\textwidth,page=2]{figures/00mse_vMultiLogi.pdf}
    \includegraphics[width=0.45\textwidth,page=2]{figures/00mse_bMultiLogi.pdf}
    \caption{$\x_{-1,i}$'s are lognormal}
  \end{subfigure}
  \begin{subfigure}{\textwidth}\centering
    \includegraphics[width=0.45\textwidth,page=3]{figures/00mse_vMultiLogi.pdf}
    \includegraphics[width=0.45\textwidth,page=3]{figures/00mse_bMultiLogi.pdf}
    \caption{$\x_{-1,i}$'s are $\mathbb{T}_3$}
  \end{subfigure}
  \begin{subfigure}{\textwidth}\centering
    \includegraphics[width=0.45\textwidth,page=4]{figures/00mse_vMultiLogi.pdf}
    \includegraphics[width=0.45\textwidth,page=4]{figures/00mse_bMultiLogi.pdf}
    \caption{$\x_{-1,i}$'s are exponential}
  \end{subfigure}
  \caption{Log of empirical variances and squared biases (the smaller the better) of subsample estimators for different sample sizes in multi-class logistic regression.}
  \label{fig:logistic1}
\end{figure}

We repeat the simulation for $R=1000$ times and calculate the empirical
mean squared error (MSE) as $R^{-1}\sum_{r=1}^R\|\hbeta_{-K}^{(r)}-\bbeta_{-K}\|^2$,
empirical variance as $R^{-1}\sum_{r=1}^R\|\hbeta_{-K}^{(r)}-\bar\bbeta_{-K}\|^2$ and
empirical squared bias as $\|\bar\bbeta_{-K}-\bbeta_{-K}\|^2$,
where $\hbeta_{-K}^{(r)}$ is the estimate at the $r$-th repetition and
$\bar\bbeta_{-K}$ is the average of $\hbeta_{-K}^{(r)}$'s. 
In each repetition, we generate the full data set and use a uniform sample of
size 400 to calculate the pilot estimates. Results are presented in
Figure~\ref{fig:logistic1}. Since the empirical MSE is the sum of the
  asymptotic variance and empirical squared bias, we only plot the empirical
  variances and squared biases. For comparison, we also implement the commonly
used IPW estimator (IPW) and a naive estimator
(naive) that is obtained by maximizing
$\sumN\delta_i\ell(\btheta;\x_i,y_i)$. We implement both the gradient
  normal (GN) based sampling probability and the LUS probability. For the LUS
  probability, with $\gamma=N/n$, the average subsample size would be smaller
  than $n$ so the comparison would not be fair. Thus we have to set $\gamma$ to
  be smaller than $N/n$ for different distributions of $\x$ in order to
  have a fair comparison. Since the sampling ratios $n/N$ considered here are
  small, this is equivalent to scale the LUS probability by numbers that are
  larger than one.

The simulation results in Figure \ref{fig:logistic1} can be summarized as
follows: 1) The variance is the dominating term in the MSE and the squared
bias is a small term for the MSCLE and the IPW estimators. 2) The bias of the
Naive estimator does not decrease to zero with the increase of the subsample
size, which suggests that the naive estimator is subject to non-negligible
biases. 3) Both the IPW estimator and our proposed MSCLE show that the variance
decreases as the sample size increases.
3) In all scenarios for both the GN probability and the LUS probability, our
proposed MSCLE has smaller MSEs than the IPW estimator, which confirms our
theory in Theorem~\ref{thm:2p}. We see that the GN probability outperforms the
LUS probability especially if the responses are balanced or if the covariate
distributions have heavy tails. This is because the LUS probability is designed
specifically to address imbalanced data and it does not take into account the
structure information in the variance matrix represented by $\x$.

\begin{figure}[htp]%
  \centering
  \begin{subfigure}{\textwidth}\centering
    \includegraphics[width=0.45\textwidth,page=1]{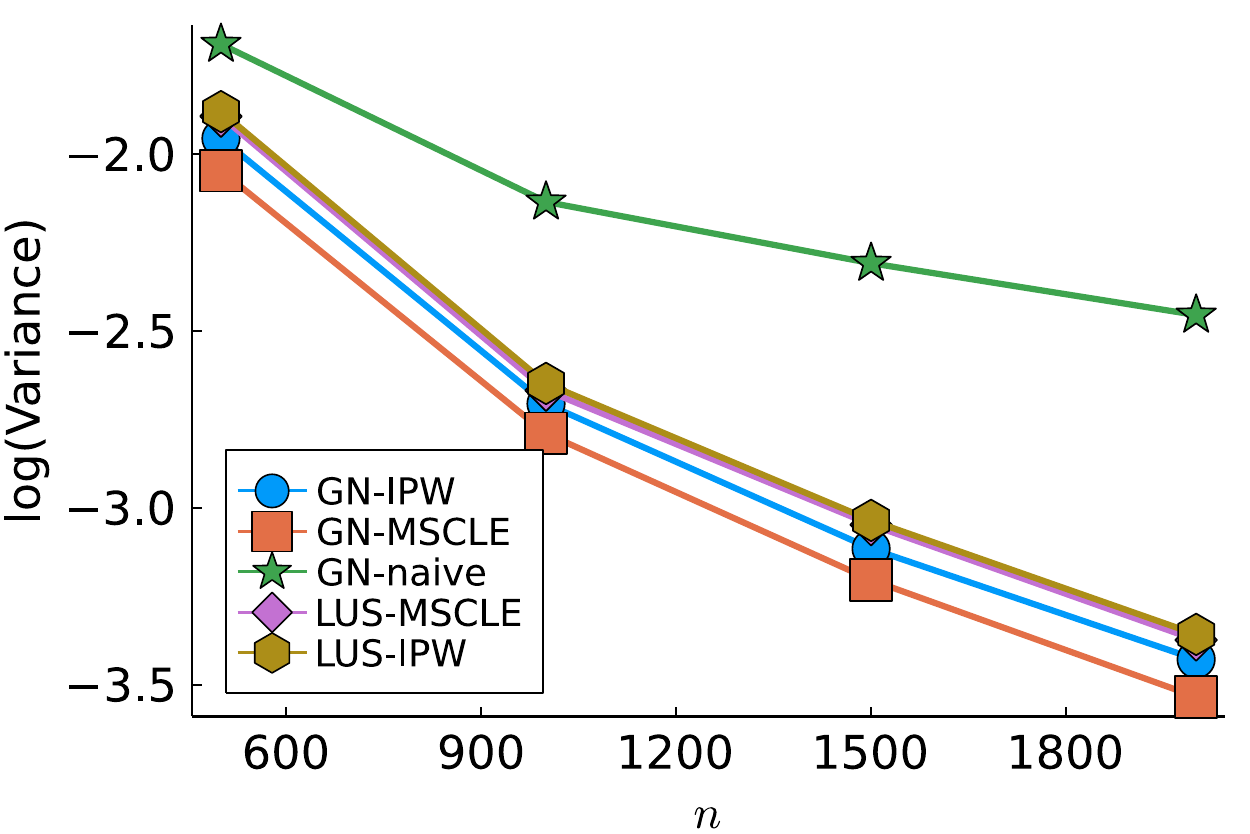}
    \includegraphics[width=0.45\textwidth,page=1]{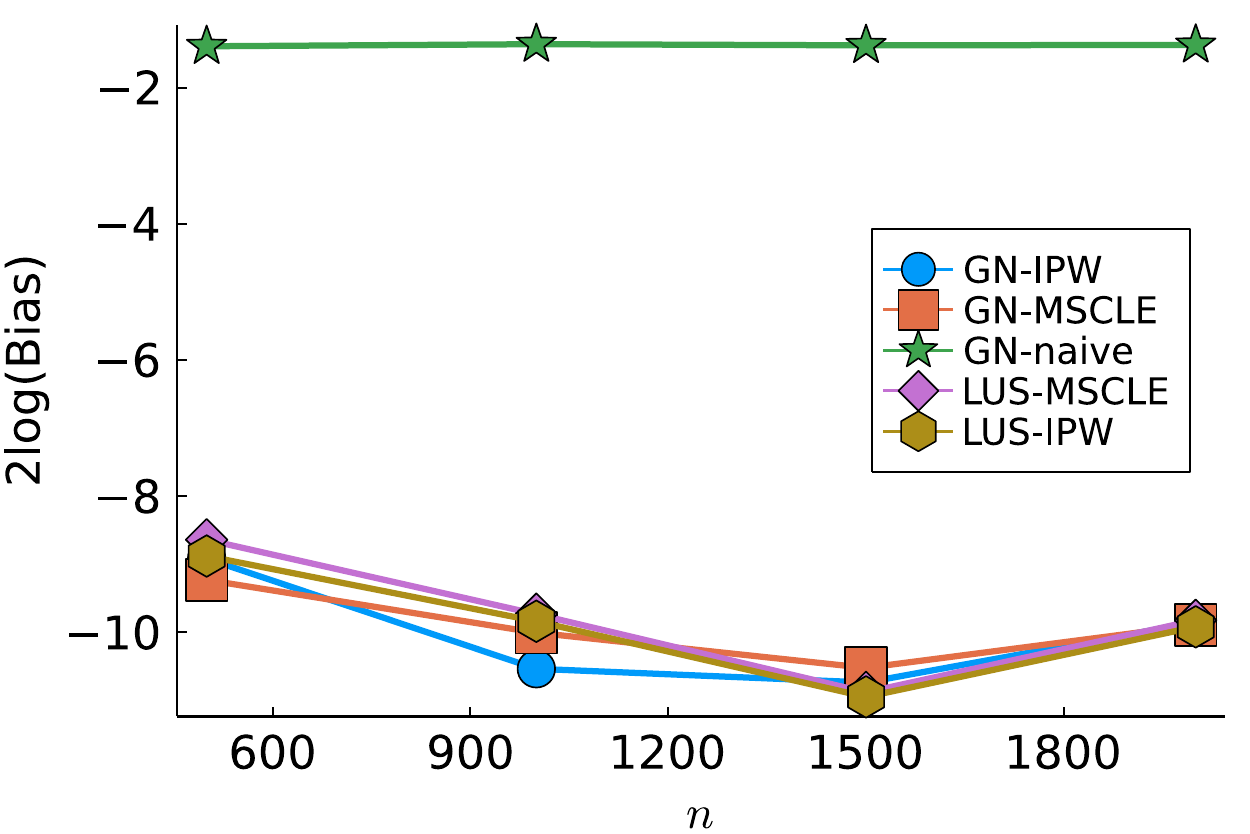}
    \caption{$\x_{-1,i}$'s are normal}
  \end{subfigure}
  \begin{subfigure}{\textwidth}\centering
    \includegraphics[width=0.45\textwidth,page=2]{figures/00mse_vMultiLogiMisPlt.pdf}
    \includegraphics[width=0.45\textwidth,page=2]{figures/00mse_bMultiLogiMisPlt.pdf}
    \caption{$\x_{-1,i}$'s are lognormal}
  \end{subfigure}
  \begin{subfigure}{\textwidth}\centering
    \includegraphics[width=0.45\textwidth,page=3]{figures/00mse_vMultiLogiMisPlt.pdf}
    \includegraphics[width=0.45\textwidth,page=3]{figures/00mse_bMultiLogiMisPlt.pdf}
    \caption{$\x_{-1,i}$'s are $\mathbb{T}_3$}
  \end{subfigure}
  \begin{subfigure}{\textwidth}\centering
    \includegraphics[width=0.45\textwidth,page=4]{figures/00mse_vMultiLogiMisPlt.pdf}
    \includegraphics[width=0.45\textwidth,page=4]{figures/00mse_bMultiLogiMisPlt.pdf}
    \caption{$\x_{-1,i}$'s are exponential}
  \end{subfigure}
  \caption{Log of empirical variances and squared biases (the smaller the better) of subsample estimators for different sample sizes in multi-class logistic regression when the pilot estimator is misspecified.}
  \label{fig:logistic2}
\end{figure}

To examine  the impact of   misspecification in the pilot estimator, we run the aforementioned simulation with the same setting except that we use a wrong pilot estimator. Specifically, for each component of the pilot estimator obtained from the pilot sample, we add to it a random number generated from a uniform distribution between 1 and 2, and then we use this misspecified pilot estimator to implement the all the estimators in comparison. Results are presented in Figure~\ref{fig:logistic2}.

From Figure~\ref{fig:logistic2}, we see that the MSE gets larger for all the methods, especially for the naive estimator and the IPW estimator. This is because with this misspecified pilot estimator the resulting informative subamples are more different from the original data distribution. The inflation on the MSE due to the misspecified pilot estimator is much smaller for the MSCLE than that for the IPW estimator. The advantage of the MSCLE over the IPW estimator become more significant with the misspecified pilot estimator. 
One reason for this behavior is that adding a random bias to the pilot estimator
$\tilde\bvtheta_{\plt}$ brings in additional variation to the estimated probabilities
$\pi_{N}(\x_i,\y_i;\tilde\bvtheta_{\plt})$'s. These probabilities are in
the denominators for the IPW estimator while they are in the expectation
and then log-transformed for the
MSCLE, so the additional variation expresses more significantly for the IPW 
estimator, resulting in a larger variance of the estimator. Another reason is that
the GN sampling probability is a version of the L-optimal probabilities for the
IPW estimator with consistent pilot
\citep{WangZhuMa2018,ai2020optimal,yu2020quasi}, while it is not optimal for the
MSCLE, so a systematic bias on the pilot has a larger negative impact on the
IPW estimator. 
We have additional numerical experiments to investigate the effect of different
levels of pilot misspecification on the final estimator. Please see the details
in Section~\ref{sec:multi-class-logistic-1} of the Appendix.

\begin{figure}[H] %
  \centering
  \begin{subfigure}{0.485\textwidth}
    \includegraphics[width=\textwidth,page=1]{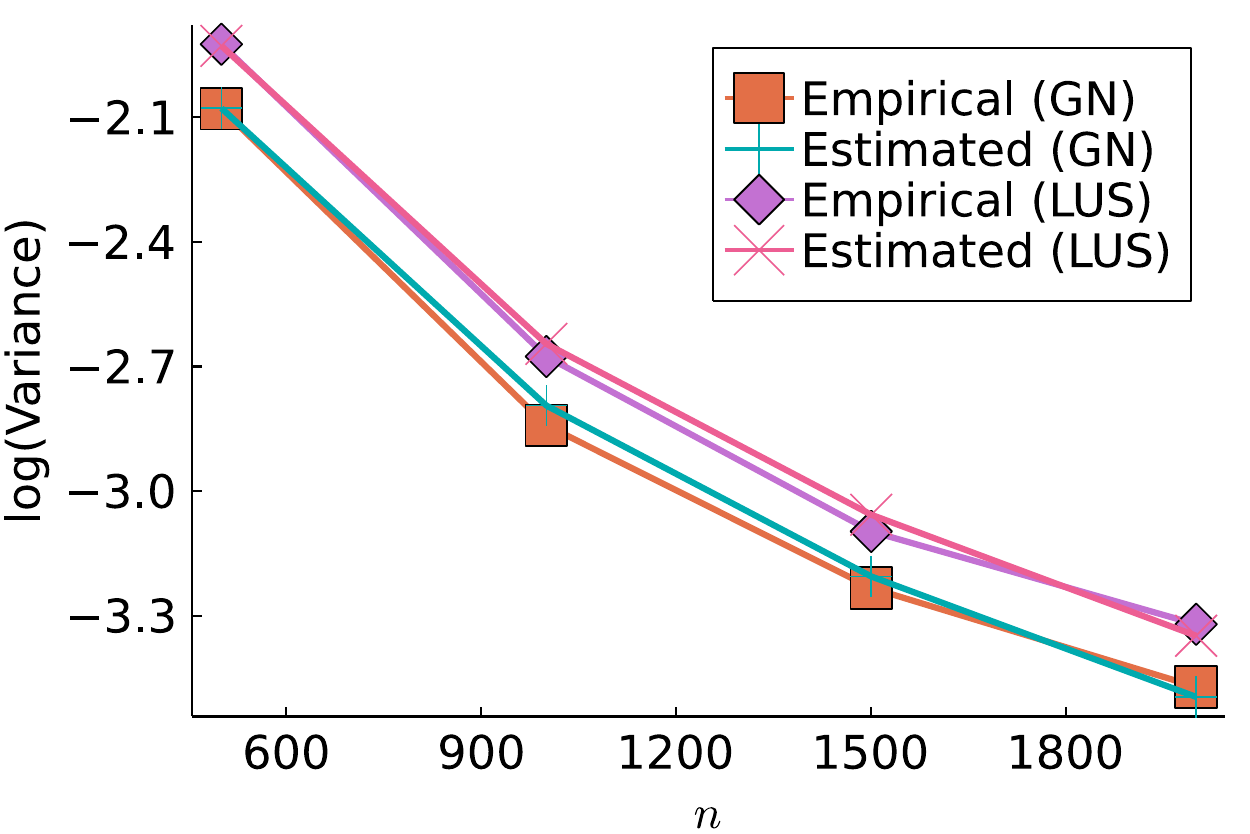}
    \caption{$\x_{-1,i}$'s are normal}
  \end{subfigure}
  \begin{subfigure}{0.485\textwidth}
\includegraphics[width=\textwidth,page=2]{figures/00VarEstMultiLogi.pdf}
    \caption{$\x_{-1,i}$'s are lognormal}
  \end{subfigure}
  \begin{subfigure}{0.485\textwidth}
\includegraphics[width=\textwidth,page=3]{figures/00VarEstMultiLogi.pdf}
    \caption{$\x_{-1,i}$'s are $\mathbb{T}_3$}
  \end{subfigure}
  \begin{subfigure}{0.485\textwidth}
\includegraphics[width=\textwidth,page=4]{figures/00VarEstMultiLogi.pdf}
    \caption{$\x_{-1,i}$'s are exponential}
  \end{subfigure}
  \caption{Log of empirical variances and average estimated variances of
    subsample MSCLE for different sample sizes in multi-class logistic
    regression. Logarithm is taken for better presentation.}
  \label{fig:logistic3}
\end{figure}

To evaluate the performance of the asymptotic variance in Theorem~\ref{thm:1p},
we also calculated the estimated variance by plug-in estimation. The results are
reported in Figure~\ref{fig:logistic3}. The empirical variances
and the estimated variances are quite close, and this is the case for both the
GN sampling probability and the LUS probability, conforming our theoretical
result in Theorem~\ref{thm:1p}.

The MSCLE needs to calculate
$\bar{\pi}_{N}(\x_i;\btheta\mid\tilde\bvtheta_{\plt})$'s and may
require more computational resources. However, this is done only on the selected
subsample, so the actual difference for computational cost will not be
large. %
To check this, we also recorded the computational costs in terms of CPU times
and memory allocations for the MSCLE and the IPW estimator.  As expected, the
time differences would not be noticeable if we present the total time of
calculating a subsample estimator and calculating
$\pi_{N}(\x_i,y_i;\tilde\bvtheta_{\plt})$'s. To see the differences more
clearly, we separated the times for calculating subsample estimators from those
of calculating $\pi_{N}(\x_i,y_i;\tilde\bvtheta_{\plt})$'s. We implemented all
the algorithm in Julia \citep{bezanson2017julia} on a Desktop running Ubuntu
20.04. We restricted all the calculations to use one thread of the CPU with a
base frequency of 2,200 megahertz and a maximum boosted frequency of 4,549
megahertz. We repeated the simulation for 100 times, and calculated the average
CPU times and memory allocations. We used a smaller number of iterations here
because the variations of the computational costs across different repetitions
are much smaller than that of the estimators.

Results for case (a) with multivariate normal covariates are reported in
Table~\ref{tab:2}. Indeed, the MSCLE takes more time than the IPW, but the
difference is very small and is negligible compared with the major time of
calculating $\pi_{N}(\x_i,y_i;\tilde\bvtheta_{\plt})$'s that both methods
require. In terms of the memory allocations, the major cost is also on
calculating $\pi_{N}(\x_i,y_i;\tilde\bvtheta_{\plt})$'s. Here IPW uses a little
more memory because we created an additional weighted covariate matrix in the
implementation to save some CPU time. The difference on memory allocations is
not significant either. Results for other covariate distributions are similar
and thus we omit them.

\begin{table}[htbp]
\caption{Average CPU times and Memory allocations for the multi-class logistic
  regression example. Here ``calPI'' is for the step of calculating sampling
  probabilities $\pi_{N}(\x_i,y_i;\tilde\bvtheta_{\plt})$'s and taking subsamples.}
\label{tab:2}
\centering
\begin{tabular}{l|rrrrc|rrrr}
\hline
  & \multicolumn{4}{c}{CPU time (millisecond)} & 
  & \multicolumn{4}{c}{Memory allocation (megabyte)} \\ 
  & $n$=500  & 1000  & 1500 & 2000  &  & $n=$500  & 1000  & 1500 & 2000   \\\hline
MSCLE    & 0.28  & 0.45  & 0.62  & 0.77  &  & 0.88   & 1.72   & 2.56   & 3.40   \\
IPW & 0.27  & 0.41  & 0.54  & 0.71  &  & 0.89   & 1.78   & 2.64   & 3.53   \\
calPI    & 41.17 & 39.18 & 38.99 & 38.19 &  & 160.38 & 160.41 & 160.44 & 160.47 \\ 
Full & \multicolumn{4}{c}{660.40} &  & \multicolumn{4}{c}{1,777.66}  \\ \hline
\end{tabular}
\label{}
\end{table}

\begin{figure}[htp]%
  \centering
  \begin{subfigure}{\textwidth}\centering
    \includegraphics[width=0.45\textwidth,page=1]{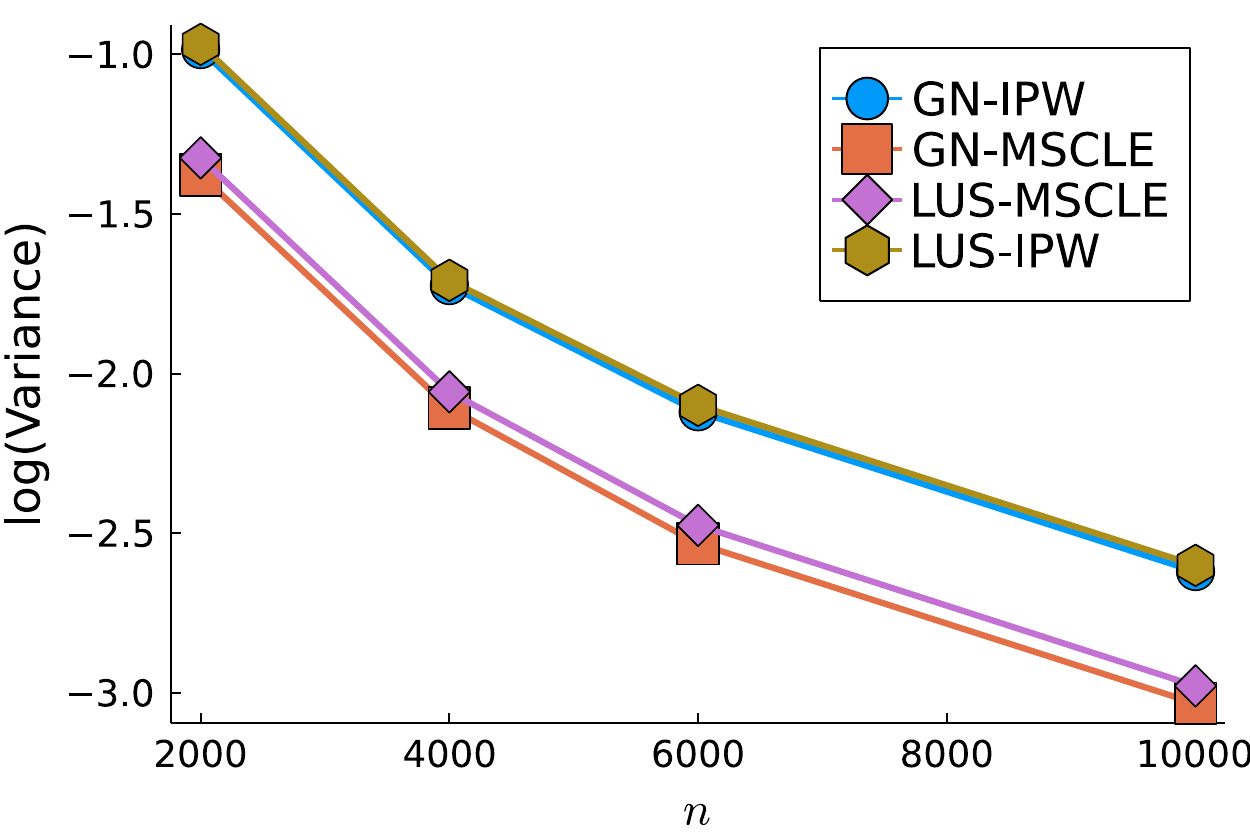}
    \includegraphics[width=0.45\textwidth,page=1]{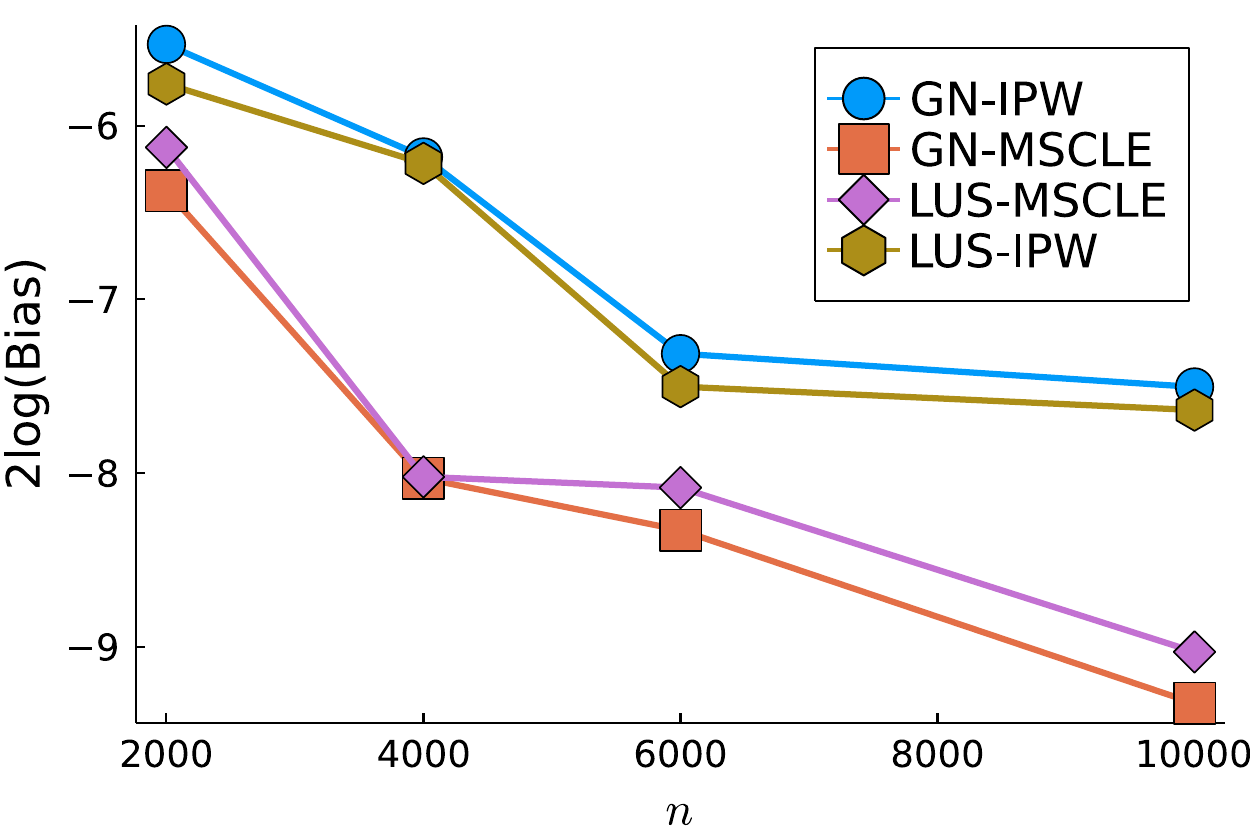}
    \caption{Low sampling ratio}
  \end{subfigure}
  \begin{subfigure}{\textwidth}\centering
    \includegraphics[width=0.45\textwidth,page=1]{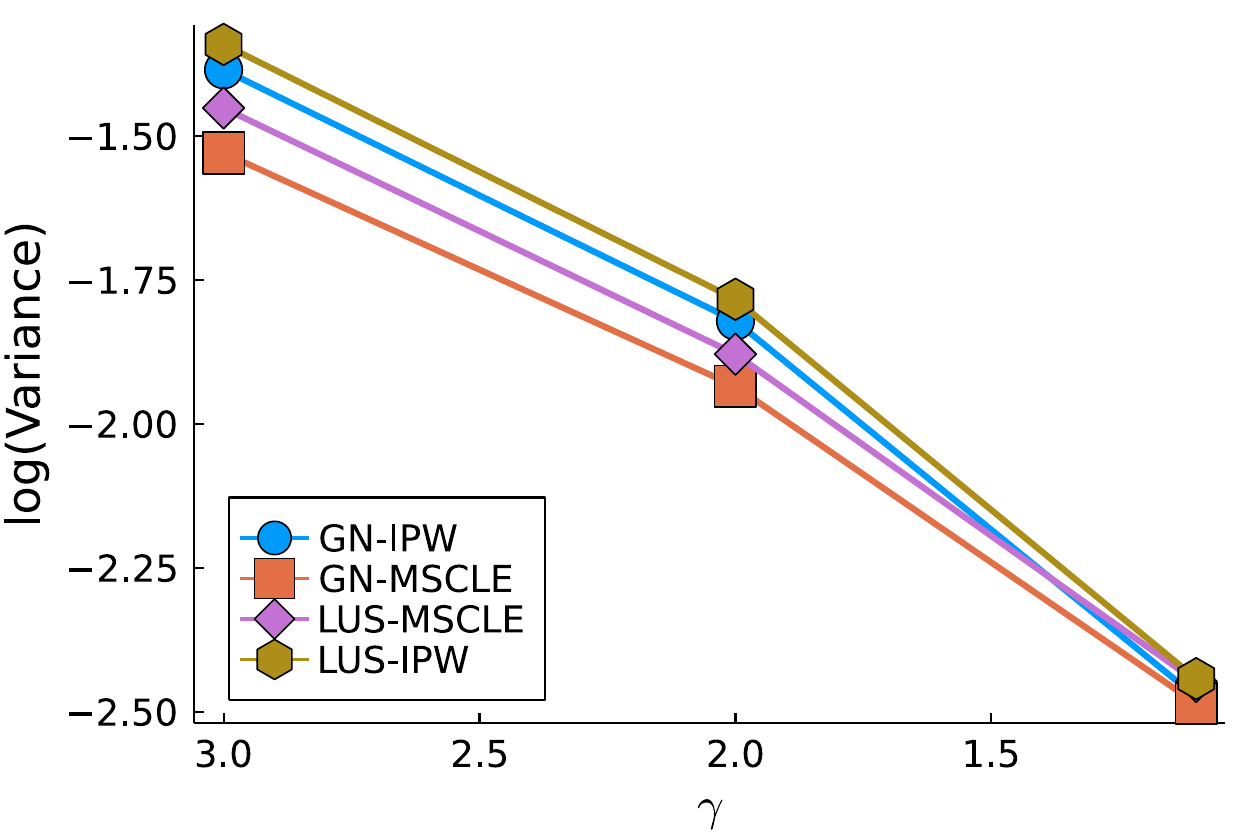}
    \includegraphics[width=0.45\textwidth,page=1]{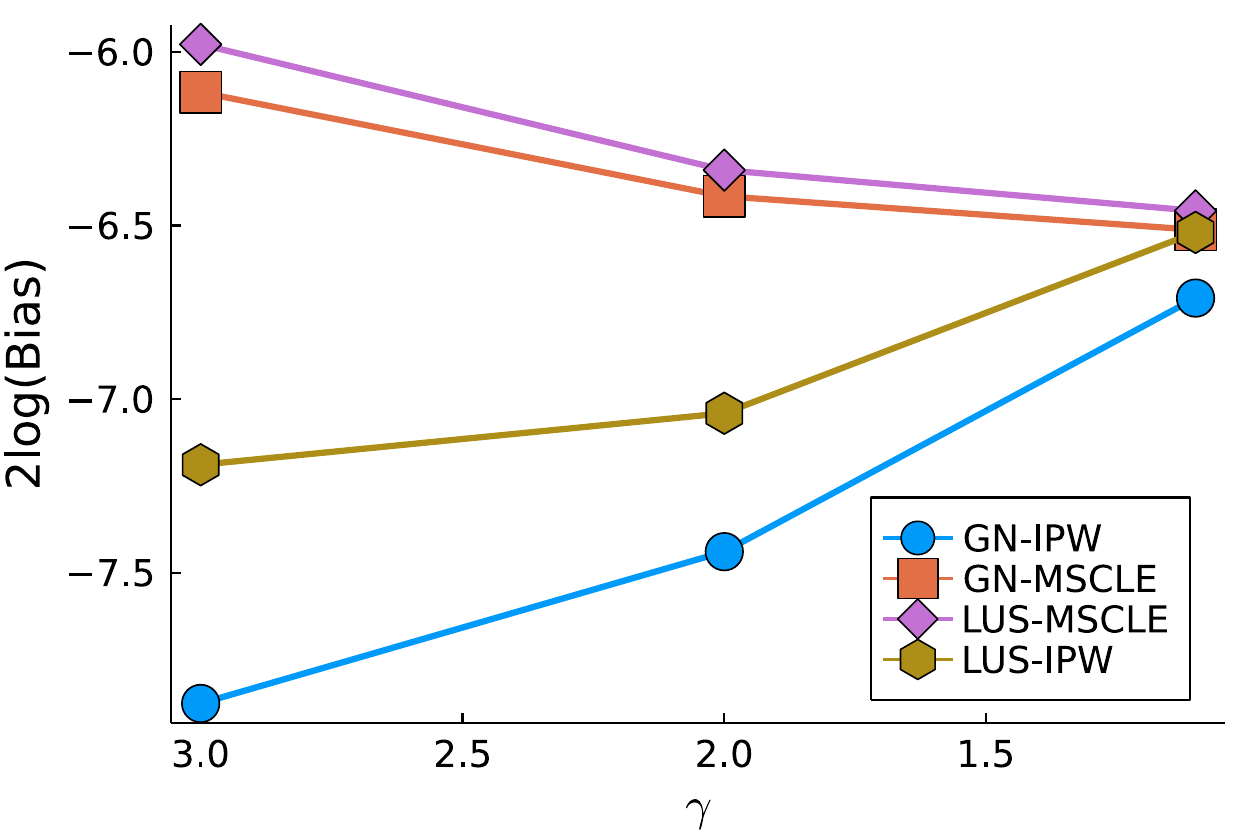}
    \caption{High sampling ratio}
  \end{subfigure}
  \caption{Log of empirical variances and squared biases (the smaller the
    better) of subsample estimators for different sample sizes in multi-class
    logistic regression with the setting of the LUS paper \citep{Han2019}.}
  \label{fig:logistic4}
\end{figure}

To exam the performance of the proposed method with a higher dimension, we also
performed experiments using a setting considered in \cite{Han2019}. Specifically,
the conditional distribution of $\x_{-1,i}$ given $\y_i=\1_k$ is
$\Nor(\bm\mu_k, \bm\Omega_k)$, for $k=1, 2$, and $3$. Here, $\bm\mu_1$ is a 20
dimensional vector with the first ten elements being one and the last ten
elements being zero; $\bm\mu_2$ is a 20 dimensional vector with the first ten
elements being zero and the last ten elements being one; $\bm\mu_2$ is a 20
dimensional vector of zeros, and $\bm\Omega_s$'s are all equal to the identity
matrix. The marginal distribution of $\y_i$ is
$\Pr(\y_i=\1_1)=\Pr(\y_i=\1_3)=0.1$ and $\Pr(\y_i=\1_2)=0.8$. We consider two
settings of sampling rates. The first setting is with
$N=10^6$ and $n=2000, 4000, 6000$ and $10000$, so the sampling rate $n/N$ is
low. For this low sampling rate, we can control the average sample size by
scaling the LUS probability without affecting its form. The second setting is
the same as that in \cite{Han2019} with $N=50000$ and $\gamma=3, 2$, and
$1.1$. For this setting the sampling rate is high and we cannot control the
average sample size for the LUS probability without changing its form. We
implemented the LUS probability with the given values of $\gamma$, and then
implemented the GN sampling probability with the actual average sample size
acquired by the LUS probability. 
Figure~\ref{fig:logistic4} presents the empirical variances and squared biases
for the two settings. The overall pattern here is similar to that with the low
dimensional covariates. We omit the results for the naive estimator here because
it is biased as seen in previous examples. 

\newpage
\subsection{Poisson regression}
\label{sec:poisson-regression-2}
In this section, we consider the Poisson regression model discussed in Section~\ref{sec:poisson-regression}. As in the previous example, we set the full data sample size $N=10^6$, and let the subsample sizes be $n=500; 1000; 1500;$ and $2000$. We let the dimension of the covariates $\x_i=(1,\x_{-1,i}\tp)\tp$'s be $d=7$ with the first element of one for the slope parameters. We set the true value of $\bbeta$ as a vector of 0.25's to generate the data, and consider the following distributions for the covariates corresponding to the slope parameters, $\x_{-1,i}$'s.
\begin{enumerate}[(a)]
\item Independent uniform distribution $\x_{-1,i}\sim\mathbb{U}(\bm{0}, \bm{1})$, where components of $\x_{-1,i}$ independently follow the standard uniform distribution. This distribution has a bounded support and it is symmetric.
\item Independent Beta distribution $\x_{-1,i}\sim\mathbb{B}(\bm{2}, \bm{5})$, where components of $\x_{-1,i}$ independently follow the Beta distribution with parameters 2 and 5. This distribution has a bounded support and it is skewed to the right.
\item Multivariate normal distribution $\x_{-1,i}\sim\Nor(\bm{0}, \bm\Omega)$, where $\bm\Omega$ is defined similarly to that defined in Section~\ref{sec:multi-class-logistic}. This distribution has an unbounded support and it is symmetric. %
\item Independent exponential distribution, $\x_{-1,i}\sim\mathbb{EXP}(2)$, where components of $\x_{-1,i}$ independently follow the exponential distribution with rate parameter 2. This distribution has an unbounded support, and it is asymmetric and positively skewed. 
\end{enumerate}
Here, the distribution in case (a) is the Case 1 used in \cite{yu2020quasi}. We also considered other cases of distributions used in their paper based on uniform distributions. The results are omitted because the relative performance of the three estimators are similar to that of case (a).

\begin{figure}[H]
  \centering
  \begin{subfigure}{\textwidth}\centering
    \includegraphics[width=0.45\textwidth,page=1]{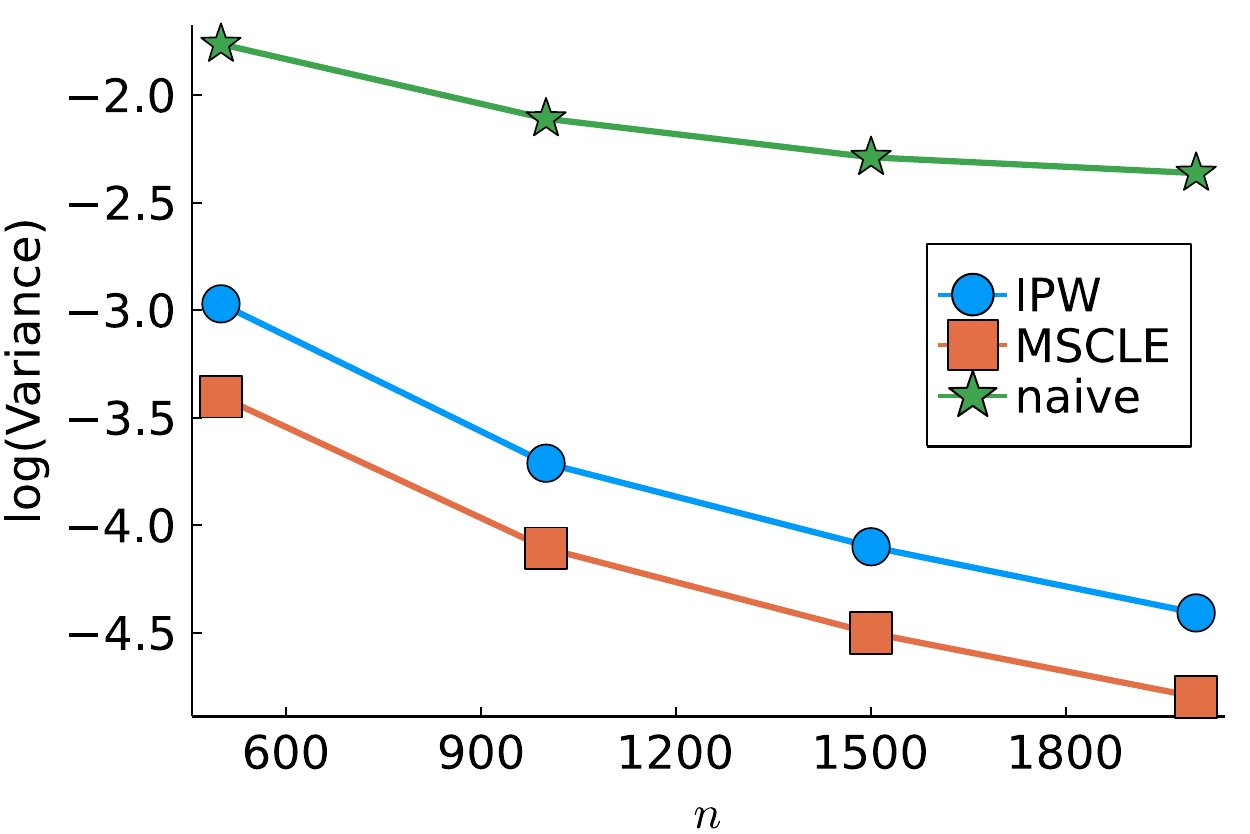}
    \includegraphics[width=0.45\textwidth,page=1]{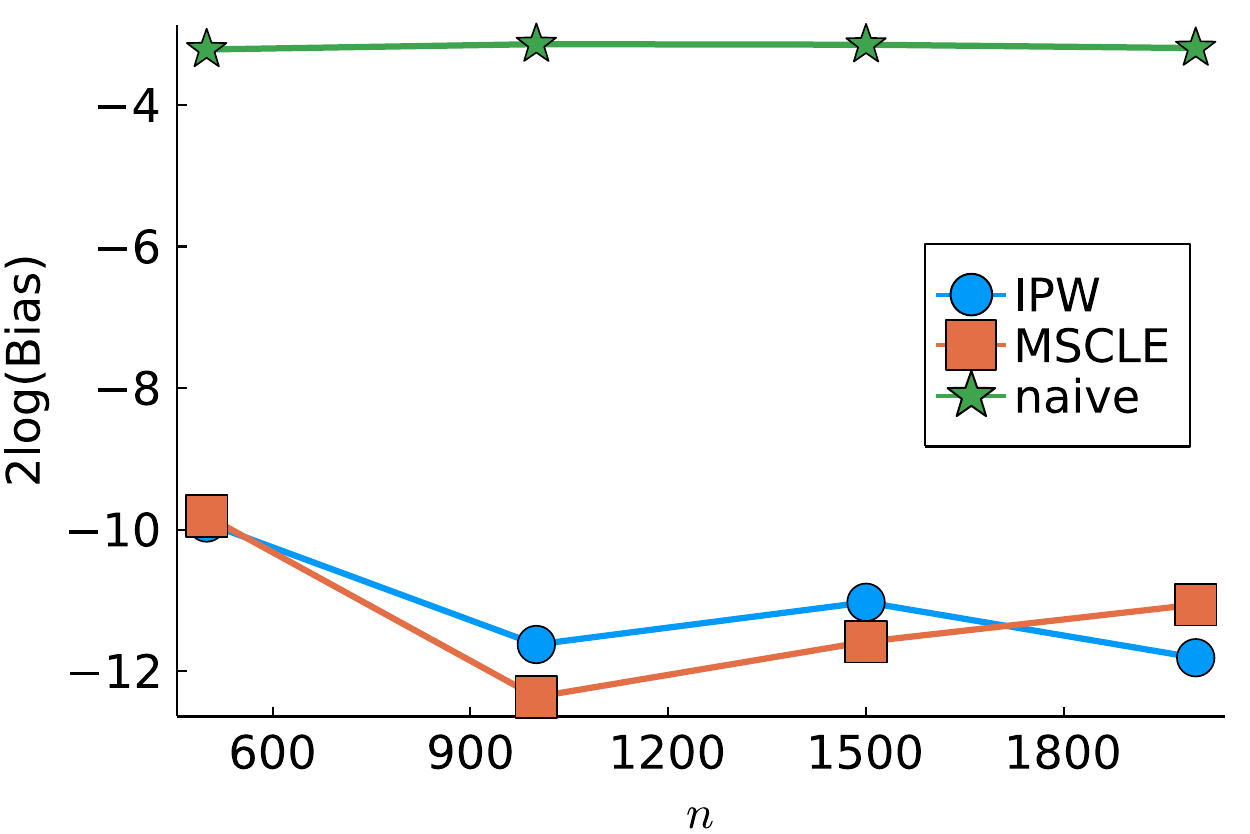}
    \caption{$\x_{-1,i}$'s are uniform}
  \end{subfigure}
  \begin{subfigure}{\textwidth}\centering
    \includegraphics[width=0.45\textwidth,page=2]{figures/00mse_vPoisson.pdf}
    \includegraphics[width=0.45\textwidth,page=2]{figures/00mse_bPoisson.pdf}
    \caption{$\x_{-1,i}$'s are Beta}
  \end{subfigure}
  \begin{subfigure}{\textwidth}\centering
    \includegraphics[width=0.45\textwidth,page=3]{figures/00mse_vPoisson.pdf}
    \includegraphics[width=0.45\textwidth,page=3]{figures/00mse_bPoisson.pdf}
    \caption{$\x_{-1,i}$'s are normal}
  \end{subfigure}
  \begin{subfigure}{\textwidth}\centering
    \includegraphics[width=0.45\textwidth,page=4]{figures/00mse_vPoisson.pdf}
    \includegraphics[width=0.45\textwidth,page=4]{figures/00mse_bPoisson.pdf}
    \caption{$\x_{-1,i}$'s are exponential}
  \end{subfigure}
  \caption{Log of empirical variances and squared biases (the smaller the better) of subsample estimators for different sample sizes in Poisson regression.}
  \label{fig:poisson1}
\end{figure}

Again, we repeat the simulation for $R=1000$ to calculate the empirical 
MSEs, variances, and squared biases, for the three estimators: the proposed
MSCLE estimator, the IPW estimator, and the naive estimator (naive). Results on
the empirical variances and squared biases are presented in Figure~\ref{fig:poisson1}. 

The simulation results in Figure~\ref{fig:poisson1} shows similar patterns of
the simulation study in Section~\ref{sec:multi-class-logistic}. The
variance is the dominating term in the MSE for both the IPW estimator and the MSCLE.
The bias of the naive estimator does not decrease with the subsample size. Both
the IPW estimator and the MSCLE have smaller variances for larger sample
sizes. The MSCLE is uniformly more efficient than the other estimators in terms
of the variances.

For the computational costs, the relative pattern is very similar to that
for the multi-class logistic regression and thus we omit the results. 

We also have additional numerical results on the effect of
pilot misspecification for Poisson regression. Please see them in
Section~\ref{sec:poisson-regression-1} of the Appendix.

\begin{figure}[htp]%
  \centering
  \begin{subfigure}{0.485\textwidth}
    \includegraphics[width=\textwidth,page=1]{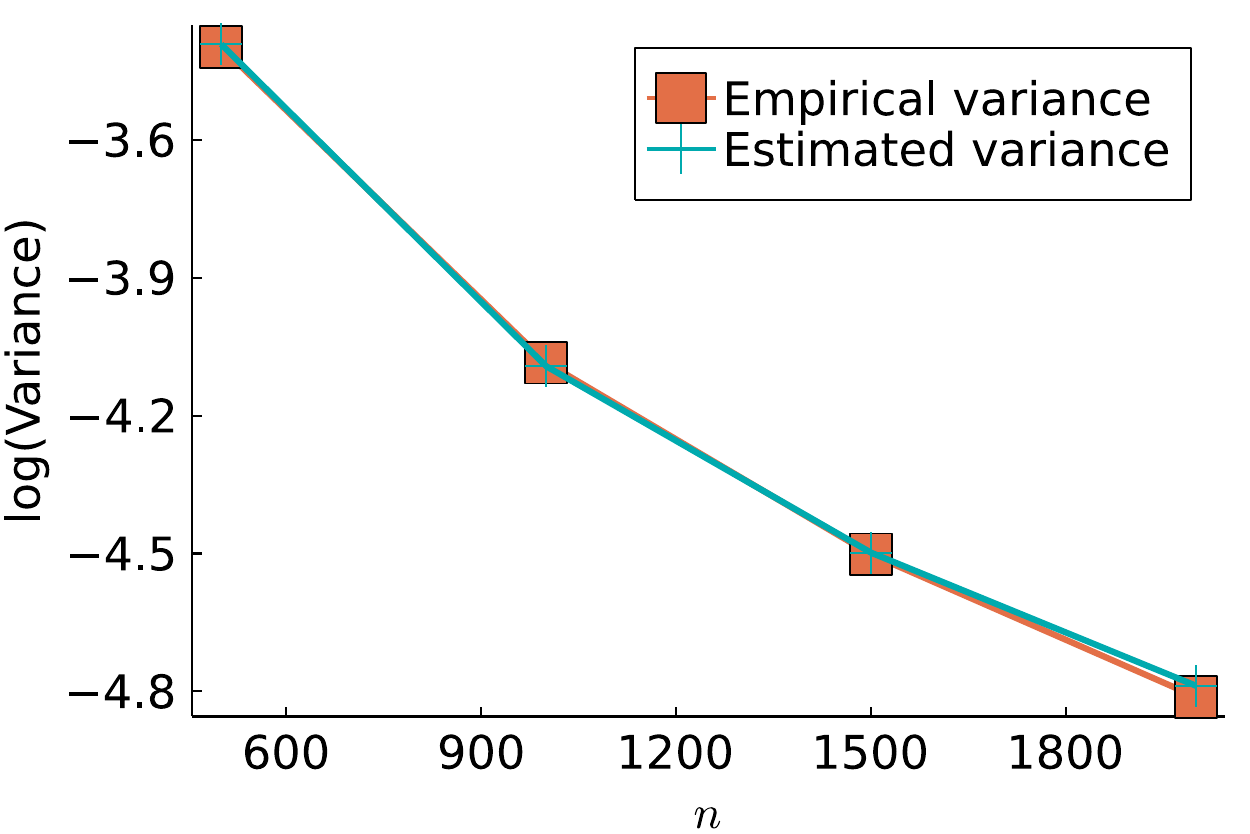}
    \caption{$\x_{-1,i}$'s are uniform}
  \end{subfigure}
  \begin{subfigure}{0.485\textwidth}
    \includegraphics[width=\textwidth,page=2]{figures/00VarEstPoisson.pdf}
    \caption{$\x_{-1,i}$'s are Beta}
  \end{subfigure}
  \begin{subfigure}{0.485\textwidth}
    \includegraphics[width=\textwidth,page=3]{figures/00VarEstPoisson.pdf}
    \caption{$\x_{-1,i}$'s are normal}
  \end{subfigure}
  \begin{subfigure}{0.485\textwidth}
    \includegraphics[width=\textwidth,page=4]{figures/00VarEstPoisson.pdf}
    \caption{$\x_{-1,i}$'s are exponential}
  \end{subfigure}
  \caption{Log of empirical variances and estimated variance of subsample
    estimators for different sample sizes in Poisson
    regression. Logarithm is taken for better presentation.}
  \label{fig:poisson3}
\end{figure}

Figure~\ref{fig:poisson3} plots the empirical variances and estimated
variances to evaluate the performance of the asymptotic
variance in Theorem~\ref{thm:1p} for Poisson regression. We see that the
empirical variances and the estimated variances are very close, showing that the
approximation based on the asymptotic distribution in Theorem~\ref{thm:1p} is
accurate.

\subsection{Real data example: cover type data}
To demonstrate the performance of the MSCLE, we applied it to the cover type
data \citep{blackard1999comparative}. This dataset contains $N=581,012$
observations on ten quantitative variables that measure geographical features
and lighting conditions. The interest is to use these variables to build a model
to predict the forest type. There are $K=7$ forest types: Spruce/Fir
($36.46\%$), Lodgepole Pine ($48.76\%$), Ponderosa Pine ($6.15\%$), Aspen
($1.63\%$), Douglas-fir ($2.99\%$), Krummholz ($3.53\%$) and Cottonwood/Willow
($0.427\%$), where the number in the parentheses is the percentage of the forest
type in the full data set. We include an intercept in the model so $d=11$ and
the dimension of unknown parameter is $(K-1)d=66$. To fit a multi-class logistic
regression model as in~(\ref{eq:16}), we use the probabilities defined in
(\ref{eq:15}) to take subsamples of sizes $n$ from the full data, and use the
subsamples to train the model. Pilot estimates are obtained from subsamples of
average size 2000 taken according to proportional subsampling
probabilities. Note that for real data the true data generating model is
unknown and any parametric model may be subject to a certain level of model
misspecification. We use the full data estimator as the ``true parameter'' and repeat the subsampling estimator for 1000 times to calculated empirical biases, variances, and MSEs.

Results are presented in Table~\ref{tab:1}. We present the squared biases so that they are directly comparable to the variances. 
For comparison, we obtained the IPW estimator (IPW) and the naive
estimator (Naive), and also implement the uniform subsampling method
(Uniform). In terms of MSE, MSCLE outperforms the IPW estimator, which
outperforms the Uniform method. %
The Naive method has small variances, but its biases are large and do not seem to converge to zero.

\begin{table}[H]
\caption{Empirical bias (squared), Variance (Var.) and the mean squared error (MSE) of the four subsample estimators for cover type data}
\label{tab:1}
\centering
\begin{tabular}{l|ccc|ccc|ccc}\hline
         & \multicolumn{3}{|c|}{$n=4000$} & \multicolumn{3}{|c|}{$n=6000$} & \multicolumn{3}{|c}{$n=8000$} \\
         \cline{2-10} 
         & Bias$^2$ & Var. & MSE & BIAS$^2$ & Var. & MSE & Bias$^2$ & Var. & MSE \\
         \hline 
MSCLE     & 25.73 &  248.27 & 273.74 & 14.53 &  120.17 & 134.58  & 12.56 & 76.87 & 89.35     \\ 
IPW & 22.97 &  325.70 & 348.35 & 6.76 & 163.67 & 170.27  & 3.16 &  104.57 & 107.62 \\ 
Naive    & 146.04 & 231.76 & 377.57  & 146.73 & 118.05 & 264.66 & 150.49 & 80.78 & 231.49     \\ \hline
\end{tabular}
\end{table}

\subsection{Real data example: the MNIST data}
\label{sec:real-data-example}
In this section, we illustrate the advantage of the MSCLE over the IPW estimator
using the famous MNIST data that is available at
\url{http://yann.lecun.com/exdb/mnist/}. The data has a training set with 60,000
instances and a testing set with 10,000 instances. Each instance is an image of
a handwritten digit with 28 by 28 greyscale pixels. The goal is to train a model
using the training set to predict the handwritten digits of \{0, 1, 2, ..., 9\}
in the testing set. We implement the convolutional neural network LeNet-5
\citep{lecun1998gradient} with Flux.jl \citep{innes:2018}, and use a subsample
of average size $n=5,000$ out of $N=60,000$ (about 8.3\% of the training data)
to train the model.

To apply the sampling probabilities presented in~(\ref{eq:15}), we use the norms
of the of the 28 by 28 greyscale pixels to replace $\|\x_i\|$'s. Note that this
does not give the GN sampling probability for this model. We adopt this
simplified version of the sampling probability because the complicated model
structure and high dimension (with 44,426 parameters) makes it difficult to calculate
the real gradient norms. The pilot probabilities $\p(\x_i,\tbeta_{\plt})$'s are
obtained from a pilot model trained with 2,000 uniformly selected instances from
the training set.

\begin{figure}[htp]%
  \centering
  \includegraphics[width=\textwidth]{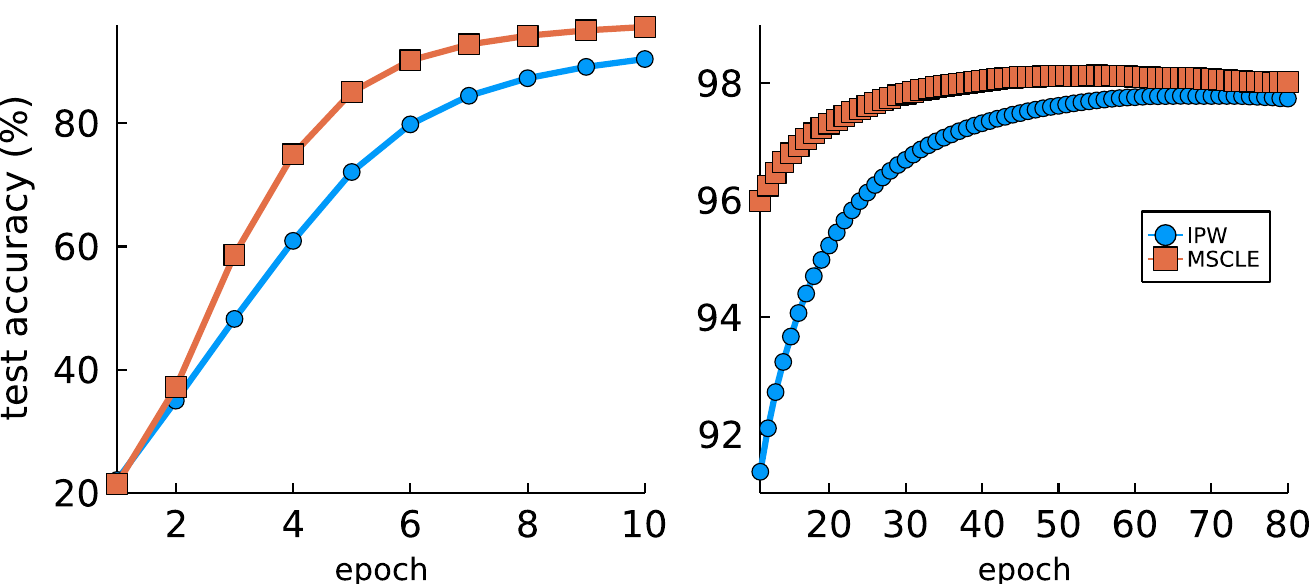}
  \caption{Classification accuracy (in percentage) on the test data against epoch in the training using subsamples of size $n=5,000$ from the MNIST data.}
  \label{fig:mnist}
\end{figure}

Due to the complexity of the model the parameters do not have explicit
interpretations, so we report the test accuracy which is the percentage of correct
classification in the testing set. Figure~\ref{fig:mnist} plots the test
accuracy against the epoch for the MSCLE and the IPW estimator. The epoch value
is the number of passes throughout the entire training set that the mini batch
gradient descent algorithm has completed. It is seen that the MSCLE outperforms
the IPW estimator uniformly across epochs.

Again, no model can be the exact data generating model for real data, so
although being a very rich model, the LeNet-5 model here may be subject to a
certain level of misspecification as well. As a matter of fact, since the number
of parameters is much larger than the subsample size, some of the regularity
assumptions %
in Section~\ref{sec:estim-param-subs} may not hold. The promising performance in
Figure~\ref{fig:mnist} indicates that the proposed MSCLE may be applicable in
more general scenarios than that restricted by the regularity assumptions in
Section~\ref{sec:estim-param-subs}.

\section{Discussion} 

Subsampling is a useful technique for handling big data. To estimate the parameters with the subsample data, the inverse probability IPW estimator have been used as a gold standard method. In this paper, we consider an alternative method of parameter estimation using the sampled conditional likelihood function. The resulting MSCLE is consistent and is more efficient than the IPW estimator. The computation for obtaining MSCLE requires computing the bias-adjustment term in the naive (complete-case) method. Explicit closed-form formula for the bias-adjustment terms are given in Section~\ref{sec:gener-line-models}. 

To obtain efficient subsampling probabilities, we need  consistent estimates of the model parameters. In this paper, we assume that an independent pilot subsample is available outside of the current sample. In this case, the theory can be developed as presented in Section~\ref{sec:estim-param-subs}. Otherwise, one can consider an adaptive estimation method which simultaneously updates the parameter estimates and the corresponding selection probabilities. Such adaptive methods introduce  additional computational burden and also theoretical challenges. In addition, the subsampling idea is closely related with reservior sampling \citep{efraimidis2006} which is a useful tool for handling streaming data.  Thus, the proposed MSCLE can be used to handle the reservior sample with unequal selection probabilities.  Such extensions will be the topics for future research. 

\acks{The authors thank the reviewers for their insightful comments and
  suggestions which greatly helped improve the paper. Wang's research was
  supported by NSF grant CCF-2105571 and a GPU grant from NVIDIA
  Corporation. Kim's research was supported by NSF grant CSSI-1931380.}

\appendix

\setcounter{equation}{0}
\renewcommand{\theequation}{A.\arabic{equation}}
\setcounter{figure}{0}
\renewcommand{\thefigure}{A.\arabic{figure}}

\section{Proofs and technical details}
\subsection{Proof of Theorem~\ref{thm:1}}
For any functions $h_1(\x,y)$ and $h_2(\x)$, assuming integrability in the following, we have
\begin{align}
  \Exp\{\delta h_1(\x,y)h_2(\x)\}\notag
  &=\Exp\big[\delta h_2(\x)\Exp\{h_1(\x,y)\mid\delta,\x\}\big]\\
  &=\Exp\Big[\pi(\x,y)h_2(\x)
    \Exp\{h_1(\x,y)\mid\x,\delta=1\}\Bigm|\x,y\Big]\notag\\
  &=\Exp\Big[\delta h_2(\x)\Exp\{h_1(\x,y)\mid\x,\delta=1\}\Bigm|\x,y\Big]
    \notag\\
  &=\Exp\big[\delta h_2(\x)\Exp\{h_1(\x,y)\mid\x,\delta=1\}\big].\label{eq:27}
\end{align}

    The estimator $\htheta_{S}$ is the maximizer of $\ell_{S}(\btheta)$ in (\ref{eq:5}), 
so $\sqrt{N}(\htheta_{S}-\btheta)$ is the maximizer of $\gamma(\eeta)=\ell_{S}(\btheta+\eeta/\sqrt{N})-\ell_{S}(\btheta)$. By Taylor's expansion,
\begin{equation}\label{eq:31}
  \gamma(\eeta)
  =\frac{1}{\sqrt{N}}\eeta\tp\dot\ell_{S}(\btheta)
    +\frac{1}{2N}\eeta\tp\ddot\ell_{S}(\btheta)\eeta+R,
\end{equation}
where $R=\op$ because Assumption~\ref{asmp:2} indicates that
\begin{equation}
  |R|\le\frac{d^3\|\eeta\|^3}{3N^{1/2}}\times\oneN\sumN B(\x_i,y_i)=\op.
\end{equation}

For $\dot\ell_{S}(\btheta)=\sumN\dot\ell_{S}(\btheta;\x_i,y_i)$, $\dot\ell_{S}(\btheta;\x_i,y_i)$ are i.i.d. random vectors. Let $h_1(\x_i,y_i)=\dot\ell(\btheta;\x_i,y_i)$ and $h_2(\x_i)=1$ in (\ref{eq:27}), we know that $\Exp\{\dot\ell_{S}(\btheta;\x_i,y_i)\}=\0$, and using (\ref{eqn:10}) we have
\begin{align}
  \Var\{\dot\ell_{S}(\btheta;\x_i,y_i)\}
  &=\Exp\Big(\delta[\dot\ell(\btheta;\x,y)
    -\Exp\{\dot\ell(\btheta;\x,y)\mid\x,\delta=1)\}]^{\otimes2}\Big)\\
  &=\Exp\Big(\delta[\dot\ell^{\otimes2}(\btheta;\x,y)
    -\Exp^{\otimes2}\{\dot\ell(\btheta;\x,y)\mid\x,\delta=1)\}]\Big)\\
  &=\Exp\Bigg(\delta\bigg[\dot\ell^{\otimes2}(\btheta;\x,y)
    -\frac{\Exp^{\otimes2}
    \{\dot\ell(\btheta;\x,y)\pi(\x,y)\mid\x\}}
    {\bar{\pi}^2(\x; \btheta)}\bigg]\Bigg)
    =\bSigma_{\btheta}.\label{eq:18}
\end{align}
Thus, from the central limit theorem,
\begin{equation}\label{eq:29}
  \frac{1}{\sqrt{N}}\dot\ell_{S}(\btheta)\cvd\Nor(\0,\bSigma_{\btheta}).
\end{equation}
Now we investigate the Hessian matrix
\begin{align}
  \oneN\ddot\ell_{S}(\btheta)
  =&\oneN\sumN\delta_i[\ddot\ell(\btheta;\x_i,y_i)
     -\Exp\{\ddot\ell(\btheta;\x_i,y_i)\mid\x_i,\delta_i=1)\}]\notag\\
  &\quad -\oneN\sumN\delta_i\bigg[\int\dot\ell(\btheta;\x_i,y_i)
     \dot f\tp(y_i\mid\x_i,\delta_i=1;\btheta)\ud y\bigg]\notag\\
  \equiv&\Delta_1-\H_N.
\end{align}
Note that $\Delta_1$ is an average of i.i.d terms, and 
from Assumption~\ref{asmp:1} and (\ref{eq:27}), we know 
that $\Exp(\Delta_1)=\0$. Thus, from the strong law of large numbers,
$\Delta_1=o(1)$ almost surely.

  Under Assumptions~\ref{asmp:1} and \ref{asmp:2}, and the fact that
$\pi(\x,y)$ is bounded above by one, derivatives can pass the integration sign
in $\int f(y\mid\x;\btheta)\pi(\x,y)\ud y$ by the dominated convergence
theorem. 
Thus, using (\ref{eq:36}) we have
\begin{align}
  &\dot f(y_i\mid\x_i,\delta_i=1;\btheta)\notag\\
  &=\frac{\dot f(y_i\mid\x_i;\btheta)\pi(\x_i,y_i)}
  {\bar{\pi}(\x_i; \btheta)}%
    -\frac{f(y_i\mid\x_i;\btheta)\pi(\x_i,y_i)
  \Exp\{\dot\ell(\btheta;\x_i,y_i)\pi(\x_i,y)\mid\x_i\}}
  {\bar{\pi}^2(\x_i; \btheta)}\notag\\
  &=\frac{\dot\ell(\btheta;\x_i,y_i)\pi(\x_i,y_i)f(y_i\mid\x_i;\btheta)}
  {\bar{\pi}(\x_i; \btheta)}%
    -\frac{f(y_i\mid\x_i;\btheta)\pi(\x_i,y_i)
  \Exp\{\dot\ell(\btheta;\x_i,y_i)\pi(\x_i,y)\mid\x_i\}}
  {\bar{\pi}^2(\x_i; \btheta)}.
\end{align}
Thus, we have that
\begin{align}
  \H_N
  &=\oneN\sumN\delta_i\bigg[\int\dot\ell(\btheta;\x_i,y_i)
     \dot f\tp(y\mid\x_i,\delta_i=1;\btheta)\ud y\bigg]\\
  &=\oneN\sumN\delta_i\bigg[
    \frac{\Exp\{\dot\ell^{\otimes2}(\btheta;\x_i,y_i)
    \pi(\x_i,y_i)\mid\x_i\}}
    {\bar{\pi}(\x_i; \btheta)}
    -\frac{\Exp^{\otimes2}\{\dot\ell(\btheta;\x_i,y_i)
    \pi(\x_i,y)\mid\x_i\}}{\bar{\pi}^2(\x_i; \btheta)}\bigg].
\end{align}
Since $\H_N$ is an average of i.i.d terms, and (\ref{eqn:10}) and (\ref{eq:27}) tells us that its expectation is $\Exp(\H_N)=\bSigma_{\btheta}$, from the strong law of large numbers, $\H_{N}\rightarrow\bSigma_{\btheta}$ almost surely. 
Therefore,
\begin{equation}\label{eq:30}
  -\oneN\ddot\ell_{S}(\btheta)=\bSigma_{\btheta}+\op.
\end{equation}

From (\ref{eq:31}), (\ref{eq:29}) and (\ref{eq:30}), applying the Basic Corollary in page 2 of \cite{hjort2011asymptotics} finishes the proof.

\subsection{Proofs of Theorem~\ref{thm:2}}

  Let $\ell_{S}(\btheta;\x,y)=\log\{f(y\mid\x,\delta=1;\btheta)\}$. Note that $\Exp\{\delta\u(\btheta;\x,y)\mid\x\}=\0$. Using
\begin{align*}
  \frac{\partial}{\partial\btheta\tp}
  \Exp\{\delta\u(\btheta;\x,y)\mid\x\}
  =\Exp\{\delta\dot{\u}(\btheta;\x,y)\mid\x\}
  +\Exp\{\delta\u(\btheta;\x,y)\dot\ell_{S}\tp(\btheta;\x,y)\mid\x\},
\end{align*}
we have
\begin{equation*}
  \Exp\{\delta\u(\btheta;\x, y)\dot\ell_{S}\tp(\btheta;\x,y)\}
  = -\Exp\{\delta\dot{\u}(\btheta; \x, y)\}
  = -\M_{\btheta}.
\end{equation*}
Hence, by calculating the variance covariance matrix
\begin{align*}
  &\Var\big\{\delta\M_{\btheta}^{-1}\u(\btheta;\x, y)
    +\delta\bSigma_{\btheta}^{-1}\dot\ell_{S}(\btheta;\x,y)\big\}\\
  &=\Var\big\{\delta\M_{\btheta}^{-1}\u(\btheta;\x,y)\big\}
    +\M_{\btheta}^{-1}\Exp\{\delta\u(\btheta;\x,y)\dot\ell_{S}\tp(\btheta;\x,y)\}
    \bSigma_{\btheta}^{-1}\\
  &\quad +\bSigma_{\btheta}^{-1} 
    \Exp\{\delta\dot\ell_{S}(\btheta;\x,y)\u\tp(\btheta;\x,y)\}(\M_{\btheta}^{-1})\tp
    +\bSigma_{\btheta}^{-1}\Var\{\delta\dot\ell_{S}(\btheta;\x,y)\}
    \bSigma_{\btheta}^{-1}\\
  &=\M_{\btheta}^{-1}\Exp\{\delta\u^{\otimes2}(\btheta;\x,y)\}(\M_{\btheta}^{-1})\tp
  -\bSigma_{\btheta}^{-1}\ge\0,
\end{align*}
we have proved (\ref{eq:59}). Now, 
we know that the equality holds if $T(\btheta):= \delta\M_{\btheta}^{-1}\u(\btheta;\x, y)
    +\delta\bSigma_{\btheta}^{-1}\dot\ell_{S}(\btheta;\x,y)$ is a constant vector.
    Since $\Exp\{T(\btheta) \}=\0$, we have 
    $T( \btheta) = \0$, or  $\u(\btheta;\x, y)=-
    \M_{\btheta}\bSigma_{\btheta}^{-1}\dot\ell_{S}(\btheta;\x,y)$. %
    This completes the proof.

\subsection{Proof of Theorem~\ref{thm:1p}}
\label{sec:proof-theor-refthm:1}

  Under Assumptions~\ref{asmp:1p} and \ref{asmp:2p}, and the fact that
  $\pi_{N}(\x,y;\bvtheta)$ is bounded above by one, we know that derivatives can
  pass the integration sign in $\int
  f(y\mid\x;\btheta)\pi_{N}(\x,y;\tilde\bvtheta_{\plt})\ud y$ by the dominated
  convergence theorem. 
  The score function for $\ell_{S}(\btheta\mid\tilde\bvtheta_{\plt})$ can be
  written as
  \begin{align}
  \dot\ell_{S}(\btheta\mid\tilde\bvtheta_{\plt})
  &=\sumN\delta_i[\dot\ell(\btheta;\x_i,y_i)
    -\int\dot\ell(\btheta;\x_i,y)f(y\mid\x_i,\delta=1;\tilde\bvtheta_{\plt})\ud y],\\
  &=\sumN\delta_i[\dot\ell(\btheta;\x_i,y_i)
    -\Exp\{\dot\ell(\btheta;\x_i,y_i)\mid\x_i,\delta=1;\tilde\bvtheta_{\plt}\}],
  \end{align}
  where
  \begin{equation}
    f(y_i\mid\x_i,\delta_i=1,\tilde\bvtheta_{\plt})=
  \frac{f(y_i\mid\x_i;\tilde\bvtheta_{\plt})\pi(\x_i,y_i)}
  {\int f(y\mid\x_i;\tilde\bvtheta_{\plt})\pi(\x_i,y)\ud y}.
  \end{equation}
Similarly to (\ref{eq:27}), for any functions $h_1(\x,y)$ and $h_2(\x)$, we have
\begin{align}
  \Exp\{\delta h_1(\x,y)h_2(\x)\mid\tilde\bvtheta_{\plt}\}\notag
  &=\Exp\Big[\delta h_2(\x)\Exp\{h_1(\x,y)\mid\delta,\x\}\Bigm|\tilde\bvtheta_{\plt}\Big]\\
  &=\Exp\Big[\pi_{N}(\x,y;\tilde\bvtheta_{\plt})h_2(\x)
    \Exp\{h_1(\x,y)\mid\x,\delta=1\}\Bigm|\x,y;\tilde\bvtheta_{\plt}\Big]\notag\\
  &=\Exp\Big[\delta h_2(\x)\Exp\{h_1(\x,y)\mid\x,\delta=1\}\Bigm|\x,y;\tilde\bvtheta_{\plt}\Big]
    \notag\\
  &=\Exp\big[\delta h_2(\x)\Exp\{h_1(\x,y)\mid\x,\delta=1\}\bigm|\tilde\bvtheta_{\plt}\big]\label{eq:44}
\end{align}

Now we prove Theorem~\ref{thm:1p}. 
    Since $\htheta_{s,\tilde\bvtheta_{\plt}}$ is the maximizer of $\ell_{S}(\btheta\mid\tilde\bvtheta_{\plt})$,  $\sqrt{n}(\htheta_{s,\tilde\bvtheta_{\plt}}-\btheta)$ is the maximizer of $\gamma(\eeta\mid\tilde\bvtheta_{\plt})=\ell_{S,\tilde\bvtheta_{\plt}}(\btheta+\eeta/\sqrt{n})-\ell_{S}(\btheta\mid\tilde\bvtheta_{\plt})$. By Taylor's expansion,
\begin{equation}\label{eq:46}
  \gamma(\eeta\mid\tilde\bvtheta_{\plt})
  =\frac{1}{\sqrt{n}}\eeta\tp\dot\ell_{S}(\btheta\mid\tilde\bvtheta_{\plt})
    +\frac{1}{2n}\eeta\tp\ddot\ell_{S}(\btheta\mid\tilde\bvtheta_{\plt})\eeta+R,
\end{equation}
where $R=\op$. We show why $R$ is a small term in probability in the following. By (\ref{eq:61}) and (\ref{eq:32}) of Assumption~\ref{asmp:2p}, $R$ satisfies that
\begin{align}
  |R|\le\frac{d^3\|\eeta\|^3}{3}\times\frac{1}{n^{3/2}}\sumN\delta_iB_{\tilde\bvtheta_{\plt}}(\x_i,y)
  \equiv\frac{d^3\|\eeta\|^3}{3}\times\Delta.
\end{align}
For any constant $\epsilon>0$, from Markov's inequality,
\begin{align}
  \Pr(\Delta>\epsilon\mid\tilde\bvtheta_{\plt})
  \le\frac{N\int\pi_{N}(\x,y;\tilde\bvtheta_{\plt})B(\x,y)f(y\mid\x;\btheta)\ud y}
  {n^{3/2}\epsilon}=\op,
\end{align}
where the last step above is because (\ref{eq:35}) and the fact that $\tilde\bvtheta_{\plt}$ is independent of the data imply that $Nn^{-1}\Exp\{\pi_{N}(\x,y;\tilde\bvtheta_{\plt})B_{\tilde\bvtheta_{\plt}}(\x,y)
  \mid\tilde\bvtheta_{\plt}\}=\Op$. 
Because a conditional probability is a bounded random variable, $\Pr(\Delta>\epsilon)=\Exp\{\Pr(\Delta>\epsilon\mid\tilde\bvtheta_{\plt})\}\rightarrow0$, indicating that $\Delta=\op$, and therefore $R=\op$.

For $\dot\ell_{S}(\btheta\mid\tilde\bvtheta_{\plt})$, it is a sum of i.i.d terms conditionally on $\tilde\bvtheta_{\plt}$. Let $h_1(\x,y)=\dot\ell(\btheta;\x,y)$ and $h_2(\x)=1$ in (\ref{eq:44}), we know that $\Exp\{\dot\ell_{S}(\btheta\mid\tilde\bvtheta_{\plt})\mid\tilde\bvtheta_{\plt}\}=\0$; and using a similar procedure to obtain~(\ref{eq:18}), we have
\begin{align}
  \onen\Var\{\dot\ell_{S}(\btheta\mid\tilde\bvtheta_{\plt})\mid\tilde\bvtheta_{\plt}\}
  &=\frac{N}{n}\Exp\Big(\delta[\dot\ell(\btheta;\x,y)
    -\Exp\{\dot\ell(\btheta;\x,y)\mid\x,\delta=1,\tilde\bvtheta_{\plt})\}]^{\otimes2}
    \Bigm|\tilde\bvtheta_{\plt}\Big)\notag \\
  &=\frac{N}{n}\Exp\Big(\delta[\dot\ell^{\otimes2}(\btheta;\x,y)
    -\Exp^{\otimes2}\{\dot\ell(\btheta;\x,y)\mid\x,\delta=1,\tilde\bvtheta_{\plt})\}]
    \Bigm|\tilde\bvtheta_{\plt}\Big)\notag \\
  &=\bSigma_{N{\tilde\bvtheta_{\plt}}}\cvp\bSigma_{\btheta,\bvtheta_p},\label{eq:52}
\end{align}
where the convergence in probability in the last step is due to Assumption~\ref{asmp:3}.
Note that $\dot\ell_{S}(\btheta;\x_i,y_i\mid\tilde\bvtheta_{\plt})=\delta_i[\dot\ell(\btheta;\x_i,y_i)
-\Exp\{\dot\ell(\btheta;\x_i,y_i)\mid\x_i,\delta_i=1,\tilde\bvtheta_{\plt})\}]$.
We check the Lindeberg-Feller condition given $\tilde\bvtheta_{\plt}$. For any $\epsilon>0$, %
\begin{align*}
  &\onen\sumN\Exp\big\{\|\dot\ell_{S}(\btheta;\x_i,y_i\mid\tilde\bvtheta_{\plt})\|^2
    I(\|\dot\ell_{S}(\btheta;\x_i,y_i\mid\tilde\bvtheta_{\plt})\|^2>n\epsilon)
    \bigm|\tilde\bvtheta_{\plt}\big\}\\
  &=\frac{N}{n}\Exp\big\{\|\dot\ell_{S}(\btheta;\x_i,y_i\mid\tilde\bvtheta_{\plt})\|^2
    I(\|\dot\ell_{S}(\btheta;\x_i,y_i\mid\tilde\bvtheta_{\plt})\|^2>n\epsilon)
    \bigm|\tilde\bvtheta_{\plt}\big\}\\
  &\le\frac{N}{n^2\epsilon}
    \Exp\big\{\|\dot\ell_{S}(\btheta;\x_i,y_i\mid\tilde\bvtheta_{\plt})\|^4
    \bigm|\tilde\bvtheta_{\plt}\big\}\\
  &=\frac{N}{n^2\epsilon}
    \Exp\big[\delta\|\dot\ell(\btheta;\x,y)
    -\Exp\{\dot\ell(\btheta;\x,y)\mid\x_i,\delta_i=1
    ;\tilde\bvtheta_{\plt})\}\|^4\bigm|\tilde\bvtheta_{\plt}\big]\\
  &\le\frac{8N}{n^2\epsilon}
    \Exp\big[\delta\|\dot\ell(\btheta;\x,y)\|^4
    +\Exp\{\|\dot\ell(\btheta;\x,y)\|^4
    \mid\x_i,\delta_i=1,\tilde\bvtheta_{\plt})\}\bigm|\tilde\bvtheta_{\plt}\big]\\
  &=\frac{16N}{n^2\epsilon}
    \Exp\{\pi_{N}(\x,y;\tilde\bvtheta_{\plt})\|\dot\ell(\btheta;\x,y)\|^4
    \bigm|\tilde\bvtheta_{\plt}\}\\
  &=\op,
\end{align*}
where the last step is because (\ref{eq:34}) in Assumption~\ref{asmp:1} and the fact that $\tilde\bvtheta_{\plt}$ is independent of the data imply that 
\begin{equation}\label{eq:49}
  \frac{N}{n}
    \Exp\{\pi_{N}(\x,y;\tilde\bvtheta_{\plt})\|\dot\ell(\btheta;\x,y)\|^4
    \bigm|\tilde\bvtheta_{\plt}\}=\Op.
\end{equation}
Thus, from the Lindeberg-Feller central limit theorem, conditionally on $\tilde\bvtheta_{\plt}$,
\begin{equation}\label{eq:47}
  \frac{1}{\sqrt{n}}\dot\ell_{S}(\btheta\mid\tilde\bvtheta_{\plt})
  \cvd\Nor(\0,\bSigma_{\btheta,\bvtheta_p}).
\end{equation}
Now we investigate the Hessian matrix. Note that
  \begin{align}
    \Exp\{\dot\ell(\btheta;\x,y)\mid\x,\delta=1,\tilde\bvtheta_{\plt}\}]
    =\int\dot\ell(\btheta;\x,y)f(y\mid\x,\delta=1,\tilde\bvtheta_{\plt})\ud y,
  \end{align}
  where $f(y\mid\x,\delta=1,\tilde\bvtheta_{\plt})$ is the density of $y$ conditional on $\x$, $\delta=1$, and $\tilde\bvtheta_{\plt}$. We have
\begin{align}
  \onen\ddot\ell_{S}(\btheta\mid\tilde\bvtheta_{\plt})
  =&\onen\sumN\delta_i\big[\ddot\ell(\btheta;\x_i,y_i)
     -\Exp\{\ddot\ell(\btheta;\x_i,y_i)
     \mid\x_i, \delta_i=1;\tilde\bvtheta_{\plt}\}\big]\notag\\
  &\quad -\onen\sumN\delta_i\bigg[\int\dot\ell(\btheta;\x_i,y_i)
     \dot f\tp(y_i\mid\x_i,\delta_i=1,\tilde\bvtheta_{\plt})\ud y\bigg]\notag\\
  \equiv&\Delta_2-\H_N^{\tilde\bvtheta_{\plt}}.
\end{align}
From~(\ref{eq:44}), we know that $\Exp(\Delta_2)=\0$. For the $j_1,j_2$-th element of $\Delta_2$, $\Delta_{2,j_1j_2}$,
\begin{align*}
  &\Exp\Big(\Delta_{2,j_1,j_2}^2\Bigm|\tilde\bvtheta_{\plt}\Big)\\
  &=\frac{1}{n}\sumN\Exp\Big(
    \delta_i\big[\ddot\ell_{j_1j_2}(\btheta;\x_i,y_i)
    -\Exp\{\ddot\ell_{j_1j_2}(\btheta;\x_i,y_i)
    \mid\x_i,\delta_i=1,\tilde\bvtheta_{\plt})\}\big]^2
    \Bigm|\tilde\bvtheta_{\plt}\Big)\\
  &\le\frac{2}{n^2}\sumN\Exp\Big(
    \delta_i\big[\|\ddot\ell(\btheta;\x_i,y_i)\|^2
    +\|\Exp\{\ddot\ell(\btheta;\x_i,y_i)\mid\x_i,\delta_i=1;
    \tilde\bvtheta_{\plt})\|^2\big]\Bigm|\tilde\bvtheta_{\plt}\Big)\\
  &\le\frac{2}{n^2}
    \sumN\Exp\Big(\delta_i[\|\ddot\ell(\btheta;\x_i,y_i)\|^2
    +\Exp\{\|\ddot\ell(\btheta;\x_i,y_i)\|^2
    \mid\x_i,\delta_i=1,\tilde\bvtheta_{\plt})]\Bigm|\tilde\bvtheta_{\plt}\Big)\\
  &=\frac{4N}{n^2}
    \Exp\Big\{\pi_{N}(\x,y;\tilde\bvtheta_{\plt})\|\ddot\ell(\btheta;\x,y)\|^2
    \Bigm|\tilde\bvtheta_{\plt}\Big\}=\op,
\end{align*}
where the last step is because (\ref{eq:33}) of Assumption~\ref{asmp:1p} and the fact that the $\tilde\bvtheta_{\plt}$ is independent of the data imply that $Nn^{-1}\Exp\big\{\pi_{N}(\x,y;\tilde\bvtheta_{\plt})\|\ddot\ell(\btheta;\x,y)\|^2
    \bigm|\tilde\bvtheta_{\plt}\big\}=\Op$.
Thus, $\Delta_2\cvp\0$.

The partial derivative of $f(y_i\mid\x_i,\delta_i=1,\tilde\bvtheta_{\plt})$ with respect to $\btheta$ is
\begin{align*}
  &\dot f(y_i\mid\x_i,\delta_i=1,\tilde\bvtheta_{\plt})\notag\\
  &=\frac{\dot f(y_i\mid\x_i;\btheta)\pi_{N}(\x_i,y_i;\tilde\bvtheta_{\plt})}
  {\Exp\{\pi_{N}(\x_i,y_i;\tilde\bvtheta_{\plt})\mid\x_i,\tilde\bvtheta_{\plt}\}}\notag\\
  &\quad-\frac{f(y_i\mid\x_i;\btheta)\pi_{N}(\x_i,y_i;\tilde\bvtheta_{\plt})
    \Exp\{\dot\ell(\btheta;\x_i,y_i)\pi(\x_i,y;\tilde\bvtheta_{\plt})
    \mid\x_i,\tilde\bvtheta_{\plt}\}}
    {\Exp^2\{\pi_{N}(\x_i,y_i;\tilde\bvtheta_{\plt})
    \mid\x_i,\tilde\bvtheta_{\plt}\}}\notag\\
  &=\frac{\dot\ell(\btheta;\x_i,y_i)\pi_{N}(\x_i,y_i;\tilde\bvtheta_{\plt})
    f(y_i\mid\x_i;\btheta)}
    {\Exp\{\pi_{N}(\x_i,y_i;\tilde\bvtheta_{\plt})\mid\x_i,\tilde\bvtheta_{\plt}\}}\notag\\
  &\quad-\frac{f(y_i\mid\x_i;\btheta)\pi_{N}(\x_i,y_i;\tilde\bvtheta_{\plt})
    \Exp\{\dot\ell(\btheta;\x_i,y_i)\pi(x_i,y;\tilde\bvtheta_{\plt})
    \mid\x_i,\tilde\bvtheta_{\plt}\}}
    {\Exp^2\{\pi_{N}(\x_i,y_i;\tilde\bvtheta_{\plt})\mid\x_i,\tilde\bvtheta_{\plt}\}}.
\end{align*}
Thus, we have that
\begin{align*}
  \H_N^{\tilde\bvtheta_{\plt}}
  &=\onen\sumN\delta_i\bigg[\int\dot\ell(\btheta;\x_i,y_i)
     \dot f\tp(y\mid\x_i,\delta_i=1,\tilde\bvtheta_{\plt})\ud y\bigg]\\
  &=\onen\sumN\delta_i\bigg[
    \frac{\Exp\{\dot\ell^{\otimes2}(\btheta;\x_i,y_i)
    \pi_{N}(\x_i,y_i;\tilde\bvtheta_{\plt})\mid\x_i,\tilde\bvtheta_{\plt}\}}
    {\Exp\{\pi_{N}(\x_i,y_i;\tilde\bvtheta_{\plt})\mid\x_i,\tilde\bvtheta_{\plt}\}}\\
  &\hspace{3cm}  -\frac{\Exp^{\otimes2}\{\dot\ell(\btheta;\x_i,y_i)
    \pi(x_i,y;\tilde\bvtheta_{\plt})\mid\x_i,\tilde\bvtheta_{\plt}\}}
    {\Exp^2\{\pi_{N}(\x_i,y_i;\tilde\bvtheta_{\plt})\mid\x_i,\tilde\bvtheta_{\plt}\}}\bigg].
\end{align*}
Now we exam the limit of $\H_N^{\tilde\bvtheta_{\plt}}$. First, the expectation satisfies that
\begin{align*}
  &\Exp(\H_N^{\tilde\bvtheta_{\plt}}\mid\tilde\bvtheta_{\plt})\\
  &=\frac{N}{n}\Exp\Bigg[
    \dot\ell^{\otimes2}(\btheta;\x,y)\pi_{N}(\x,y;\tilde\bvtheta_{\plt})
    -\frac{\Exp^{\otimes2}\{\dot\ell(\btheta;\x,y)
    \pi_{N}(\x,y;\tilde\bvtheta_{\plt})\mid\x,\tilde\bvtheta_{\plt}\}}
    {\Exp\{\pi_{N}(\x,y;\tilde\bvtheta_{\plt})\mid\x,\tilde\bvtheta_{\plt}\}}
    \Biggm|\tilde\bvtheta_{\plt}\Bigg]
  =\bSigma_{N{\tilde\bvtheta_{\plt}}}.
\end{align*}
By the continuous mapping theorem, Assumption~\ref{asmp:3p}, and the fact that  $\tilde\bvtheta_{\plt}$ is independent of the data, we know that $\Exp(\H_N^{\tilde\bvtheta_{\plt}}\mid\tilde\bvtheta_{\plt})=\bSigma_{N{\tilde\bvtheta_{\plt}}}\cvp\bSigma_{\btheta,\bvtheta_p}$.

For the $j$-th diagonal element of $\H_{N,jj}$, using~(\ref{eq:44}), we have
\begin{align*}
  &\Var(\H_{N,jj}^{\tilde\bvtheta_{\plt}}\mid\tilde\bvtheta_{\plt})\\
  &=\frac{N}{n^2}\Var\Big[
    \delta\Exp\{\dot\ell_j^2(\btheta;\x,y)\mid\x,\delta=1,\tilde\bvtheta_{\plt}\}
    -\delta\Exp^2\{\dot\ell_j(\btheta;\x,y)\mid\x,\delta=1,\tilde\bvtheta_{\plt}\}
    \Bigm|\tilde\bvtheta_{\plt}\Big]\\
  &\le\frac{4N}{n^2}\Exp\Big[\delta\Exp^2\{\|\dot\ell(\btheta;\x,y)\|^2
    \mid\x,\delta=1,\tilde\bvtheta_{\plt}\}\Big]\\
  &\le\frac{4N}{n^2}\Exp\Big[\delta\Exp\{\|\dot\ell(\btheta;\x,y)\|^4
    \mid\x,\delta=1,\tilde\bvtheta_{\plt}\}\Big]\\
  &=\frac{4N}{n^2}\Exp\{\delta\|\dot\ell(\btheta;\x,y)\|^4\mid\tilde\bvtheta_{\plt}\}\\
  &=\frac{4N}{n^2}\Exp\{\pi_{N}(\x,y;\tilde\bvtheta_{\plt})\|\dot\ell(\btheta;\x,y)\|^4
    \mid\tilde\bvtheta_{\plt}\}=\op,
\end{align*}
where the last step is from (\ref{eq:49}). Thus, noting that $\H_N^{\tilde\bvtheta_{\plt}}\ge\0$ and $\bSigma_{N{\tilde\bvtheta_{\plt}}}=\bSigma_{\btheta,\bvtheta_p}+\op$, we know $\H_{N}^{\tilde\bvtheta_{\plt}}=\bSigma_{\btheta,\bvtheta_p}+\op$. Therefore,
\begin{equation}\label{eq:48}
  -\onen\ddot\ell_{S}(\btheta\mid\tilde\bvtheta_{\plt})=\bSigma_{\btheta,\bvtheta_p}+\op
\end{equation}

From (\ref{eq:46}), (\ref{eq:47}) and (\ref{eq:48}), applying the Basic Corollary in page 2 of \cite{hjort2011asymptotics} gives that for any constant vector $c$, 
\begin{equation*}
  \Pr\big\{\sqrt{n}\bSigma_{\btheta,\bvtheta_p}(\htheta_{s,\tilde\bvtheta_{\plt}}-\btheta)
  \le c \bigm|\tilde\bvtheta_{\plt}\big\} \cvp \Phi(c),
\end{equation*}
where $\Phi$ is the multivariate standard normal distribution function. Since a probability is bounded, this implies that
\begin{equation*}
  \Pr\big\{\sqrt{n}\bSigma_{\btheta,\bvtheta_p}(\htheta_{s,\tilde\bvtheta_{\plt}}-\btheta)
  \le c\} \cvp \Phi(c),
\end{equation*}
and this finishes the proof.

\subsection{Proof of Theorem~\ref{thm:2p}}

  Note that $\ell_{S}(\btheta;\x,y\mid\tilde\bvtheta_{\plt})=\log\{f(y\mid\x,\delta=1,\tilde\bvtheta_{\plt};\btheta)\}$. Taking partial derivatives of $\Exp\{\delta\u_{\tilde\bvtheta_{\plt}}(\btheta;\x,y)\mid\x,\tilde\bvtheta_{\plt}\}=\0$, we have
\begin{align*}
  \0=\Exp\{\delta\dot{\u}_{\tilde\bvtheta_{\plt}}(\btheta;\x,y)
  \mid\x,\tilde\bvtheta_{\plt}\}
  +\Exp\{\delta\u_{\tilde\bvtheta_{\plt}}(\btheta;\x,y)\dot\ell_{S}\tp(\btheta;\x,y\mid\tilde\bvtheta_{\plt})
  \mid\x,\tilde\bvtheta_{\plt}\},
\end{align*}
which implies
\begin{align*}
  \Exp\{\delta\u_{\tilde\bvtheta_{\plt}}(\btheta;\x,y)
  \dot\ell_{S}\tp(\btheta;\x,y\mid\tilde\bvtheta_{\plt})
  \mid\tilde\bvtheta_{\plt}\}
  &= -\Exp\{\delta\dot{\u}_{\tilde\bvtheta_{\plt}}(\btheta;\x,y)
  \mid\tilde\bvtheta_{\plt}\}.%
\end{align*}
Hence, by calculating the variance covariance matrix conditional on $\tilde\bvtheta_{\plt}$,
\begin{align*}
  &\Var\big\{\delta\M_{\btheta,\bvtheta_p}^{-1}
    \u_{\tilde\bvtheta_{\plt}}(\btheta;\x,y)
    +\delta\bSigma_{\btheta,\bvtheta_p}^{-1}
    \dot\ell_{S}(\btheta;\x,y\mid\tilde\bvtheta_{\plt})
    \bigm|\tilde\bvtheta_{\plt}\big\}\\
  &=\Var\big\{\delta\M_{\btheta,\bvtheta_p}^{-1}
    \u_{\tilde\bvtheta_{\plt}}(\btheta;\x,y)\bigm|\tilde\bvtheta_{\plt}\big\}
    +\M_{\btheta,\bvtheta_p}^{-1}\Exp\{\delta\u_{\tilde\bvtheta_{\plt}}(\btheta;\x,y)
    \dot\ell_{S}\tp(\btheta;\x,y)\mid\tilde\bvtheta_{\plt}\}
    \bSigma_{\btheta,\bvtheta_p}^{-1}\\
  &\quad +\bSigma_{\btheta,\bvtheta_p}^{-1} 
    \Exp\big\{\delta\dot\ell_{S}(\btheta;\x,y\mid\tilde\bvtheta_{\plt})
    \u_{\tilde\bvtheta_{\plt}}\tp(\btheta;\x,y)\bigm|\tilde\bvtheta_{\plt}\big\}
    (\M_{\btheta,\bvtheta_p}^{-1})\tp\\
  &\quad +\bSigma_{\btheta,\bvtheta_p}^{-1}\Var\{\delta\dot\ell_{S}(\btheta;\x,y
    \mid\tilde\bvtheta_{\plt})\mid\tilde\bvtheta_{\plt}\}\bSigma_{\btheta,\bvtheta_p}^{-1}\\
  &=nN^{-1}\Big\{\M_{\btheta,\bvtheta_p}^{-1} \V_{N,\btheta,\bvtheta} (\M_{\btheta,\bvtheta_p}^{-1})\tp
    -\M_{\btheta,\bvtheta_p}^{-1}\M_{N,\tilde\bvtheta_{\plt}}
    \bSigma_{\btheta,\bvtheta_p}^{-1}\\
  &\hspace{3cm}-\bSigma_{\btheta,\bvtheta_p}^{-1}\M_{N,\tilde\bvtheta_{\plt}}\tp
    (\M_{\btheta,\bvtheta_p}^{-1})\tp
    +\bSigma_{\btheta,\bvtheta_p}^{-1} \bSigma_{N\tilde\bvtheta_{\plt}}
    \bSigma_{\btheta,\bvtheta_p}^{-1}\Big\}\{1+\op\}\\
  &=nN^{-1}\big\{\M_{\btheta,\bvtheta_p}^{-1} \V_{\btheta,\bvtheta} (\M_{\btheta,\bvtheta_p}^{-1})\tp
    - \bSigma_{\btheta,\bvtheta_p}^{-1}\big\}\{1+\op\},
\end{align*}
where the last step is from the continuous mapping theorem. Thus, 
\begin{equation*}
  \M_{\btheta,\bvtheta_p}^{-1} \V_{\btheta,\bvtheta} (\M_{\btheta,\bvtheta_p}^{-1})\tp
    - \bSigma_{\btheta,\bvtheta_p}^{-1}
    =n^{-1}N\Var(\delta\M_{\btheta,\bvtheta_p}^{-1}T_{N\tilde\bvtheta_{\plt}}
    \mid\tilde\bvtheta_{\plt})\{1+\op\},
\end{equation*}
where
\begin{equation*}
T_{N\tilde\bvtheta_{\plt}}=
    \u_{\tilde\bvtheta_{\plt}}(\btheta;\x,y)
    +\M_{\btheta,\bvtheta_p}\bSigma_{\btheta,\bvtheta_p}^{-1}
    \dot\ell_{S}(\btheta;\x,y\mid\tilde\bvtheta_{\plt}).
\end{equation*}
Letting $N\rightarrow\infty$, we know that $\M_{\btheta,\bvtheta_p}^{-1} \V_{\btheta,\bvtheta} (\M_{\btheta,\bvtheta_p}^{-1})\tp\ge \bSigma_{\btheta,\bvtheta_p}^{-1}$ in the Loewner ordering, and the equality holds if $n^{-1}N\Exp(\delta\|T_{N\tilde\bvtheta_{\plt}}\|^2\mid\tilde\bvtheta_{\plt})=\op$ because
$\|\Var(\delta\M_{\btheta,\bvtheta_p}^{-1}T_{N\tilde\bvtheta_{\plt}}
\mid\tilde\bvtheta_{\plt})\|
\le\|\M_{\btheta,\bvtheta_p}^{-1}\|^2\|\Var(\delta T_{N\tilde\bvtheta_{\plt}}
\mid\tilde\bvtheta_{\plt})\|
\le\|\M_{\btheta,\bvtheta_p}^{-1}\|^2\|\Exp(\delta\|T_{N\tilde\bvtheta_{\plt}}\|^2
    \mid\tilde\bvtheta_{\plt})\|$.

\subsection{Technical details for examples}
\label{sec:techn-deta-exampl}

\subsubsection{Derivation of equation~(\ref{eq:14}) for multi-class logistic regression}
\label{sec:deriv-equat-refeq:14}

\noindent{\bf Derivation for the specific $\pi_{N}(\x_i,\y_i;\tilde\bvtheta_{\plt})$ in
  (\ref{eq:14})}: \\
The sampled data log-likelihood function for $\bbeta$ (up to an additive constant) is
\begin{equation*}
  \ell_{S}(\bbeta\mid\tilde\bvtheta_{\plt})=\sumN\delta_i\bigg[
  (\y_i\otimes\x_i)\tp\bbeta
  -\log\bigg\{\sum_{l=1}^K\exp(\x_i\tp\bbeta_l)\bigg\}
  -\log\Exp\{\|\y_i-\p(\x_i,\tbeta_{\plt})\|\mid\x_i,\tbeta_{\plt}\}\bigg].
\end{equation*}
By direct calculation, we know that the score function is  
\begin{equation*}
  \dot\ell_{S}(\bbeta\mid\tilde\bvtheta_{\plt})
  =\sumN\delta_i\bigg[\y_i
    -\frac{\Exp\{\y_i\|\y_i-\p(\x_i,\tbeta_{\plt})\|\mid\x_i,\tbeta_{\plt}\}}
    {\Exp\{\|\y_i-\p(\x_i,\tbeta_{\plt})\|\mid\x_i,\tbeta_{\plt}\}}\bigg]\otimes\x_i.
\end{equation*}
For $k=1,...,K$, when $y_{i,k}=1$, 
\begin{align*}
  \|\y_i-\p(\x_i,\tbeta_{\plt})\|^2
  &=\{1-p_{k}(\x_i,\tbeta_{\plt})\}^2+\sum_{l\neq k}p_{l}^2(\x_i,\tbeta_{\plt})\\
  &=1-2p_{k}(\x_i,\tbeta_{\plt})+\sum_{l=1}^Kp_{l}^2(\x_i,\tbeta_{\plt})
  =\exp(2\tilde{g}_{i,k})\bigg\{\sum_{l=1}^Kp_{l}^2(\x_i,\tbeta_{\plt})\bigg\}.
\end{align*}
Using the above expression, we obtain the following expectations.
\begin{align*}
  \Exp\{y_{i,k}\|\y_i-\p(\x_i,\tbeta_{\plt})\|\mid\x_i,\tbeta_{\plt}\}
    &=p_{k}(\x_i,\bbeta)\exp(\tilde{g}_{i,k})
      \bigg\{\sum_{l=1}^Kp_{l}^2(\x_i,\tbeta_{\plt})\bigg\}^{1/2},
  \quad\text{and}\\
  \Exp\{\|\y_i-\p(\x_i,\tbeta_{\plt})\|\mid\x_i,\tbeta_{\plt}\}
    &=\bigg\{\sum_{k=1}^Kp_{k}(\x_i,\bbeta)\exp(\tilde{g}_{i,k})\bigg\}
      \bigg\{\sum_{l=1}^Kp_{l}^2(\x_i,\tbeta_{\plt})\bigg\}^{1/2}.
\end{align*}
Thus, the ratio is 
\begin{equation*}
  \frac{\Exp\{y_{i,k}\|\y_i-\p(\x_i,\tbeta_{\plt})\|\mid\x_i\}}
    {\Exp\{\|\y_i-\p(\x_i,\tbeta_{\plt})\|\mid\x_i\}}
  =\frac{p_{k}(\x_i,\bbeta)\exp(\tilde{g}_{i,k})}{\sum_{k=1}^Kp_{k}(\x_i,\bbeta)\exp(\tilde{g}_{i,k})}
  =\frac{\exp(\x_i\tp\bbeta_k+\tilde{g}_{i,k})}
  {\sum_{l=1}^K\exp(\x_i\tp\bbeta_l+g_{i,l})}=\tilde{p}_{i,k}^g.
\end{equation*}

\noindent{\bf Derivation for general $\pi_{N}(\x_i,\y_i;\tilde\bvtheta_{\plt})$}: \\
From the facts of $\bar{\pi}_{N}(\x_i;\btheta\mid
\tilde\bvtheta_{\plt})=\sum_{l=1}^K\pi_{N}(\x_i,\1_l;\tilde\bvtheta_{\plt})p_{l}(\x_i,\bbeta)$
and $\partial p_{l}(\x_i,\btheta)/\partial\btheta = p_{l}(\x_i,\btheta)\{\1-\p(\x_i,\bbeta))\otimes\x_i$, we know that
\begin{align}
  \frac{\partial{\log\bar{\pi}}_{N}(\x_i;\btheta\mid \tilde\bvtheta_{\plt})}
  {\partial\btheta}
  &=\frac{\sum_{l=1}^K\pi_{N}(\x_i,\1_l;\tilde\bvtheta_{\plt})
    \partial p_{l}(\x_i,\btheta)/\partial\btheta}
    {\sum_{l=1}^K\pi_{N}(\x_i,\1_l;\tilde\bvtheta_{\plt})p_{l}(\x_i,\bbeta)}\\
  &=\frac{\sum_{l=1}^K\pi_{N}(\x_i,\1_l;\tilde\bvtheta_{\plt})
    p_{l}(\x_i,\btheta)\{\1_l-\p(\x_i,\bbeta))\otimes\x_i}
    {\sum_{l=1}^K\pi_{N}(\x_i,\1_l;\tilde\bvtheta_{\plt})p_{l}(\x_i,\bbeta)}\\
  &=\frac{\sum_{l=1}^K\pi_{N}(\x_i,\1_l;\tilde\bvtheta_{\plt})p_{l}(\x_i,\btheta)\1_l}
    {\sum_{l=1}^K\pi_{N}(\x_i,\1_l;\tilde\bvtheta_{\plt})p_{l}(\x_i,\bbeta)}\otimes\x_i
    - \p(\x_i,\bbeta)\otimes\x_i.
\end{align}
Note that with $\tilde{g}_{i,k}=\log\{\pi_{N}(\x_i,\1_k;\tilde\bvtheta_{\plt})\}$, the $k$-th element of
\begin{equation}
\frac{\sum_{l=1}^K\pi_{N}(\x_i,\1_l;\tilde\bvtheta_{\plt})p_{l}(\x_i,\btheta)\1_l}
    {\sum_{l=1}^K\pi_{N}(\x_i,\1_l;\tilde\bvtheta_{\plt})p_{l}(\x_i,\bbeta)}
\end{equation}
is
\begin{equation}
  \frac{\pi_{N}(\x_i,\1_k;\tilde\bvtheta_{\plt})p_{k}(\x_i,\btheta)}
  {\sum_{l=1}^K\pi_{N}(\x_i,\1_l;\tilde\bvtheta_{\plt})p_{l}(\x_i,\bbeta)}
  =\frac{\pi_{N}(\x_i,\1_k;\tilde\bvtheta_{\plt})e^{\x_i\tp\bbeta_k}}
  {\sum_{l=1}^K\pi_{N}(\x_i,\1_l;\tilde\bvtheta_{\plt})e^{\x_i\tp\bbeta_l}}
  =\frac{e^{\x_i\tp\bbeta_k+\tilde{g}_{i,k}}}
  {\sum_{l=1}^Ke^{\x_i\tp\bbeta_l+\tilde{g}_{i,l}}}
  =\tilde{p}_{i,k}^g.
\end{equation}
Thus the specific expression of (\ref{eq:70}) gives
\begin{align}
  \dot\ell(\btheta;\x_i,y_i)
  - \frac{\partial{\log\bar{\pi}}_{N}(\x_i;\btheta\mid \tilde\bvtheta_{\plt})}
  {\partial\btheta}
  =\{\y_i-\p_i^g(\bbeta,\tbeta_{\plt})\}\otimes\x_i,
\end{align}
which finishes the proof.

\subsubsection{Derivations of equations~(\ref{eq:20}), (\ref{eq:21}), and (\ref{eq:22}) for Poisson regression}
\label{sec:deriv-equat-refeq:20}

For non-negative integers $m$ and $k$, denote
\begin{equation*}
q(m,k)=\sum_{y=0}^m\frac{y^ke^{-\mu}\mu^{y}}{y!}.  
\end{equation*}
We have
\begin{align*}
  q(m,k):=&\sum_{y=0}^m\frac{e^{-\mu}\mu^{y}}{y!}y^k
  =\mu\sum_{y=1}^m\frac{e^{-\mu}\mu^{y-1}}{(y-1)!}y^{k-1}\\
  =&\mu\sum_{y=0}^{m-1}\frac{e^{-\mu}\mu^{y}}{y!}(y+1)^{k-1}
  =\sum_{l=0}^{k-1}\binom{k-1}{l}
     \mu\sum_{y=0}^{m-1}\frac{e^{-\mu}\mu^{y}}{y!}y^{l}\\
  =&\mu\sum_{l=0}^{k-1}\binom{k-1}{l}q(m-1,l).
\end{align*}
Thus, we know that
\begin{align*}
  q(m,0)&=F(m;\mu),\\
  q(m,1)&%
          =\mu F(m-1;\mu),\\
  q(m,2)&%
          =\mu F(m-1;\mu)+\mu^2F(m-2;\mu),\\
  q(m,3)&%
  =\mu F(m-1;\mu)+3\mu^2F(m-2;\mu)+\mu^3F(m-3;\mu).
\end{align*}

Now, we derive the expectations. First,
\begin{align*}
  \Exp(|y_i-\tmu_i|\mid\x_i;\tbeta_{\plt})
  &=\sum_{y_i=0}^\infty\frac{e^{-\mu_i}\mu_i^{y_i}}{y_i!}|y_i-\tmu_i|
   =2\sum_{y_i=0}^{m_i}\frac{e^{-\mu_i}\mu_i^{y_i}}{y_i!}(\tmu_i-y_i)
    +\mu_i-\tmu_i\\
  &=2\tmu_iq(m_i,0)
    -2q(m_i,1)+\mu_i-\tmu_i\\
  &=2\tmu_iF(m_i;\mu_i)-2\mu_iF(m_i-1;\mu_i)+\mu_i-\tmu_i\\
  &=2(\tmu_i-\mu_i)F(m_i-1;\mu_i)
    +2\tmu_if(m_i;\mu_i)+\mu_i-\tmu_i.
\end{align*}
Second,
\begin{align*}
  &\Exp(y_i|y_i-\tmu_i|\mid\x_i;\tbeta_{\plt})
  =\sum_{y_i=0}^\infty\frac{e^{-\mu_i}\mu_i^{y_i}}{y_i!}(y_i|y_i-\tmu_i|)\\
  &=2\sum_{y_i=0}^{m_i}\frac{e^{-\mu_i}\mu_i^{y_i}}{y_i!}
    (\tmu_i y_i-y_i^2)+\Exp(y_i^2-y_i\tmu_i\mid\x_i)\\
  &=2\tmu_iq(m_i,1)-2q(m_i,2)+\mu_i+\mu_i^2-\mu_i\tmu_i\\
  &=2\tmu_i\mu_iF(m_i-1;\mu_i)-2\{\mu_iF(m_i-1;\mu_i)+\mu_i^2F(m_i-2;\mu_i)\}
    +\mu_i+\mu_i^2-\mu_i\tmu_i\\
  &=2\mu_i(\tmu_i-1)F(m_i-1;\mu_i)-2\mu_i^2F(m_i-2;\mu_i)
    +\mu_i+\mu_i^2-\mu_i\tmu_i.%
\end{align*}
Third,
\begin{align*}
  &\Exp(y_i^2|y_i-\tmu_i|\mid\x_i;\tbeta_{\plt})
  =\sum_{y_i=0}^\infty\frac{e^{-\mu_i}\mu_i^{y_i}}{y_i!}(y_i^2|y_i-\tmu_i|)\\
  &=2\sum_{y_i=0}^{m_i}\frac{e^{-\mu_i}\mu_i^{y_i}}{y_i!}
    (\tmu_i y_i^2-y_i^3)+\Exp(y_i^3-y_i^2\tmu_i\mid\x_i)\\
  &=2\tmu_iq(m_i,2)-2q(m_i,3)
    +\mu_i+3\mu_i^2+\mu_i^3-\tmu_i(\mu_i+\mu_i^2)\\
  &=2\tmu_i\{\mu_iF(m_i-1;\mu_i)+\mu_i^2F(m_i-2;\mu_i)\}\\
  &\quad-2\{\mu_iF(m_i-1;\mu_i)+3\mu_i^2F(m_i-2;\mu_i)+\mu_i^3F(m_i-3;\mu_i)\}\\
  &\quad+\mu_i+3\mu_i^2+\mu_i^3-\tmu_i(\mu_i+\mu_i^2)\\
  &=2\mu_i(\tmu_i-1)F(m_i-1;\mu_i)+2\mu_i^2(\tmu_i-3)F(m_i-2;\mu_i)
    -2\mu_i^3F(m_i-3;\mu_i)\\
  &\quad+\mu_i+3\mu_i^2+\mu_i^3-\tmu_i(\mu_i+\mu_i^2).
\end{align*}

\section{Additional numerical experiments on pilot misspecifications}
\label{sec:addit-numer-exper}
We carried out additional numerical experiments to investigate the effect of the
magnitude of pilot misspecification on different subsample estimators.

\subsection{multi-class logistic regression}
\label{sec:multi-class-logistic-1}
We used exactly the same setup as in Section~\ref{sec:multi-class-logistic} with
$n_0=400$ and $n=2,000$ in this experiment. The only difference is that the
pilot is set to be the true parameter plus $\lambda$ times a vector of
ones. With this setting, $\lambda$ controls the level of pilot misspecification
with $\lambda=0$ corresponding to the consistent pilot. We used the same
approach used in Section~\ref{sec:multi-class-logistic} to calculate the
empirical MSEs, variances, and squared biases. The relative pattern between the
variances and the squared biases are the same as in
Section~\ref{sec:multi-class-logistic}, so we report the empirical MSEs here
only. We omit the results for the native estimator due to its high bias.

\begin{figure}[htp]%
  \centering
  \begin{subfigure}{0.485\textwidth}
    \includegraphics[width=\textwidth,page=1]{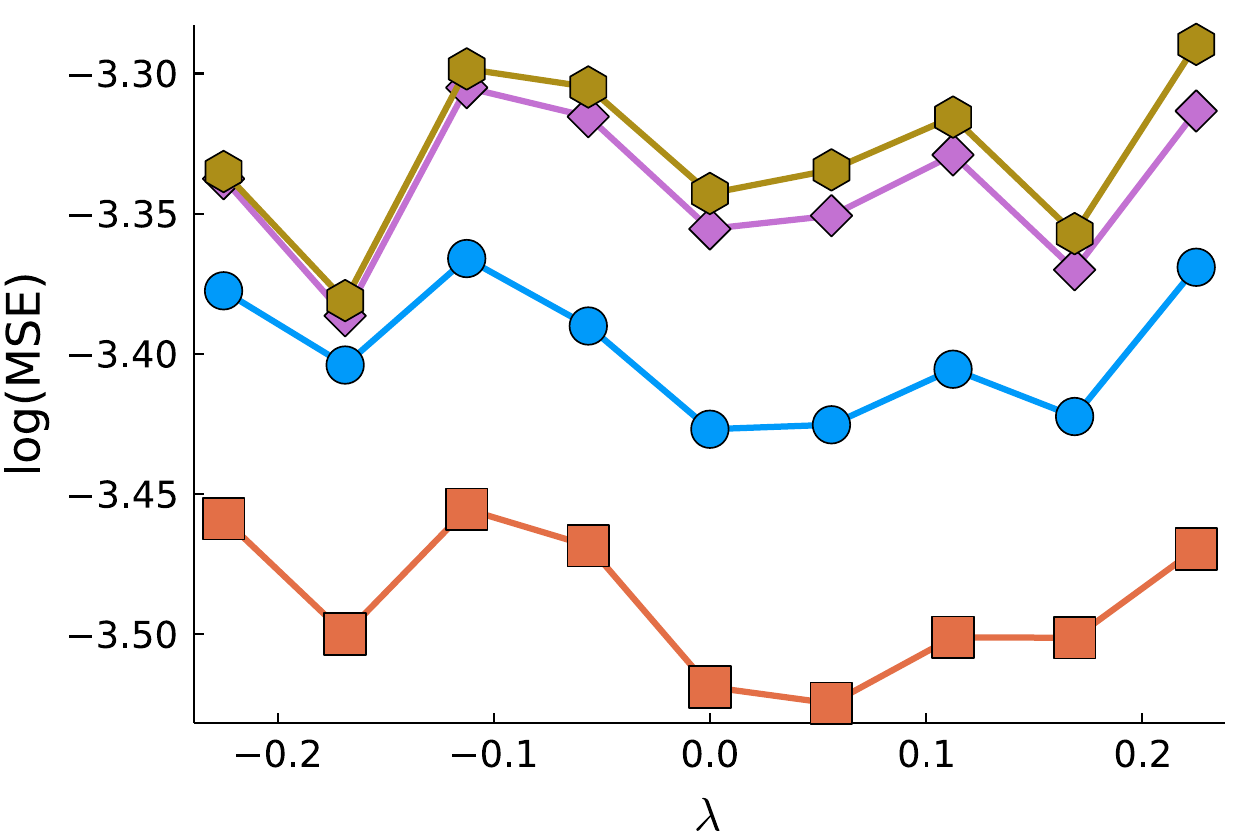}
    \caption{$\x_{-1,i}$'s are normal}
  \end{subfigure}
  \begin{subfigure}{0.485\textwidth}
    \includegraphics[width=\textwidth,page=2]{figures/00mseMultiLogiMisPltDiff.pdf}
    \caption{$\x_{-1,i}$'s are lognormal}
  \end{subfigure}
  \begin{subfigure}{0.485\textwidth}
    \includegraphics[width=\textwidth,page=3]{figures/00mseMultiLogiMisPltDiff.pdf}
    \caption{$\x_{-1,i}$'s are $\mathbb{T}_3$}
  \end{subfigure}
  \begin{subfigure}{0.485\textwidth}
    \includegraphics[width=\textwidth,page=4]{figures/00mseMultiLogiMisPltDiff.pdf}
    \caption{$\x_{-1,i}$'s are exponential}
  \end{subfigure}
  \caption{Log of empirical MSEs of subsample estimators with $n_0=400$
    and $n=2,000$ in multi-class logistic regression when the misspecified pilot
    estimator is set to be the true parameter plus a vector of $\lambda$'s.}
  \label{fig:logistic5}
\end{figure}

Figure~\ref{fig:logistic5} reports the results. The MSCLE dominates the
corresponding IPW estimator for both GN and LUS probabilities uniformly in
$\lambda$. For MSCLE, it has good performance with a consistent pilot
($\lambda=0$), but its best performance may not always be achieved in this
case. The reason is that $\Sigma_{\btheta,\bvtheta_p}^{-1}$ in
Theorem~\ref{thm:1p} may not necessarily be minimized at the true parameter.
Actually, for the MSCLE, there is no general
solution to $\bvtheta_p$ so that $\Sigma_{\btheta,\bvtheta_p}^{-1}$ is
minimized.

Note that $\Sigma_{\btheta,\bvtheta_p}^{-1}$ depends on $\bvtheta_p$ through
$\pi_{N}(\x,y;\bvtheta_p)$ so the problem of find the optimal $\bvtheta_p$ is
essentially finding the optimal sampling probability $\pi_{N}(\x,y;\bvtheta_p)$
for the MSCLE. The problem is complicated even with a much simplified scenario
of noninformative subsampling so that $\pi_{N}(\x,y;\bvtheta_p)$ does not dependent
on $y$, i.e., $\pi_{N}(\x,y;\bvtheta_p)=\pi_{N}(\x;\bvtheta_p)$.  In this
scenario, the problem of determining the optimal
$\pi_{N}(\x;\bvtheta_p)$ is called optimal design of experiments \citep[see,
e.g., ][]{kiefer1959,pukelsheim2006optimal}, and the optimal
$\pi_{N}(\x;\bvtheta_p)$ is binary with possible values of zero and one
\citep{PronzatoWang2020}. This topic is beyond the scope our investigation
because the primary focus of this study is to propose an improved estimator for
informative subsamples.

For the IPW estimator, we see that the optimal performance is not achieved with $\lambda=0$ for
cases (b) and (d). This is because the GM sampling probabilities
in~(\ref{eq:15}) minimize the variance of a certain linear function of the IPW
estimator which is different from the MSE. Nevertheless, for the IPW estimator,
the optimal sampling probabilities exist under some optimality criteria such as
the A- and L- optimality, and it requires the pilot to be consistent to achieve
the optimal variance. The proposed MSCLE has a smaller variance matrix than the
weighted estimator. Thus, we can only conclude that the asymptotic variance of
the MSCLE has a smaller upper bound with a consistent pilot than with a
misspecified pilot.

\subsection{Poisson regression}
\label{sec:poisson-regression-1}

We used exactly the same setup and procedure as in
Section~\ref{sec:poisson-regression-2} to generate data and calculate the
empirical MSEs, variances, and squared biases for $n_0=400$ and $n=2,000$. The
only difference is that the pilot is set to be the true parameter plus $\lambda$
times a vector of ones, so that $\lambda$ controls the level of pilot
misspecification with $\lambda=0$ corresponding to the consistent pilot. Again,
we report the empirical MSEs here only and omit the results for the native
estimator.

\begin{figure}[htp]%
  \centering
  \begin{subfigure}{0.485\textwidth}
    \includegraphics[width=\textwidth,page=1]{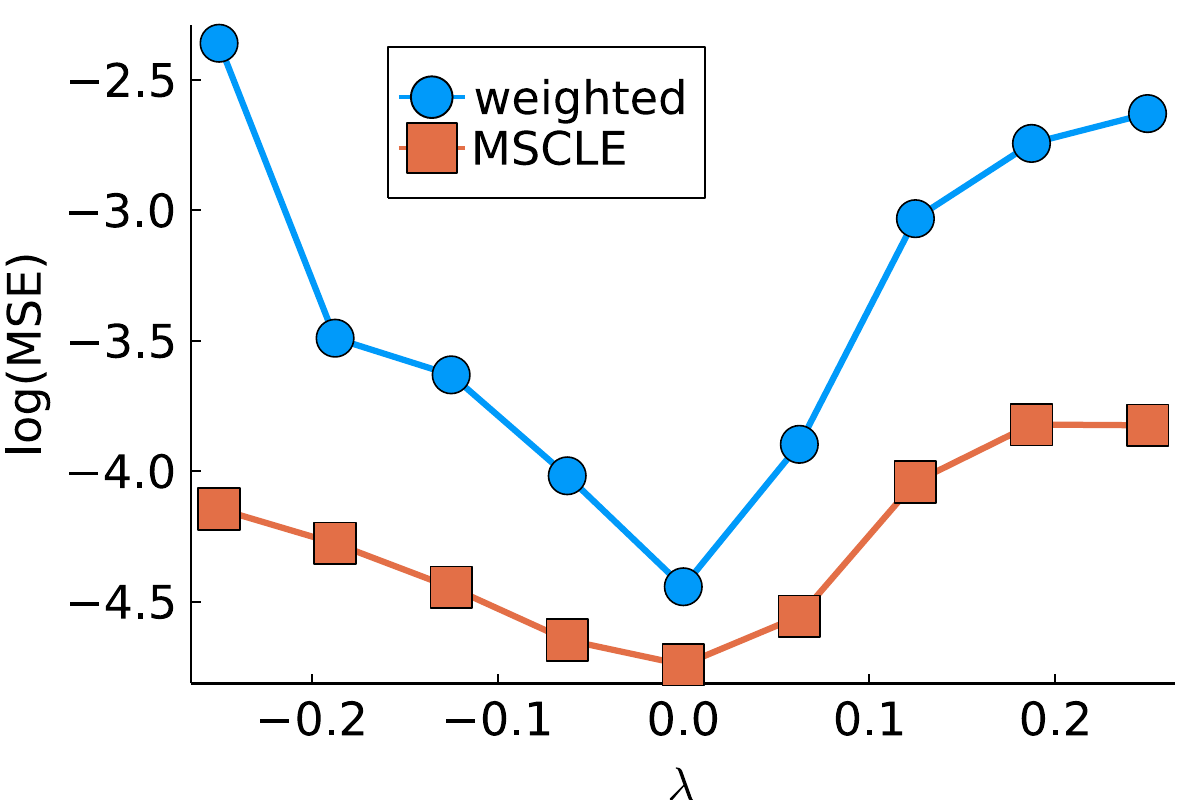}
    \caption{$\x_{-1,i}$'s are uniform}
  \end{subfigure}
  \begin{subfigure}{0.485\textwidth}
    \includegraphics[width=\textwidth,page=2]{figures/00msePoissonMisPltDiff.pdf}
    \caption{$\x_{-1,i}$'s are Beta}
  \end{subfigure}
  \begin{subfigure}{0.485\textwidth}
    \includegraphics[width=\textwidth,page=3]{figures/00msePoissonMisPltDiff.pdf}
    \caption{$\x_{-1,i}$'s are normal}
  \end{subfigure}
  \begin{subfigure}{0.485\textwidth}
    \includegraphics[width=\textwidth,page=4]{figures/00msePoissonMisPltDiff.pdf}
    \caption{$\x_{-1,i}$'s are exponential}
  \end{subfigure}
  \caption{Log of empirical MSEs of subsample estimators with $n_0=400$
    and $n=2,000$ in Poisson regression when the misspecified pilot
    estimator is set to be the true parameter plus a vector of $\lambda$'s.}
  \label{fig:poisson5}
\end{figure}

Figure~\ref{fig:poisson5} reports the results. The MSCLE dominates the
corresponding IPW estimator uniformly, and both methods achieve the best
performance for all the four cases with a consistent pilot ($\lambda=0$) in this
example.

We also consider the case of misspecified pilot estimator using the same
procedure considered in Section~\ref{sec:multi-class-logistic}. The resulting
impacts are similar to those observed in Section~\ref{sec:multi-class-logistic} so we omit the results.

\bibliography{ref}

\begin{thebibliography}{50}
\providecommand{\natexlab}[1]{#1}
\providecommand{\url}[1]{\texttt{#1}}
\expandafter\ifx\csname urlstyle\endcsname\relax
  \providecommand{\doi}[1]{doi: #1}\else
  \providecommand{\doi}{doi: \begingroup \urlstyle{rm}\Url}\fi

\bibitem[Ai et~al.(2021)Ai, Yu, Zhang, and Wang]{ai2020optimal}
Mingyao Ai, Jun Yu, Huiming Zhang, and HaiYing Wang.
\newblock Optimal subsampling algorithms for big data regressions.
\newblock \emph{Statistica Sinica}, 31\penalty0 (2):\penalty0 749--772, 2021.
\newblock \doi{10.5705/ss.202018.0439}.

\bibitem[Bezanson et~al.(2017)Bezanson, Edelman, Karpinski, and
  Shah]{bezanson2017julia}
Jeff Bezanson, Alan Edelman, Stefan Karpinski, and Viral~B Shah.
\newblock Julia: A fresh approach to numerical computing.
\newblock \emph{SIAM review}, 59\penalty0 (1):\penalty0 65--98, 2017.
\newblock URL \url{https://doi.org/10.1137/141000671}.

\bibitem[Blackard and Dean(1999)]{blackard1999comparative}
Jock~A Blackard and Denis~J Dean.
\newblock Comparative accuracies of artificial neural networks and discriminant
  analysis in predicting forest cover types from cartographic variables.
\newblock \emph{Computers and electronics in agriculture}, 24\penalty0
  (3):\penalty0 131--151, 1999.

\bibitem[Breslow and Holubkov(1997)]{bre97}
N.E. Breslow and R.~Holubkov.
\newblock Maximum likelihood estimation of logistic regression parameters under
  two-phase, outcome-dependent sampling.
\newblock \emph{Journal of the Royal Statistical Society: Series B},
  59:\penalty0 447--461, 1997.

\bibitem[Chambers and Skinner(2003)]{cha03}
R.~L. Chambers and C.~J. Skinner, editors.
\newblock \emph{Analysis of Survey Data}.
\newblock John Wiley {\&} Sons, Chichester, U.~K., 2003.

\bibitem[Cox(1969)]{cox69}
D.~R. Cox.
\newblock Some sampling problems in technology.
\newblock In U.L. Johnson and H.~Smith, editors, \emph{New Developments in
  Survey Sampling}, pages 506--527. Wiley Interscience, 1969.

\bibitem[Dhillon et~al.(2013)Dhillon, Lu, Foster, and Ungar]{dhillon2013new}
Paramveer Dhillon, Yichao Lu, Dean~P Foster, and Lyle Ungar.
\newblock New subsampling algorithms for fast least squares regression.
\newblock In \emph{Advances in Neural Information Processing Systems}, pages
  360--368, 2013.

\bibitem[Drineas et~al.(2012)Drineas, Magdon-Ismail, Mahoney, and
  Woodruff]{Drineas:12}
P.~Drineas, M.~Magdon-Ismail, M.W. Mahoney, and D.P. Woodruff.
\newblock Faster approximation of matrix coherence and statistical leverage.
\newblock \emph{Journal of Machine Learning Research}, 13:\penalty0 3475--3506,
  2012.

\bibitem[Drineas et~al.(2006)Drineas, Mahoney, and
  Muthukrishnan]{drineas2006sampling}
Petros Drineas, Michael~W Mahoney, and S~Muthukrishnan.
\newblock Sampling algorithms for $l_2$ regression and applications.
\newblock In \emph{Proceedings of the seventeenth annual ACM-SIAM symposium on
  Discrete algorithm}, pages 1127--1136. Society for Industrial and Applied
  Mathematics, 2006.

\bibitem[Efraimidis and Spirakis(2006)]{efraimidis2006}
Pavlos~S. Efraimidis and Paul~G. Spirakis.
\newblock Weighted random sampling with a reservoir.
\newblock \emph{Information Processing Letters}, 97\penalty0 (5):\penalty0
  181--185, 2006.

\bibitem[Fithian and Hastie(2014)]{fithian2014local}
William Fithian and Trevor Hastie.
\newblock Local case-control sampling: Efficient subsampling in imbalanced data
  sets.
\newblock \emph{Annals of statistics}, 42\penalty0 (5):\penalty0 1693, 2014.

\bibitem[Han et~al.(2020)Han, Tan, Yang, and Zhang]{Han2019}
Lei Han, Kean~Ming Tan, Ting Yang, and Tong Zhang.
\newblock Local uncertainty sampling for large-scale multiclass logistic
  regression.
\newblock \emph{The Annals of Statistics}, 48\penalty0 (3):\penalty0
  1770--1788, 2020.

\bibitem[Hesterberg(1995)]{hesterberg1995weighted}
Tim Hesterberg.
\newblock Weighted average importance sampling and defensive mixture
  distributions.
\newblock \emph{Technometrics}, 37\penalty0 (2):\penalty0 185--194, 1995.

\bibitem[Hjort and Pollard(2011)]{hjort2011asymptotics}
Nils~Lid Hjort and David Pollard.
\newblock Asymptotics for minimisers of convex processes.
\newblock \emph{arXiv preprint arXiv:1107.3806}, 2011.

\bibitem[Hsieh et~al.(1985)Hsieh, Manski, and McFadden]{hsi85}
D.A. Hsieh, C.F. Manski, and D.~McFadden.
\newblock Estimation of response probabilities from augmented retrospective
  observations.
\newblock \emph{Journal of the American Statistical Association}, 80:\penalty0
  651--662, 1985.

\bibitem[Hu and Lawless(1996)]{hu96}
X.J. Hu and J.F. Lawless.
\newblock Estimation from truncated lifetime data with supplementary
  information on covariates and censoring times.
\newblock \emph{Biometrika}, 83:\penalty0 651--662, 1996.

\bibitem[Hu and Lawless(1997)]{hu97}
X.J. Hu and J.F. Lawless.
\newblock Pseudo likelihood estimation in a class of problems with
  response-related missing covariates.
\newblock \emph{Canadian Journal of Statistics}, 25:\penalty0 125--142, 1997.

\bibitem[Innes(2018)]{innes:2018}
Mike Innes.
\newblock Flux: Elegant machine learning with julia.
\newblock \emph{Journal of Open Source Software}, 2018.
\newblock \doi{10.21105/joss.00602}.

\bibitem[Kalbfleisch and Lawless(1988)]{kal88}
J.D. Kalbfleisch and J.F. Lawless.
\newblock Likelihood analysis of multi-state models for disease incidence and
  mortality.
\newblock \emph{Statistics in Medicine}, 7:\penalty0 149--160, 1988.

\bibitem[Kiefer(1959)]{kiefer1959}
J.~Kiefer.
\newblock Optimum experimental designs.
\newblock \emph{Journal of the Royal Statistical Society. Series B},
  21\penalty0 (2):\penalty0 272--319, 1959.

\bibitem[Kim et~al.(2006)Kim, Navarro, and Fuller]{kim2006}
J.K. Kim, A.~Navarro, and W.A. Fuller.
\newblock Replicate variance estimation after multi-phase stratified sampling.
\newblock \emph{Journal of the American Statistical Association}, 101:\penalty0
  312--320, 2006.

\bibitem[LeCun et~al.(1998)LeCun, Bottou, Bengio, and
  Haffner]{lecun1998gradient}
Yann LeCun, L{\'e}on Bottou, Yoshua Bengio, and Patrick Haffner.
\newblock Gradient-based learning applied to document recognition.
\newblock \emph{Proceedings of the IEEE}, 86\penalty0 (11):\penalty0
  2278--2324, 1998.

\bibitem[Ma et~al.(2015)Ma, Mahoney, and Yu]{PingMa2014-JMLR}
P.~Ma, M.W. Mahoney, and B.~Yu.
\newblock A statistical perspective on algorithmic leveraging.
\newblock \emph{Journal of Machine Learning Research}, 16:\penalty0 861--911,
  2015.

\bibitem[McWilliams et~al.(2014)McWilliams, Krummenacher, Lucic, and
  Buhmann]{mcwilliams2014fast}
Brian McWilliams, Gabriel Krummenacher, Mario Lucic, and Joachim~M Buhmann.
\newblock Fast and robust least squares estimation in corrupted linear models.
\newblock In \emph{Advances in Neural Information Processing Systems}, pages
  415--423, 2014.

\bibitem[Nie et~al.(2018)Nie, Wiens, and Zhai]{nie2018minimax}
Rui Nie, Douglas~P Wiens, and Zhichun Zhai.
\newblock Minimax robust active learning for approximately specified regression
  models.
\newblock \emph{Canadian Journal of Statistics}, 46\penalty0 (1):\penalty0
  104--122, 2018.

\bibitem[Owen and Zhou(2000)]{owen2000safe}
Art Owen and Yi~Zhou.
\newblock Safe and effective importance sampling.
\newblock \emph{Journal of the American Statistical Association}, 95\penalty0
  (449):\penalty0 135--143, 2000.

\bibitem[Pfeffermann et~al.(1998)Pfeffermann, Krieger, and Rinnot]{pfe98}
D.~Pfeffermann, A.~M. Krieger, and Y.~Rinnot.
\newblock Parametric distributions of complex survey data under informative
  probability sampling.
\newblock \emph{Statistica Sinica}, 8:\penalty0 1087--1114, 1998.

\bibitem[Pronzato and Wang(2021)]{PronzatoWang2020}
Luc Pronzato and HaiYing Wang.
\newblock Sequential online subsampling for thinning experimental designs.
\newblock \emph{Journal of Statistical Planning and Inference}, 212:\penalty0
  169 -- 193, 2021.
\newblock ISSN 0378-3758.
\newblock \doi{10.1016/j.jspi.2020.08.001}.
\newblock URL
  \url{http://www.sciencedirect.com/science/article/pii/S0378375820300999}.

\bibitem[Pukelsheim(2006)]{pukelsheim2006optimal}
Friedrich Pukelsheim.
\newblock \emph{Optimal design of experiments}.
\newblock SIAM, 2006.

\bibitem[Qin(2017)]{qin2017biased}
Jing Qin.
\newblock \emph{Biased sampling, over-identified parameter problems and
  beyond}.
\newblock Springer, 2017.

\bibitem[Rubin(1976)]{rubin1976}
D.~B. Rubin.
\newblock Inference and missing data.
\newblock \emph{Biometrika}, 63\penalty0 (3):\penalty0 581--592, 1976.

\bibitem[Saegusa and Wellner(2013)]{saegusa2013}
T.~Saegusa and J.~A. Wellner.
\newblock Weighted likelihood estimation under two-phase sampling.
\newblock \emph{The Annals of Statistics}, 41:\penalty0 269--295, 2013.

\bibitem[Scott and Wild(1986)]{scott1986fitting}
A.~J. Scott and C.~J. Wild.
\newblock Fitting logistic models under case-control or choice based sampling.
\newblock \emph{Journal of the Royal Statistical Society. Series B},
  48\penalty0 (2):\penalty0 170--182, 1986.

\bibitem[Scott and Wild(1991)]{sco91}
A.J. Scott and C.J. Wild.
\newblock Fitting logistic regression models in stratified case-control
  studies.
\newblock \emph{Biometrics}, 47:\penalty0 497--510, 1991.

\bibitem[Scott and Wild(1997)]{sco97}
A.J. Scott and C.J. Wild.
\newblock Fitting regression models to case-control data by maximum likelihood.
\newblock \emph{Biometrika}, 84:\penalty0 57--71, 1997.

\bibitem[Shen et~al.(2021)Shen, Chen, and Yu]{shen2021surprise}
Xinwei Shen, Kani Chen, and Wen Yu.
\newblock Surprise sampling: Improving and extending the local case-control
  sampling.
\newblock \emph{Electronic Journal of Statistics}, 15\penalty0 (1):\penalty0
  2454--2482, 2021.

\bibitem[Ting and Brochu(2018)]{influence2018}
Daniel Ting and Eric Brochu.
\newblock Optimal subsampling with influence functions.
\newblock In S.~Bengio, H.~Wallach, H.~Larochelle, K.~Grauman, N.~Cesa-Bianchi,
  and R.~Garnett, editors, \emph{Advances in Neural Information Processing
  Systems}, volume~31. Curran Associates, Inc., 2018.
\newblock URL
  \url{https://proceedings.neurips.cc/paper/2018/file/57c0531e13f40b91b3b0f1a30b529a1d-Paper.pdf}.

\bibitem[Wang(2019)]{wang2019more}
HaiYing Wang.
\newblock More efficient estimation for logistic regression with optimal
  subsamples.
\newblock \emph{Journal of Machine Learning Research}, 20\penalty0
  (132):\penalty0 1--59, 2019.

\bibitem[Wang(2020)]{Wang2020RareICML}
HaiYing Wang.
\newblock Logistic regression for massive data with rare events.
\newblock In Hal~Daumé III and Aarti Singh, editors, \emph{Proceedings of the
  37th International Conference on Machine Learning}, volume 119 of
  \emph{Proceedings of Machine Learning Research}, pages 9829--9836. PMLR,
  13--18 Jul 2020.
\newblock URL \url{http://proceedings.mlr.press/v119/wang20a.html}.

\bibitem[Wang and Ma(2021)]{wang2021optimal}
Haiying Wang and Yanyuan Ma.
\newblock Optimal subsampling for quantile regression in big data.
\newblock \emph{Biometrika}, 108\penalty0 (1):\penalty0 99--112, 2021.

\bibitem[Wang et~al.(2018)Wang, Zhu, and Ma]{WangZhuMa2018}
HaiYing Wang, Rong Zhu, and Ping Ma.
\newblock Optimal subsampling for large sample logistic regression.
\newblock \emph{Journal of the American Statistical Association}, 113\penalty0
  (522):\penalty0 829--844, 2018.

\bibitem[Wang et~al.(2019)Wang, Yang, and Stufken]{WangYangStufken2019}
HaiYing Wang, Min Yang, and John Stufken.
\newblock Information-based optimal subdata selection for big data linear
  regression.
\newblock \emph{Journal of the American Statistical Association}, 114\penalty0
  (525):\penalty0 393--405, 2019.
\newblock \doi{10.1080/01621459.2017.1408468}.

\bibitem[Wang et~al.(2021)Wang, Zhang, and Wang]{WangZhangWang2021}
HaiYing Wang, Aonan Zhang, and Chong Wang.
\newblock Nonuniform negative sampling and log odds correction with rare events
  data.
\newblock In \emph{Proceedings of The 35 Conference on Neural Information
  Processing Systems ({NeurIPS} 2021).}, Proceedings of Machine Learning
  Research. PMLR, 2021.

\bibitem[Wang et~al.(2022)Wang, Zou, and Wang]{WangZou2019}
Jing Wang, Jiahui Zou, and HaiYing Wang.
\newblock Sampling with replacement vs poisson sampling: A comparative study in
  optimal subsampling.
\newblock \emph{IEEE Transactions on Information Theory}, 68\penalty0
  (10):\penalty0 6605--6630, 2022.
\newblock \doi{10.1109/TIT.2022.3176955}.

\bibitem[Whittemore(1997)]{whi97}
A.S. Whittemore.
\newblock Multistage sampling designs and estimating equations.
\newblock \emph{Journal of the Royal Statistical Society. Series B},
  59:\penalty0 589--602, 1997.

\bibitem[Wild(1991)]{wil91}
C.J. Wild.
\newblock Fitting prospective regression models to case-control data.
\newblock \emph{Biometrika}, 78:\penalty0 705--717, 1991.

\bibitem[Yang et~al.(2015)Yang, Zhang, Jin, and Zhu]{yang2015explicit}
Tianbao Yang, Lijun Zhang, Rong Jin, and Shenghuo Zhu.
\newblock An explicit sampling dependent spectral error bound for column subset
  selection.
\newblock In \emph{International Conference on Machine Learning}, pages
  135--143, 2015.

\bibitem[Yao and Wang(2019)]{yao2019optimal}
Yaqiong Yao and HaiYing Wang.
\newblock Optimal subsampling for softmax regression.
\newblock \emph{Statistical Papers}, 60\penalty0 (2):\penalty0 235--249, 2019.

\bibitem[Yao and Wang(2021)]{YaoWang2021JDS}
Yaqiong Yao and HaiYing Wang.
\newblock A review on optimal subsampling methods for massive datasets.
\newblock \emph{Journal of Data Science}, 19\penalty0 (1):\penalty0 151–172,
  2021.

\bibitem[Yu et~al.(2022)Yu, Wang, Ai, and Zhang]{yu2020quasi}
Jun Yu, HaiYing Wang, Mingyao Ai, and Huiming Zhang.
\newblock Optimal distributed subsampling for maximum quasi-likelihood
  estimators with massive data.
\newblock \emph{Journal of the American Statistical Association}, 117\penalty0
  (537):\penalty0 265--276, 2022.
\newblock \doi{10.1080/01621459.2020.1773832}.
\newblock URL \url{https://doi.org/10.1080/01621459.2020.1773832}.

\end{thebibliography}

\end{document}